\documentclass[final]{siamltex}
\usepackage{graphicx, amssymb, amsmath}
\usepackage{color}
\usepackage{soul} 
\usepackage{booktabs}
\usepackage{boxedminipage}
\usepackage{amsfonts}
\usepackage{algorithmic}
\usepackage{algorithm,booktabs}
\usepackage[mathscr]{eucal}
\newtheorem{dfn}{Definition}

\newcommand{\R}{\mathbb{R}}

\begin{document}

\title{Nonlinear model order reduction via\\ Dynamic Mode Decomposition \thanks{J. N. Kutz acknowledges support from the Air Force Office of Scientific Research (FA9550-15-1-0385). J. N. Kutz would like to thank Steven Brunton, Joshua Proctor, Bing Brunton, Matthew Williams, Jonathan Tu and Clancy Rowley for invaluable discussions related to the dynamic mode decomposition and Koopman operator theory.}}

\author{
Alessandro Alla\footnote{University of Hamburg, Department of Mathematics, 
 Bundesstras{\ss}e, 55, 20146, Hamburg, Germany
\texttt{alessandro.alla@uni-hamburg.de}}, 
J. Nathan Kutz \footnote{University of Washington, Department of Mathematics, Seattle, WA 98195, United States
\texttt{kutz@uw.edu}}
%
%
}

\maketitle

\begin{abstract}
We propose a new technique for obtaining reduced order models for nonlinear dynamical systems.
Specifically, we advocate the use of the recently developed Dynamic Mode Decomposition (DMD), an equation-free method, to approximate the nonlinear term. DMD is a spatio-temporal matrix decomposition of a data matrix that
correlates spatial features while simultaneously associating the activity with periodic temporal behavior.
With this decomposition, one can obtain a fully reduced dimensional surrogate model and avoid the evaluation of the nonlinear term in the online stage. This allows for an impressive speed up of the computational cost, and, at the same time, accurate approximations of the problem.  We present a suite of numerical tests to illustrate our approach and to show the effectiveness of the method in comparison to existing approaches. 
\end{abstract}



\maketitle

\section{Introduction}
\label{Section1}
\setcounter{section}{1}
\setcounter{equation}{0}
\setcounter{theorem}{0}
\setcounter{algorithm}{0}
\renewcommand{\theequation}{\arabic{section}.\arabic{equation}}

Reduced-order models (ROMs) are of growing importance in scientific computing as they
provide a principled approach to approximating high-dimensional PDEs with low-dimensional models.
Indeed, the dimensionality reduction provided by ROMs help reduce the computational complexity and time needed to solve large-scale, engineering systems~\cite{rom_book, karen1}, enabling simulation based scientific studies not possible even a decade ago. One of the primary challenges in producing the low-rank dynamical system is efficiently projecting the nonlinearity of the governing PDEs (inner products)~\cite{BMNP04, CS10} on to the proper orthogonal decomposition (POD)~\cite{Lumley:1970, HLBR_turb,Vol11} basis. This fact was recognized early on in the ROM community, and methods such as gappy POD~\cite{gap1, gap2, karni} where proposed to more efficiently enable the task. More recently, the empirical interpolation method (EIM)~\cite{BMNP04}, and the simplified discrete empirical interpolation method (DEIM)~\cite{CS10} for the proper orthogonal decomposition (POD)~\cite{Lumley:1970, HLBR_turb,Vol11}, have provided a computationally efficient method for discretely (sparsely) sampling and evaluating the nonlinearity. These broadly used and highly-successful methods ensure that the computational complexity of ROMs scale favorably with the rank of the approximation, even for complex nonlinearities. As an alternative to the EIM/DEIM architecture, we propose using the recently developed Dynamic Mode Decomposition (DMD) for producing low-rank approximations of the PDE nonlinearities. DMD provides a decomposition of data into spatio-temporal modes that correlates the data across spatial features (like POD), but also associates the correlated data to unique temporal Fourier modes, allowing for a computationally efficient regression of the nonlinear terms to a least-square fit linear dynamics approximation. We demonstrate that the POD-DMD method produces a viable ROM architecture, scaling favorable in computational efficiency relative to POD-DEIM.

ROMs are fundamentally based upon dimensionality-reduction techniques, one of the most common of which is the 
proper orthogonal decomposition~\cite{Lumley:1970, HLBR_turb} for selecting an optimal low-dimensional basis for projecting the PDE dynamics. Indeed, the POD architecture is ubiquitous
across many engineering disciplines. Some of its variants are alternatively referred to as principal components analysis (PCA,~\cite{Pearson:1901}), the Karhunen--Lo\`eve (KL) decomposition, empirical orthogonal functions (EOF,~\cite{lorenzMITTR56}), and/or the Hotelling transform~\cite{hotellingJEdPsy33}. The success of the POD method is based upon the observation that many
nonlinear dynamical systems (PDEs) often exhibit low-dimensional phenomena, so that the majority of variance/energy is optimally contained in a small number of modes computed from a singular value decomposition (SVD). One can select a POD basis by a pre-determined cut-off value, such as when the modal basis contain $99\%$ of the variance, so that only the first $\ell$-modes ($\ell$-rank truncation) are kept, or by a more principled truncation such as the optimal hard-threshold value for systems with well-characterized white noise~\cite{gavish}. The truncated POD modes are then used as the basis modes for a Galerkin expansion (Galerkin-POD) for projecting the dynamics onto a rank-$\ell$ dynamical system~\cite{HLBR_turb, Kutz:2013}. 

The POD-Galerkin method has been widely used in the scientific computing community. The primary challenge in producing the low-rank dynamical system is efficiently projecting the nonlinearity (inner products) to the POD basis, leading to numerous innovations in the ROM community for interpolating the projection. Starting with the gappy POD technique~\cite{gap1, gap2, karni}, sparse sampling was recognized early on as an effective method for approximating the nonlinear inner products. The EIM and DEIM methods proposed an algorithm for improving the greedy selection of discrete spatial points for producing an interpolated approximation of the nonlinear terms. 
This ensures that the computational cost of evaluating the nonlinearity remains proportional to the rank of the reduced POD basis. The DEIM approach combines projection with interpolation by selecting interpolation indices  to specify  an interpolation-based projection for a nearly optimal $\ell_2$ subspace approximating the nonlinearity. The EIM/DEIM are not the only methods developed to reduce the complexity of evaluating nonlinear terms, see for instance the missing point estimation (MPE,~\cite{mpe}), {{{{\em best points} method~\cite{patera}, or the so-called GNAT gappy POD~\cite{Carlberg:2013} method}}. However, they have been successful in a large number of diverse applications and models~\cite{CS10}. In any case, the MPE, gappy POD, and EIM/DEIM all use a small selected set of spatial grid points to avoid evaluation of the expensive inner products required to evaluate nonlinear terms 

An alternative to these sparse sampling techniques for evaluating the nonlinear inner products
is the DMD method. DMD provides a decomposition of data into spatio-temporal modes that correlates the data across spatial features (like POD), but also associates the correlated data to unique temporal Fourier modes. 
More precisely, DMD computes a regression of the sampled data to a best fit (least-squares) linear, constant-coefficient, system of differential equations. This spatio-temporal regression allows us to directly project the nonlinear terms in the PDE to its best-fit time dynamics. Like EIM/DEIM, it requires a singular value decomposition to generate the approximation. We demonstrate that the DMD provides a highly efficient approximation for ROMs, 
performing a much more rapid evaluation of the nonlinear terms in comparison to the EIM/DEIM methods. At its core, the DMD method can be thought of as an ideal combination of spatial dimensionality-reduction techniques, such as POD, with Fourier Transforms in time. It also allows for further innovations that integrate the DMD with key concepts from multi-resolution analysis~\cite{Kutz:mrdmd} and sparsity/compression~\cite{brunton:cs}, allowing one to potentially generalize
the proposed method to multi-scale physics problems at greatly improved speeds.

The structure of the paper is as follows. In Section \ref{Section2} we recall the Proper Orthogonal Decomposition method and the Discrete Empirical Interpolation method applied to a general dynamical system. Section \ref{Section3} explaines the Dynamic Mode Decomposition and compares the DMD method used as equation-free or as a Galerkin projection method. The coupling between POD and DMD is explained in Section \ref{Section4}. Finally, numerical tests are presented in Section \ref{Section5}. Throughout the paper we use the following notation:  all matrices and vectors are in bold letters. The basis functions are denoted by the matrix $\bf {\bf \Psi}$ with different superscripts denoting how we computed the basis, e.g. ${\bf {\bf \Psi}}^{\mbox{\tiny POD}}$ represents the basis functions from the POD method. The rank of the POD basis functions is $\ell$, whereas the rank of the nonlinear term is $k$.


\section{Problem Formulation}
\label{Section2}
\setcounter{section}{2}
\setcounter{equation}{0}
\setcounter{theorem}{0}
\setcounter{algorithm}{0}
\renewcommand{\theequation}{\arabic{section}.\arabic{equation}}

In what follows, we consider a system of ordinary differential equations:
\begin{equation}
\left\{ \begin{array}{l}\label{ode}
{\bf M}\dot{{\bf y}}(t)={\bf A}{\bf y}(t)+{\bf f}(t,{\bf y}(t)),\;\; t\in(0,T]\\
{\bf y}(0)={\bf y_0},\\
\end{array} \right.
\end{equation}
where ${\bf y_0}\in\R^n$ is a given initial data, ${\bf M, A}\in \R^{n\times n}$ given matrices and ${\bf f}:[0,T]\times\R^n\rightarrow\R^n$ a continuous function in both arguments and locally Lipschitz-type with respect to the second variable. It is well--known that under these assumptions there exists an unique solution for \eqref{ode}.
 
 This wide class of problems arises in many applications, especially from the numerical approximation of partial differential equations. In such cases, the dimension of the problem $n$ is the number of spatial grid points used from discretization and it can be very large. The solution of system \eqref{ode} may be very expensive and therefore it might be useful to simplify the complexity of the problem by means of reduced order modeling techniques.

\subsection{The POD method and reduced-order modeling}
\label{Section2.1}
One popular method for reducing the complexity of the system is the so-called Proper Orthogonal Decomposition (POD). The idea was proposed by Sirovich~\cite{Sir87} and is detailed here for completeness. We build an equidistant grid in time with constant step size $\Delta t$. Let $t_0:=0<t_1<t_2<\ldots<t_m\leq T$ with $t_j=j\Delta t,\;\;j=0,\ldots,m$. Let us assume we know the exact solution of (\ref{ode}) on the time grid points $t_j$, $j\in \{1,\ldots,m\}$. Our aim is to determine a POD basis of rank $\ell\ll n$ to describe the set of data collected in time by solving the following minimization problem:
\begin{equation}\label{pbmin}
\min_{ {\boldsymbol{\psi}}_1,\ldots,{\boldsymbol{\psi}}_\ell\in\R^n} \sum_{j=1}^m \alpha_j\left\|{\bf y}(t_j)-\sum_{i=1}^\ell \langle {\bf y}(t_j),{\boldsymbol{\psi}}_i\rangle{\boldsymbol{\psi}}_i\right\|_W^2\quad \mbox{such that }\langle {\boldsymbol{\psi}}_i,{\boldsymbol{\psi}}_j\rangle=\delta_{ij},
\end{equation}
where the coefficients $\alpha_j$ are non-negative and ${\bf y}(t_j)$ are the so called {\em snapshots}, e.g. the solution of \eqref{ode} at a given time $t_j$.  Additionally, we assume ${\bf y}(t_j)\in V$ for a suitable Hilbert space $V$. The norm, here and in the sequel of the section, can be interpreted as the weighted norm such that $\langle {\bf u}, {\bf v}\rangle_{\tiny {\bf W}}={\bf u}^T{\bf W}{\bf v}$ and $\|\cdot\|^2=\langle\cdot,\cdot\rangle_{\tiny {\bf W}}$
where ${\bf W}\in\R^{n\times n}$ is a positive definite weighting matrix.

Solving (\ref{pbmin}) we look for  an orthonormal basis $\{{\boldsymbol{\psi}}_i\}_{i=1}^\ell$ which minimizes the distance between the sequence ${\bf y}(t_j)$ with respect to its projection onto this unknown basis. The matrix ${\bf Y}$ contains the collection of snapshots ${\bf y}(t_j)$ as columns. It is useful to look for $\ell\ll\min\{m,n\}$ in order to reduce the dimension of the problem considered.
The solution of (\ref{pbmin}) is given by the singular value decomposition of the snapshots matrix ${\bf W}^{1/2}{\bf Y}={\bf{\bf \Psi}}{\bf\Sigma} {\bf V}^T$, where we consider the first $\ell-$ columns $\{{\boldsymbol{\psi}}_i\}_{i=1}^\ell,$ of the orthogonal matrix ${\bf {\bf \Psi}}$ and set  ${\bf {\bf \Psi}}^{\mbox{\tiny POD}}={\bf W}^{-1/2} {\bf \Psi}$.

To concretely apply the POD method, the choice of the truncation parameter $\ell$ plays a crucial role. There are no a-priori estimates which guarantee the ability to build a coherent reduced model, but one can focus on heuristic considerations, introduced by Sirovich \cite{Sir87}, so as to have the following ratio close to one:
\begin{equation}\label{indicator}
\mathcal{E}(\ell)=\dfrac{\sum\limits_{i=1}^\ell\sigma^2_i}{\sum\limits_{i=1}^d\sigma^2_i}.
\end{equation}
This indicator is motivated by the fact that the error in \eqref{pbmin} is given by the singular values we neglect:
\begin{equation}\label{err_pod}
\sum_{j=1}^m \alpha_j\left\|{\bf y}(t_j)-\sum_{i=1}^\ell \langle {\bf y}(t_j),{\boldsymbol{\psi}}_i\rangle{\boldsymbol{\psi}}_i\right\|_{\bf W}^2=\sum_{i=\ell+1}^d \sigma_i^2,
\end{equation}
where $d$ is the rank of the snapshot matrix ${\bf Y}$. We note that the error \eqref{err_pod} is strictly related to the computation of the snapshots and it is not related to the reduced dynamical system.

Let us assume that we have computed the POD basis functions ${\bf {\bf \Psi}}^{\mbox{\tiny POD}}=\{{\boldsymbol{\psi}}_j\}_{j=1}^\ell\in\R^{n\times \ell}$ of rank $\ell$ for the problem \eqref{ode}, we make the following projection of the dynamics:
\begin{equation}\label{pod_ans}
{\bf y}(t)\approx{\bf {\bf \Psi}}^{\mbox{\tiny POD}} {\bf y^\ell}(t),
\end{equation}
where ${\bf y}^\ell(t)$ are functions from $[0,T]$ to $\R^\ell.$ We note that we are working with a Galerkin-type projection where we consider only few basis functions whose support is non-local, unlike Finite Element basis functions. The reduced solution ${\bf y}^\ell(t)\in V^\ell\subset V$ where $V^\ell=\mbox{span}\{{\boldsymbol{\psi}}_1,\ldots,{\boldsymbol{\psi}}_\ell\}$.

Inserting the projection assumption \eqref{pod_ans} into the full model \eqref{ode}, and making use of the orthogonality of the POD basis functions, the reduced model takes the following form: 
\begin{equation}\label{pod_sys}
\left\{\begin{array}{l}
{\bf M}^\ell\dot{{\bf y}}^\ell(t)={\bf A}^\ell {\bf y}^\ell(t)+{\bf {\bf \Psi}}^Tf(t,{\bf {\bf \Psi}} {\bf y}^\ell(t))\\
{\bf y}^\ell(0)={\bf y_0^\ell}
\end{array}\right.
\end{equation}
where $({\bf M}^\ell)_{ij}=\langle {\bf M}{\boldsymbol{\psi}}_i,{\boldsymbol{\psi}}_j\rangle, ({\bf A}^\ell)_{ij}=\langle {\bf A}{\boldsymbol{\psi}}_i,{\boldsymbol{\psi}}_j\rangle \in\R^{\ell\times\ell}$ and ${\bf y_0^\ell}=({\bf {\bf \Psi}}^{\mbox{\tiny POD}})^T{\bf y}_0\in\R^\ell$. We also note that ${\bf M}^\ell, {\bf A}^\ell\in\R^{\ell\times\ell}$.
The system (\ref{pod_sys}) is achieved following a Galerkin projection where the basis functions are computed by the POD method given by \eqref{ode}. If the dimension of the system is $\ell\ll n$, then a significant dimensionality reduction is accomplished. We note that an error analysis for $\|y(t)-{\bf {\bf \Psi}}^{\mbox{\tiny POD}}y^\ell(t)\|$ can be found in \cite{KV02}.

\subsection{Discrete Empirical Interpolation Method}
\label{Section2.2}


For the results in this section we closely follow the presentation in \cite{Vol11}. The ROM introduced in (\ref{pod_sys}) is a nonlinear system where the significant challenge with the POD--Galerkin approach is the computational complexity associated with the evaluation of the nonlinearity. To illustrate this issue, we consider the nonlinearity in (\ref{pod_sys}): 
$$ {\bf F}(t,{\bf y}^\ell(t))=({\bf {\bf \Psi}}^{\mbox{\tiny POD}})^T {\bf f}(t,{\bf {\bf \Psi}}^{\mbox{\tiny POD}} {\bf y}^\ell(t))=\langle {\bf f}(t,{	\bf y}(t)), {\bf {\bf \Psi}}^{\mbox{\tiny POD}}\rangle.$$
To compute this inner product, the variable ${\bf y}^\ell(t)\in\R^\ell$ is first expanded to an $n-$dimensional vector ${\bf {\bf \Psi}}^{\mbox{\tiny POD}} {\bf y}^\ell(t)\in\R^n$, then the nonlinearity ${\bf f}(t,{\bf {\bf \Psi}}^{\mbox{\tiny POD}} {\bf y}^\ell(t))$ is evaluated and, at the end, we return back to the reduced-order model. This is computationally expensive since it implies that the evaluation of the nonlinear term requires computing the full, high-dimensional model, and therefore the reduced model is not independent of the full dimension $n.$ We note that, for simplicity, we dropped the weighted inner product.

To avoid this computationally expensive, high-dimensional, evaluation the {\em Empirical Interpolation Method} (EIM, \cite{BMNP04}) and {\em Discrete Empirical Interpolation Method} (DEIM, \cite{CS10}) were introduced. We note that DEIM is built upon EIM:  the two methods are essentially equivalent and are based on a POD approach combined with a greedy algorithm. DEIM is the tensorial matricial version of EIM, and it is used here due to the nature of our time-dependent problem. The interested reader is referred to \cite{CS10} for further information. 

The computation of the POD basis functions for the nonlinear part are related to the set of the snapshots ${\bf f}(t_j,{\bf y}(t_j))$ where ${\bf y}(t_j)$ is already computed from \eqref{ode}. We denote with ${\bf U}\in\R^{n\times k}$ the POD basis function of rank $k$ of the nonlinear part. The DEIM approximation of ${\bf  f}(t,{\bf y}(t))$ is as follows
$${\bf f}^{\mbox{\tiny DEIM}}(t,{\bf y}^{\mbox{\tiny DEIM}}(t))={\bf U}({\bf S}^T {\bf U})^{-1} {\bf f}(t,{\bf y}^{\mbox{\tiny DEIM}}(t))$$
where ${\bf S}\in\R^{n\times k}$ and ${\bf y}^{\mbox{\tiny DEIM}}(t)={\bf S}^T{\bf {\bf \Psi}}^{\mbox{\tiny POD}}{\bf y}^\ell(t)$. The matrix ${\bf S}$ is the interpolation point where the nonlinearity is evaluated and the selection of its points is made according to an LU decomposition algorithm with pivoting~\cite{CS10}, or following the QR decomposition with pivoting~\cite{DG15}. We note that, here, we suppose that the maths ${\bf S}$ can be moved into the nonlinearity. 

Let us define ${\bf {\bf \Psi}}^{\mbox{\tiny DEIM}}:={\bf U}({\bf S}^T {\bf U})^{-1}$. We note that ${\bf {\bf \Psi}}^{\mbox{\tiny DEIM}}\in\R^{n\times k}$. Therefore the reduced nonlinearity may be expressed as:
$$({\bf {\bf \Psi}}^{\mbox{\tiny POD}})^T{\bf f}^{\mbox{\tiny DEIM}}(t,{\bf y}^{\mbox{\tiny DEIM}}(t))=({\bf {\bf \Psi}}^{\mbox{\tiny POD}})^T{\bf {\bf \Psi}}^{\mbox{\tiny DEIM}} {\bf f}(t,{\bf y}^{\mbox{\tiny DEIM}})$$
 where we only select a small (sparse) number of rows of ${\bf {\bf \Psi}}^{\mbox{\tiny POD}}{\bf y}^\ell(t)$.
As for the computational expense, the matrices $${\bf S}^T{\bf \Psi}^{\mbox{\tiny POD}}\in\R^{k\times\ell},\,({\bf S}^T {\bf U})^{-1}\in\R^{k\times k} \hbox{ and } ({\bf \Psi}^{\mbox{\tiny POD}})^T{\bf \Psi}^{\mbox{\tiny DEIM}}\in\R^{\ell\times k}$$ can all be precomputed. All the precomputed quantities are independent of the full dimension $n.$ Additionally, during the iteration process the nonlinearity needs only to be evaluated at the $k$ interpolation points since ${\bf S}^T{\bf \Psi} {\bf y}^\ell(t)\in\R^k.$ Typically the dimension $k$ is much smaller than the full dimension. This allows the reduced-order model to be completely independent of the full dimension as follows:
\begin{equation}\label{pod_sysdeim}
\left\{\begin{array}{l}
{\bf M}^\ell\dot{{\bf y}}^\ell(t)={\bf A}^\ell {\bf y}^\ell(t)+({\bf \Psi}^{\mbox{\tiny POD}})^T{\bf \Psi}^{\mbox{\tiny DEIM}} {\bf f}(t,{\bf y}^{\mbox{\tiny DEIM}})\\
{\bf y}^\ell(0)={\bf y_0}^\ell.
\end{array}\right.
\end{equation}
We note that the only difference with respect to \eqref{pod_sys} is the low-rank approximation of the nonlinear term. The error between ${\bf f}(t,{\bf y}(t))$ and its DEIM approximation $f^{\mbox{\tiny DEIM}}$ is given by
$$\|{\bf f}-{\bf f}^{\mbox{\tiny DEIM}}\|_2\leq  c\|({\bf I}-{\bf UU}^{T})f\|_2\quad \,\,\, \mbox{with} \,\,\, c=\|({\bf S}^T{\bf U})^{-1}\|_2$$
where different error performance is achieved depending on the selection of the interpolation points in $S$ as shown in \cite{DG15}.

\section{Dynamic Mode Decomposition}
\label{Section3}
\setcounter{section}{3}
\setcounter{equation}{0}
\setcounter{theorem}{0}
\setcounter{algorithm}{0}
\renewcommand{\theequation}{\arabic{section}.\arabic{equation}}

DMD is an {\em equation-free}, data-driven method capable of providing accurate assessments of the spatio-temporal coherent structures in a given complex system, or short-time future estimates of such a systems. It traces its origins to pioneering work of Bernard Koopman in 1931~\cite{koopman}, whose work was revived by Igor Mezi\'c and co-workers starting in 2004~\cite{Mezic2004, Mezic2005, mezic2}. Koopman theory is a dynamical systems tool that provides information about a nonlinear dynamical system via an associated infinite-dimensional linear system. Specifically, it provides a characterization that is readily interpretable in terms of standard methods of dynamical systems. Defining it as a data-driven algorithm, Schmid~\cite{DMD0, DMD1} proposed the DMD architecture for modeling complex flows, with Rowley et al.~\cite{DMD4} showing quickly thereafter that the DMD method is actually a special case of Koopman theory.

Given the connection between DMD and Koopman theory~\cite{Mezic2004, Mezic2005, DMD4}, we begin by defining the
Koopman operator:\\

\begin{dfn}[Koopman Operator~\cite{koopman}]
For a dynamical system:
\begin{equation}  \label{eq:Ukoop}
  \frac{d{\bf y}}{dt} = \bf{N}({\bf y}),
\end{equation}
where ${\bf y}\in {\mathcal{M}}$, an $n$-dimensional manifold. 
The Koopman operator $\mathcal{ K}$ acts
on a set of scalar observable functions ${g}:  \mathcal{ M} \rightarrow \mathbb{C}$  so that 
\begin{equation}
  \mathcal{ K} {g}({\bf y})  = {g}\left( {\bf N}({\bf y}) \right) \, . \label{eq:koop}
\end{equation}
\end{dfn}


This shows that the Koopman operator is a {\em linear} operator that acts on scalar functions $g$. In a general setting, the Koopman operator can act on a set of observables $g_j$ that are denoted by components of the vector ${\bf g}$. But as already mentioned, the DMD is a specific realization of the Koopman theory. Specifically, the observable is the state space itself so that one considers the linear observable ${\bf g}({\bf y})={\bf y}$.

When considering such a linear observable, the DMD algorithm determines the Koopman eigenvalues and modes directly from data. Specifically, one can use the recent formal definition of the DMD method~\cite{DMD5}:\\

\begin{dfn}[Dynamic Mode Decomposition \cite{DMD5}] 
Suppose we have a dynamical system (\ref{eq:Ukoop}) and two sets of data 
%
\begin{equation}
  {\bf Y} \!=\! \begin{bmatrix}
\vline & \vline & & \vline \\
{\bf y}(t_0) & {\bf y}(t_1) & \cdots & {\bf y}(t_{m-1})\\
\vline & \vline & & \vline
\end{bmatrix}, \hspace{0.1in} {\bf Y}' \!=\! \begin{bmatrix}
\vline & \vline & & \vline \\
{\bf y}(t_1) & {\bf y}(t_2) & \cdots & {\bf y}(t_m)\\
\vline & \vline & & \vline
\end{bmatrix}
\end{equation}
%
with ${\bf y}(t_j)$ an initial condition to (\ref{eq:Ukoop}) and ${\bf y}(t_{j+1})$ it corresponding output
after some prescribed evolution time $\tau$ with there being $m$ initial conditions considered. The DMD modes
are eigenvectors of 
\begin{equation}
  {\bf A}_{\bf y} = {\bf Y}' {\bf Y}^\dag
  \label{eq:newDMD}
\end{equation}
where $\dag$ denotes the Moore-Penrose pseudoinverse.
\end{dfn}\\

%
\noindent
The definition of DMD thus yields the matrix ${\bf A}_{\bf y}$, which is a finite dimensional approximation of
the Koopman operator for a linear observable.

The definition of DMD produces a regression procedure whereby the
data snapshots in time are used to produce the best-fit linear dynamical system for the data ${\bf Y}$.
The DMD procedure thus constructs the proxy, approximate linear evolution
\begin{equation}
\frac{d \tilde{\bf y}}{dt} = {\bf A}_{\bf y} \tilde{\bf y} \label{eq:dA}
\end{equation}
with $\tilde{\bf y}(0) = \tilde{\bf y}_0$ and whose solution is
\begin{equation}
\tilde{\bf y}(t)=  \sum_{i=1}^k b_k {\boldsymbol{\psi}}_i \exp(\omega_i t) \, ,
\label{eq:omegaj}
\end{equation}
where ${\boldsymbol{\psi}}_i$ and $\omega_i$ are the eigenfunctions  and eigenvalues 
of the matrix ${\bf A}_{\bf y}$.
The ultimate goal in the DMD algorithm is to optimally construct the matrix ${\bf A}_{\bf y}$ so
that the true and approximate solution remain optimally close in a least-square sense, i.e.  
$\| {\bf y}(t) - \tilde{\bf y}(t) \| \ll 1$. 
Of course, the optimality of the approximation holds only over the sampling window where ${\bf A}_{\bf y}$ is constructed, but the approximate solution
can be used to not only make future state predictions, but also to decompose the dynamics 
into various time-scales since the $\omega_k$ are prescribed. Moreover, the DMD typically makes use
of low-rank structure so that the total number of modes, $k\ll n$, allows for dimensionality reduction of the
dynamical system.

In effect, the least-square regression of the nonlinear dynamical system to the linear system (\ref{eq:dA})
allows us to approximate the governing equation (\ref{ode}) in the following manner:
\begin{equation}
 \begin{array}{l}\label{ode2}
{\bf M}\dot{{\bf y}}(t)={\bf A}{\bf y}(t)+{\bf f}(t,{\bf y}(t)) \approx {\bf A}{\bf y}(t) + {\bf A}_{\bf y} {\bf y}(t)\\
\end{array} 
\end{equation}
where the DMD algorithm constructs the matrix ${\bf A}_{\bf y}$ approximating the nonlinearity over
the snapshots collected.

In practice, the matrix ${\bf A}_{\bf y}$ is, in general, highly ill-conditioned and when the state dimension $n$ is large, the aforementioned matrix may be even intractable to analyze directly. Instead, DMD circumvents the eigendecomposition of ${\bf A}_{\bf y}$ by considering a rank-reduced representation in terms of a POD-projected matrix $\tilde{\bf A}_{\bf y}$. The DMD algorithm proceeds as follows~\cite{DMD5}:\\

\begin{enumerate}
\item First, take the SVD of ${\bf Y}$:
\begin{equation}
  {\bf Y} = {\bf U} {\bf \Sigma} {\bf V}^*,
  \label{eq:kry3}
\end{equation}
where $*$ denotes the conjugate transpose, 
${\bf U}\in {\mathbb{C}}^{n\times k}$, ${\bf \Sigma}\in {\mathbb{C}}^{k\times k}$
and ${\bf V}\in {\mathbb{C}}^{m-1\times k}$. Here $k$ is the rank of the reduced SVD
approximation to ${\bf Y}$. The left singular vectors ${\bf U}$ are POD modes.

The SVD reduction in (\ref{eq:kry3}) could also be exploited at this stage in the algorithm 
to perform a low-rank truncation of the data. Specifically, if low-dimensional structure is
present in the data, the singular values of ${\bf \Sigma}$ will decrease sharply to zero with
perhaps only a limited number of dominant modes. A principled way to truncate noisy data would be to use the
recent hard-thresholding algorithm of Gavish and Donoho~\cite{gavish}.

\item Next, compute $\tilde{\bf A}_{\bf y}$, the $k\times k$ projection of the full matrix ${\bf A}_{\bf y}$ onto POD modes:
\begin{eqnarray}
{\bf A}_{\bf y} &=& {\bf Y}' {\bf V}\boldsymbol{\Sigma}^{-1}{\bf U}^*\nonumber\\
\Longrightarrow\quad\tilde{\bf A}_{\bf y} &=& {\bf U}^*{\bf A}_{\bf y}{\bf U} = {\bf U}^*{\bf Y}' {\bf V}\boldsymbol{\Sigma}^{-1}.
\end{eqnarray}
\item Compute the eigendecomposition of $\tilde{\bf A}_{\bf y}$:
\begin{equation}
\tilde{\bf A}_{\bf y}{\bf W} = {\bf W}\boldsymbol{\Lambda},
\end{equation}
where columns of ${\bf W}$ are eigenvectors and $\boldsymbol{\Lambda}$ is a diagonal matrix containing the corresponding eigenvalues $\lambda_i$.
\item Finally, we may reconstruct eigendecomposition of ${\bf A}_{\bf y}$ from ${\bf W}$ and $\boldsymbol{\Lambda}$. In particular, the eigenvalues of ${\bf A}_{\bf y}$ are given by $\boldsymbol{\Lambda}$ and the eigenvectors of ${\bf A}_{\bf y}$ (DMD modes) are given by columns of $\bf {\bf \Psi}$:
\begin{equation}
{\bf {\bf \Psi}} = {\bf Y}' {\bf V}\boldsymbol{\Sigma}^{-1}{\bf W}.
\label{eq:DMDjtu}
\end{equation}
\end{enumerate}
Note that Eq.~\eqref{eq:DMDjtu} from~\cite{DMD5} differs from the formula ${\bf {\bf \Psi}}={\bf U}{\bf W}$ from~\cite{DMD1}, although these will tend to converge if ${\bf Y}$ and ${\bf Y}'$ have the same column spaces. As a peuso-algorithm, it can be summarized in 
Algorithm \ref{Alg_DMD}


\begin{algorithm}
\caption{Exact DMD}
\label{Alg_DMD}
\begin{algorithmic}[1]
\REQUIRE Snapshots $\{{\bf y}(t_0),\ldots,{\bf y}(t_m)\}$,
\STATE Set ${\bf Y}=[{\bf y}(t_0),\dots, {\bf y}(t_{m-1})]$ and $Y'=[{\bf y}(t_1),\dots, {\bf y}(t_m)]$,
\STATE Compute the SVD of ${\bf Y}$, ${\bf Y}={\bf U}{\bf \Sigma}{\bf  V}^T$
\STATE Define $\tilde{{\bf A}}_{\bf y}:={\bf U}^*{\bf Y}'{\bf V}{\bf \Sigma}^{-1}$
\STATE Compute eigenvalues and eigenvectors of $\tilde{{\bf A}}_{\bf y} {\bf W}={\bf W}{\bf \Lambda}$.
\STATE Set ${\bf \Psi}^{\mbox{\tiny DMD}}={\bf Y}'{\bf V}{\bf \Sigma}^{-1}{\bf W}$ \\
\end{algorithmic}
\end{algorithm}

\subsection{Applications of the DMD method}
In this section, we propose two different applications of the DMD method to ROMs. Our first application concerns the interpolation of a parametrized function which is compared with the DEIM approach. The second one is related to the reduction of dynamical systems and considers the DMD method as a Galerkin projection strategy.

\paragraph{\bf Test 1: Interpolation of parametrized functions} 
Let us consider the following  nonlinear parametrized functions:
\begin{equation}
s(x;\mu)=(1-x)\cos(3\pi\mu(x+1))e^{-(1+x)\mu}
\end{equation}
where $s:\Omega\times\mathcal{D}\rightarrow\R$, $x\in\Omega=[-1,1]$ and $\mu\in\mathcal{D}=[1,\pi]$. This nonlinear function is from \cite{CS10}.  Let us discretize the space domain $[x_1,\ldots,x_n]\in\R^n$ with $x_i$ equidistant in $\Omega$. With compact notation we define $\bf{f}:\mathcal{D}\rightarrow\R^n$ by
\begin{equation}
{\bf f}(\mu)=[s(x_1;\mu),\ldots,s(x_n;\mu)]\in\R^n
\end{equation}
for $\mu\in\mathcal{D}$. This example uses $51$ snapshots of ${\bf f}(\mu_j)$ to build the POD basis functions where $\mu_i$ are equidistributed point in $[1,\pi]$ and $n=101$. 
Figure \ref{par_1} shows the behavior of the function with $\mu=\{1.17, 3.1\}$ and the decay of the singular values of the snapshots set.

The purpose of this subsection is to show that the DMD might also be used as interpolation method as it is shown in Figure \ref{par_3}. As we can see DMD is able to reconstruct the parametrized functions in $\mu=\{1.17, 3.1\}$ which is not included in the snapshot set.
If we look more closely into this approximations and compare it with the DEIM interpolation method we can see that the DMD method is always really faster then DEIM (Figure \ref{par_2} left) and the error, at the beginning, is comparable in fact up to the first 10 modes we have same error. The error is computed with respect to the Frobenius norm.

\begin{figure}[t]
\centering
\includegraphics[scale=0.3]{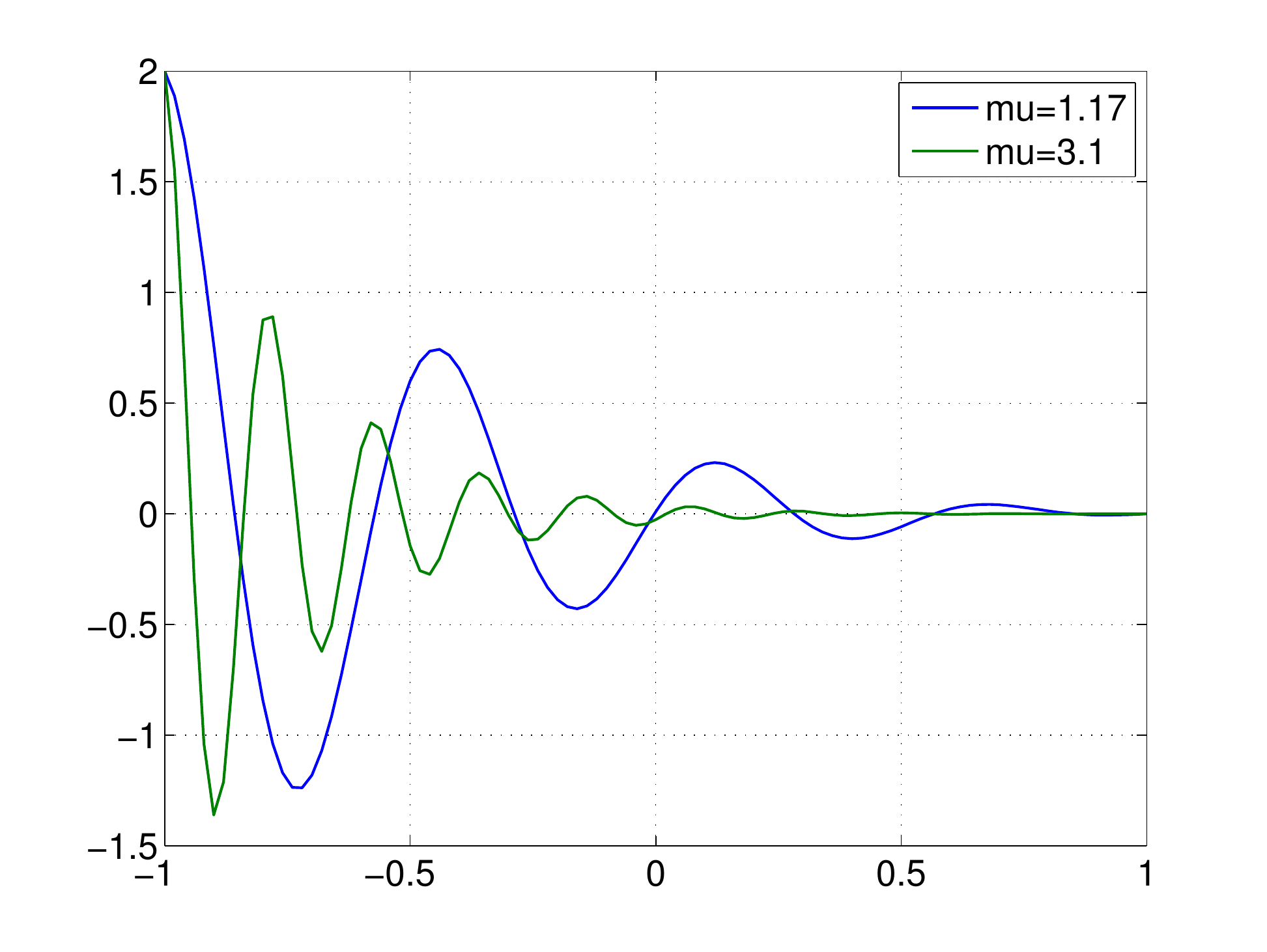}\includegraphics[scale=0.3]{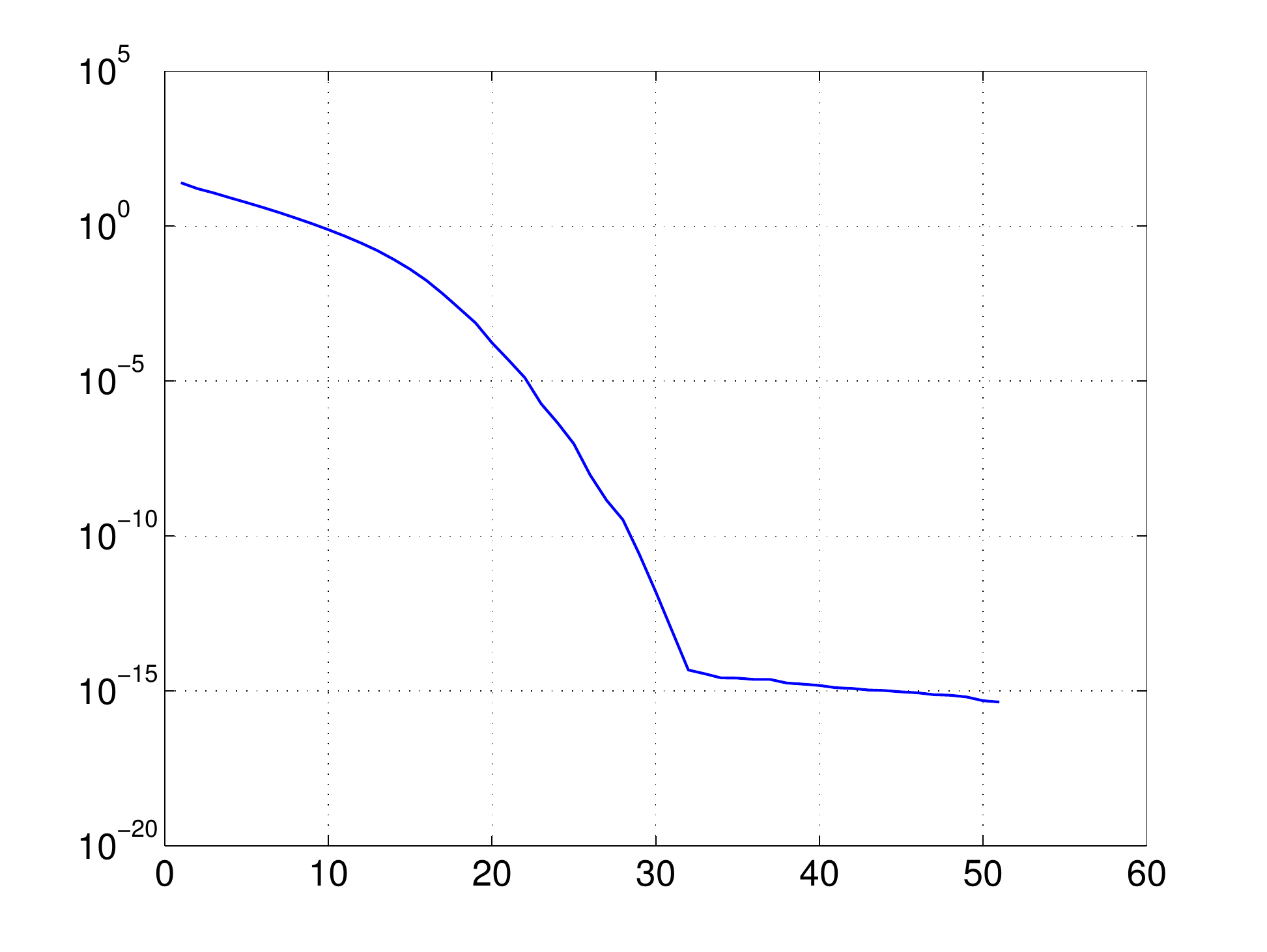}
\caption{Test 1: Plot of ${\bf f}(\mu)$ for $\mu=\{1.17, 3.1\}$ (left) and singular values of ${\bf f}(\mu)$ (right)}
\label{par_1}
\end{figure}

\begin{figure}[t]
\centering
\includegraphics[scale=0.3]{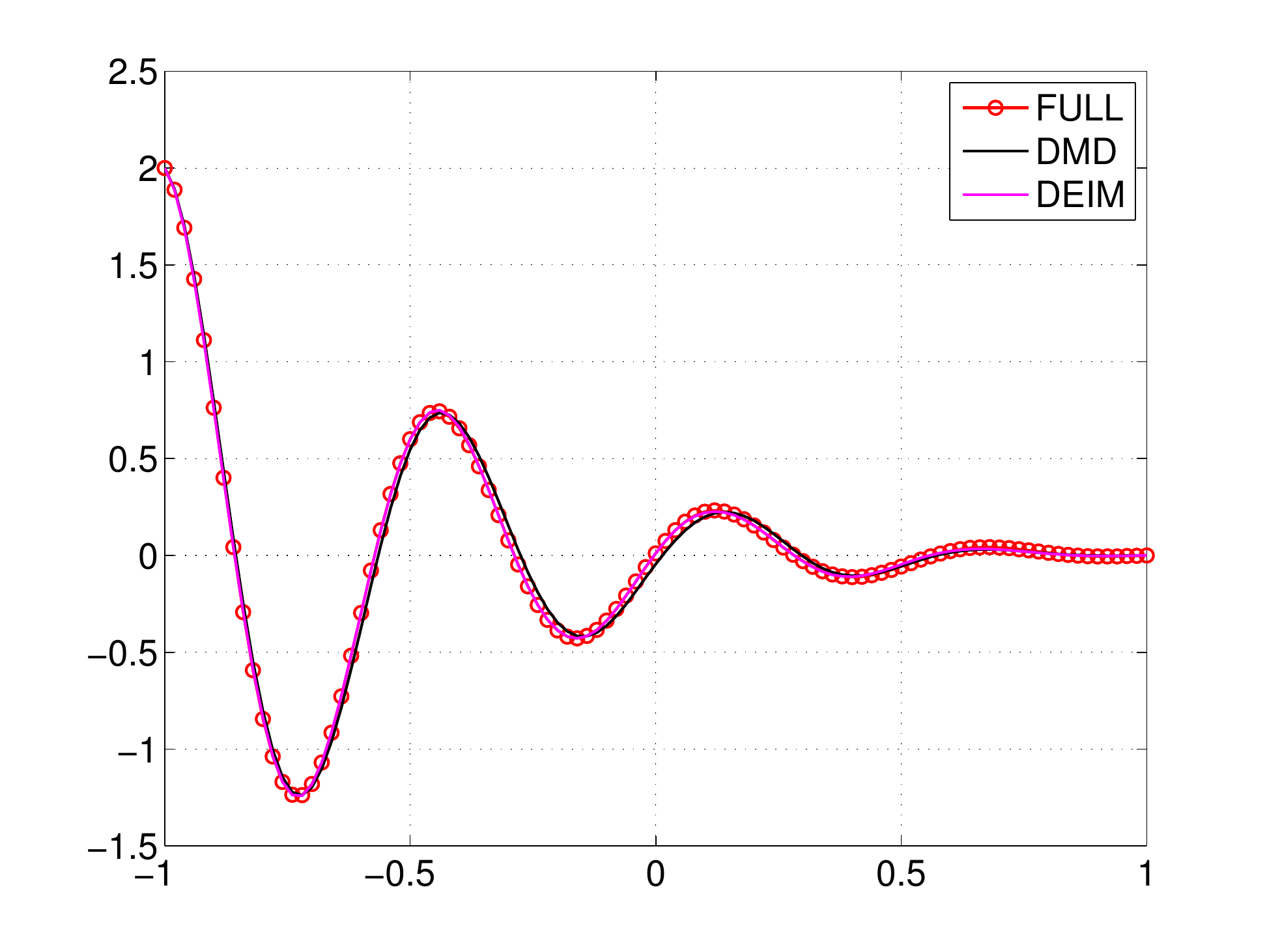}\includegraphics[scale=0.3]{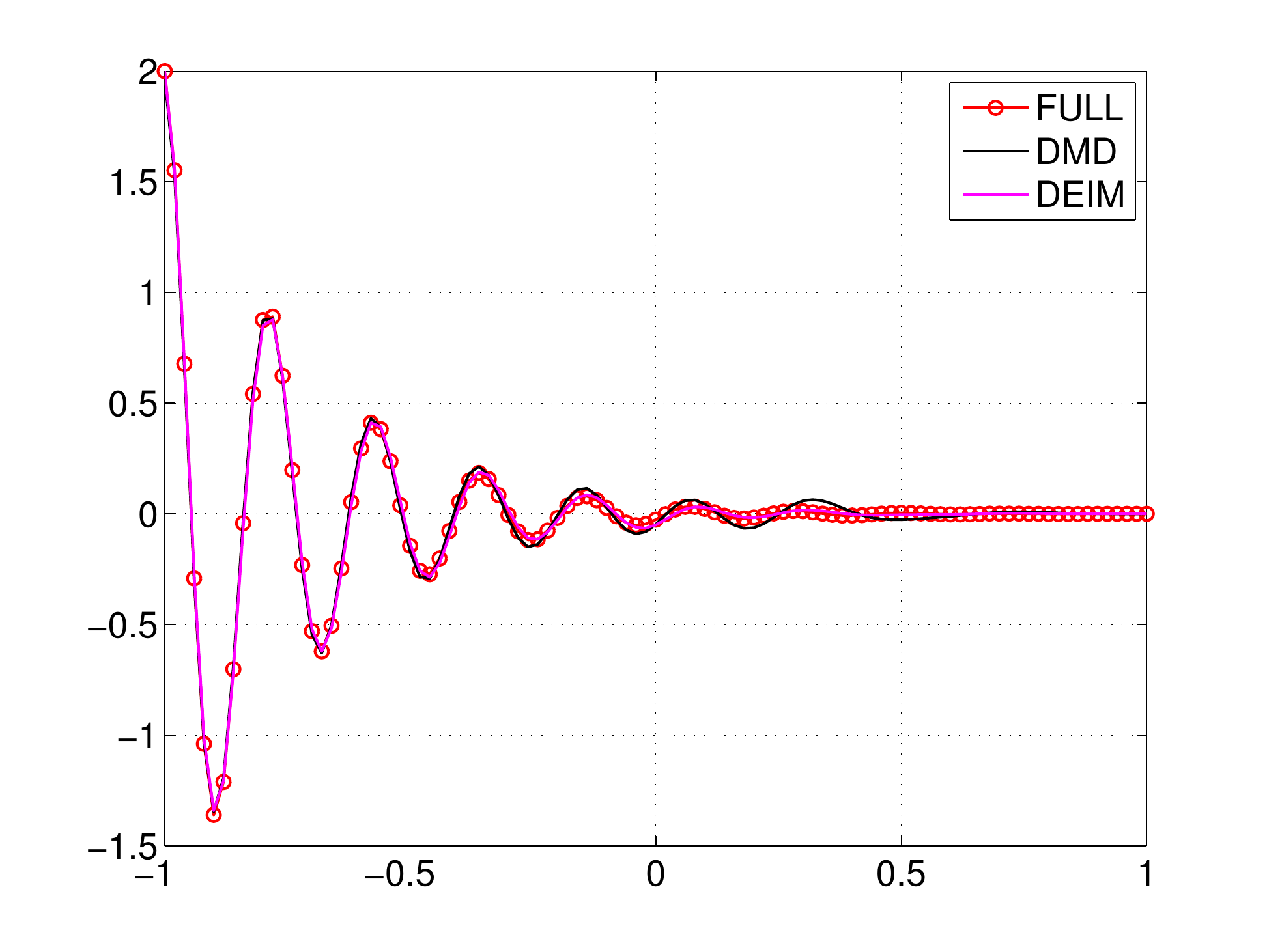}
\caption{Test 1: DMD and DEIM interpolation with $\mu=1.17$ (left) and $\mu=3.1$ (right).}
\label{par_3}
\end{figure}

\begin{figure}[t]
\centering
\includegraphics[scale=0.3]{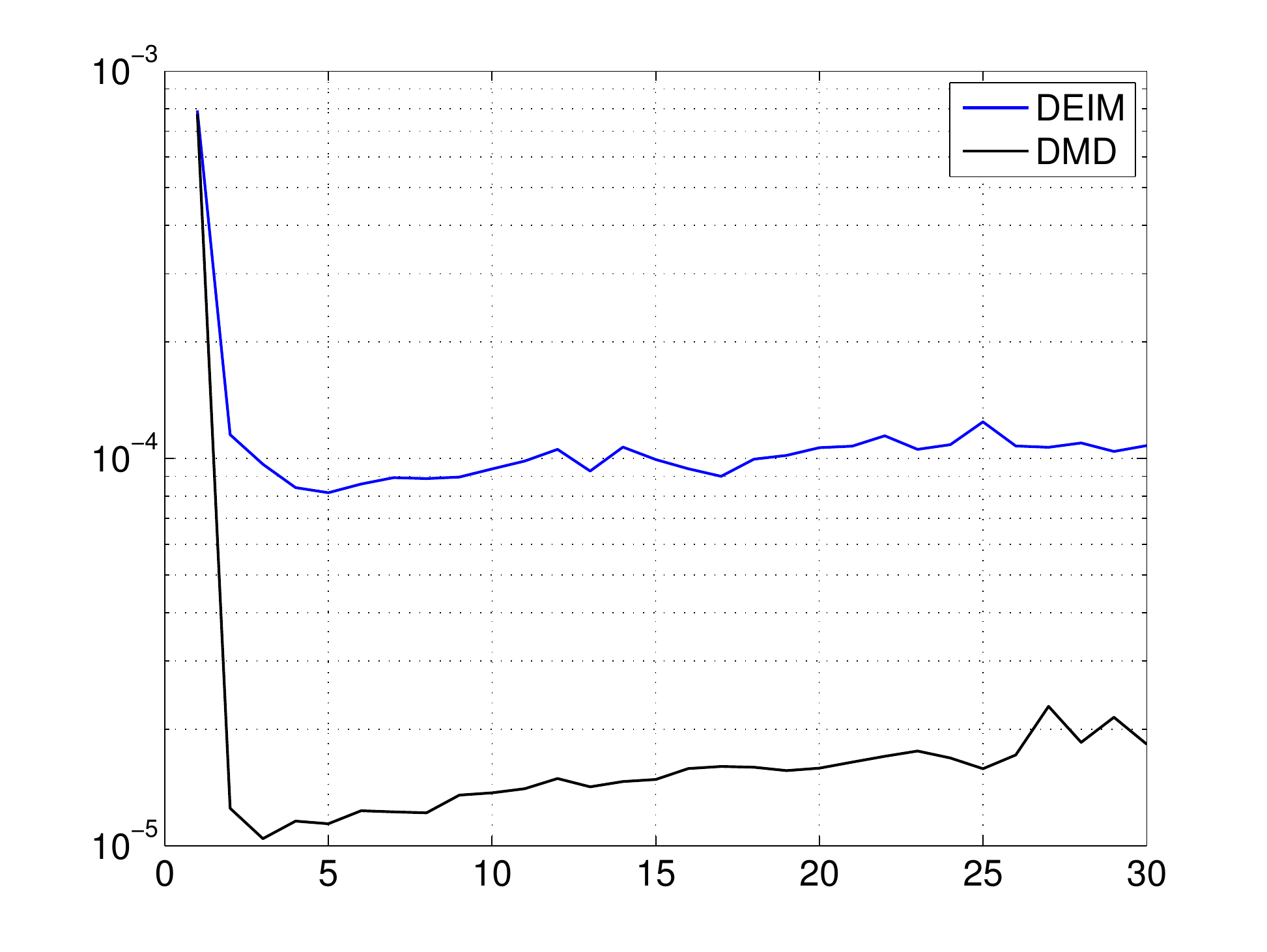}\includegraphics[scale=0.3]{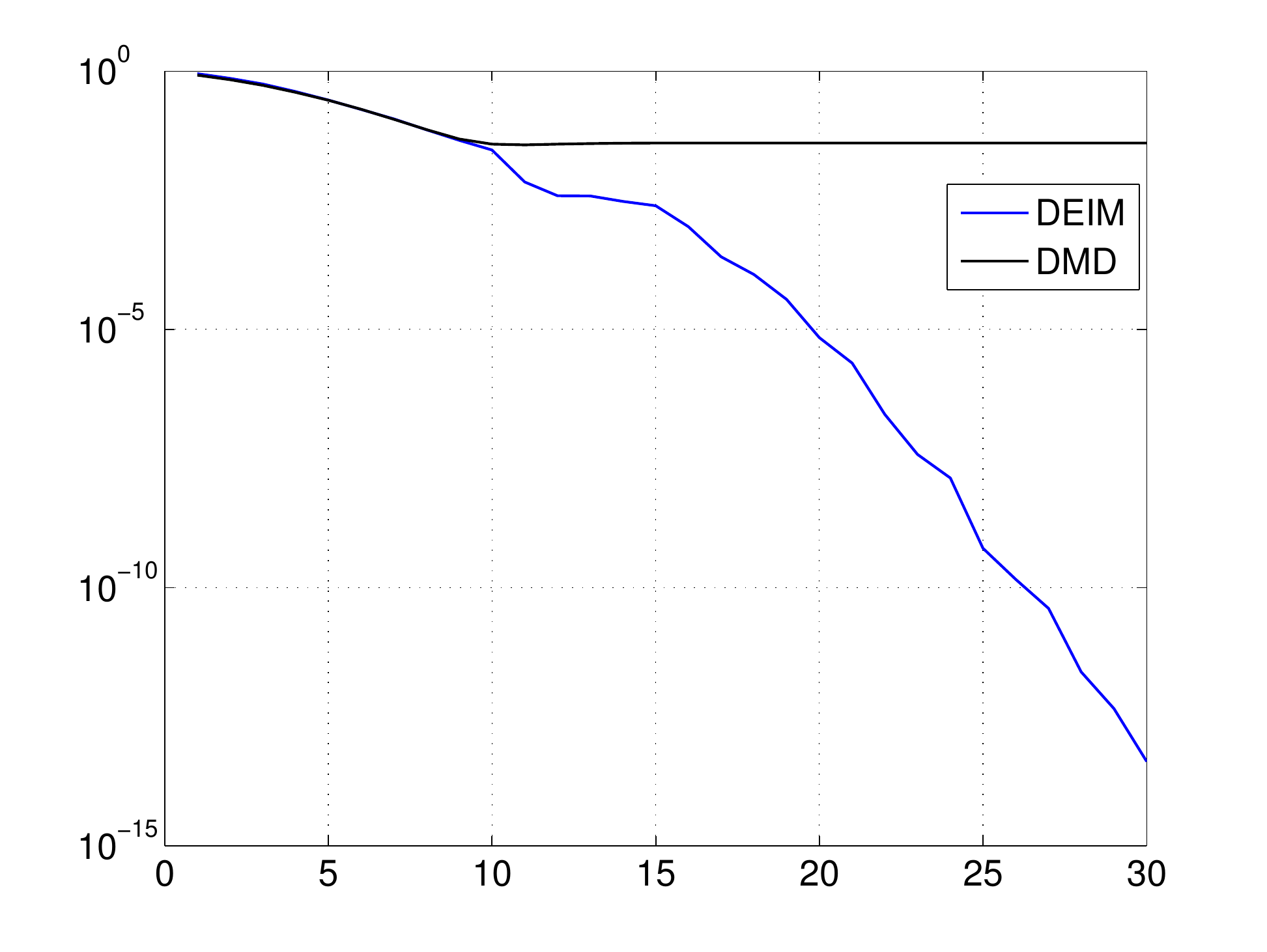}
\caption{Test 1: CPU time (left), relative error (right)}
\label{par_2}
\end{figure}

\paragraph{\bf Test 2: DMD-Galerkin approximation}
Although DMD is a well-known equation-free method, it also works in a Galerkin projection framework. In this subsection, we compare the performance of the POD method when the DMD method is integrated with the Galerkin method. We note that DMD basis function in the DMD-Galerkin projection are computed following Algorithm \ref{Alg_DMD} and then orthonormalized.
Let us consider the following linear advection-diffusion equation:
\begin{equation}\label{test:adv}
\begin{aligned}
y_t(x,t)+\theta y_x(x,t)&=0&& (x,t)\in[a,b]\times[0,T],\\
y(x,0)&=y_0(x)&& x\in [a,b],\\
y(a,t)&=0=y(b,t)&& t\in [0,T],
\end{aligned}
\end{equation}
where $a=0, b=4,T=3, \theta=1, y_0(x)=\sin(\pi x)$ if $0\leq x\leq 1$ and $0$ elsewhere. In order to lead \eqref{test:adv} to our general formulation \eqref{ode} we utilize a Finite Difference discretization with a spatial step $\Delta x=0.01$. The dimension of the problem is $n=399$. We note that in this case the mass matrix ${\bf M}$ is the identity matrix. In order to apply POD and DMD, we need to compute the snapshot set which is given by the temporal discretization of \eqref{ode} with an implicit Euler scheme and a temporal step size $\Delta t=0.01.$ The solution of equation \eqref{test:adv} builds the snapshot set and it is visualized on the left-side of Figure \ref{fig1:adv}. We also show the decay of the singular values of the snapshot matrix on the right of Figure \ref{fig1:adv}.

\begin{figure}[t]
\centering
\includegraphics[scale=0.3]{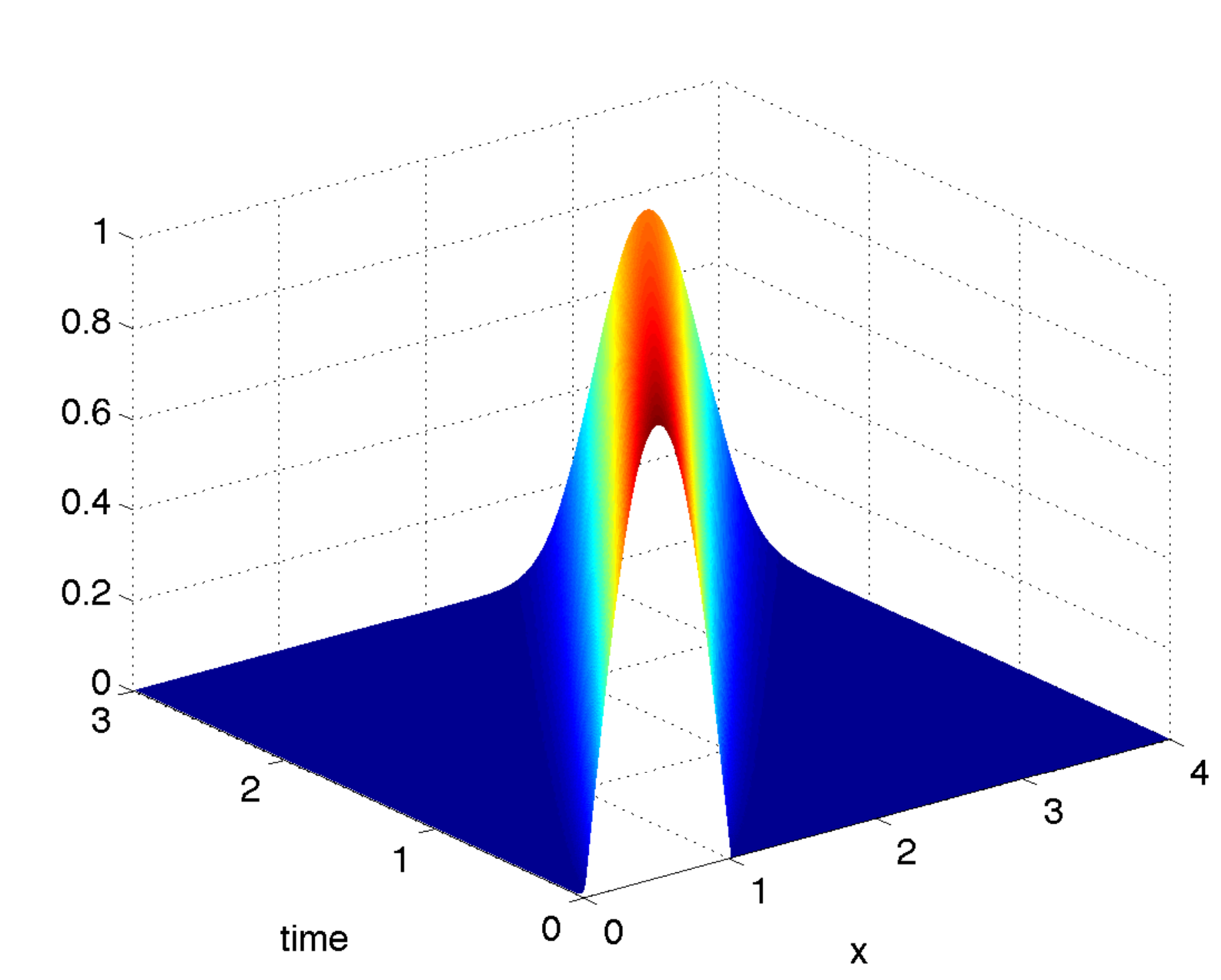}\includegraphics[scale=0.3]{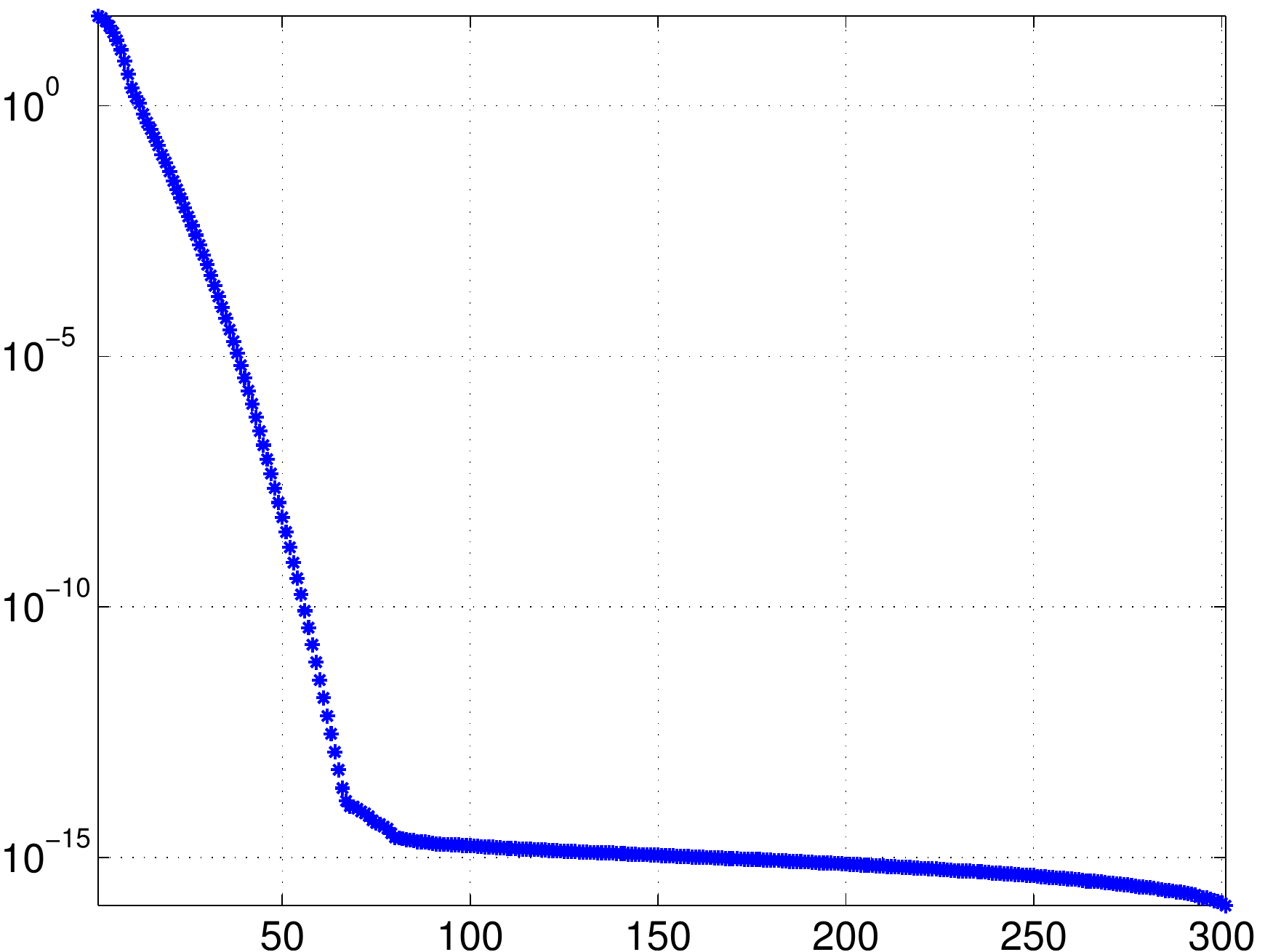}
\caption{Test 2: Solution of equation \eqref{test:adv} (left) and singular values (right)}
\label{fig1:adv}
\end{figure}

The POD-Galerkin has already been explained in Section \ref{Section2}. The DMD-Galerkin approach assumes that our solution can be written as ${\bf y}(t)\approx {\bf \Psi}^{\mbox{\tiny DMD}}{\bf y}^{\mbox{\tiny DMD}}(t)$. This assumption is very similar to \eqref{pod_ans}, but considers different basis functions. The reduced problem has the same form of \eqref{pod_sys}.
Figure \ref{fig1:comp} shows the results of model order reduction with POD (top), with DMD considered as a Galerkin projection method (middle), and DMD as a equation-free method. The first column refers to approximations of rank $5$, the second of rank $10$ and the third of rank $15$. As expected, if we increase the rank of the basis functions, we can easily see that the approximation gets better and better. It is well-known that advection dominated problems have a high variability during time evolution and it is difficult to capture the dynamics with only a few basis functions.

\begin{figure}[t]
\begin{center}
\includegraphics[scale=0.22]{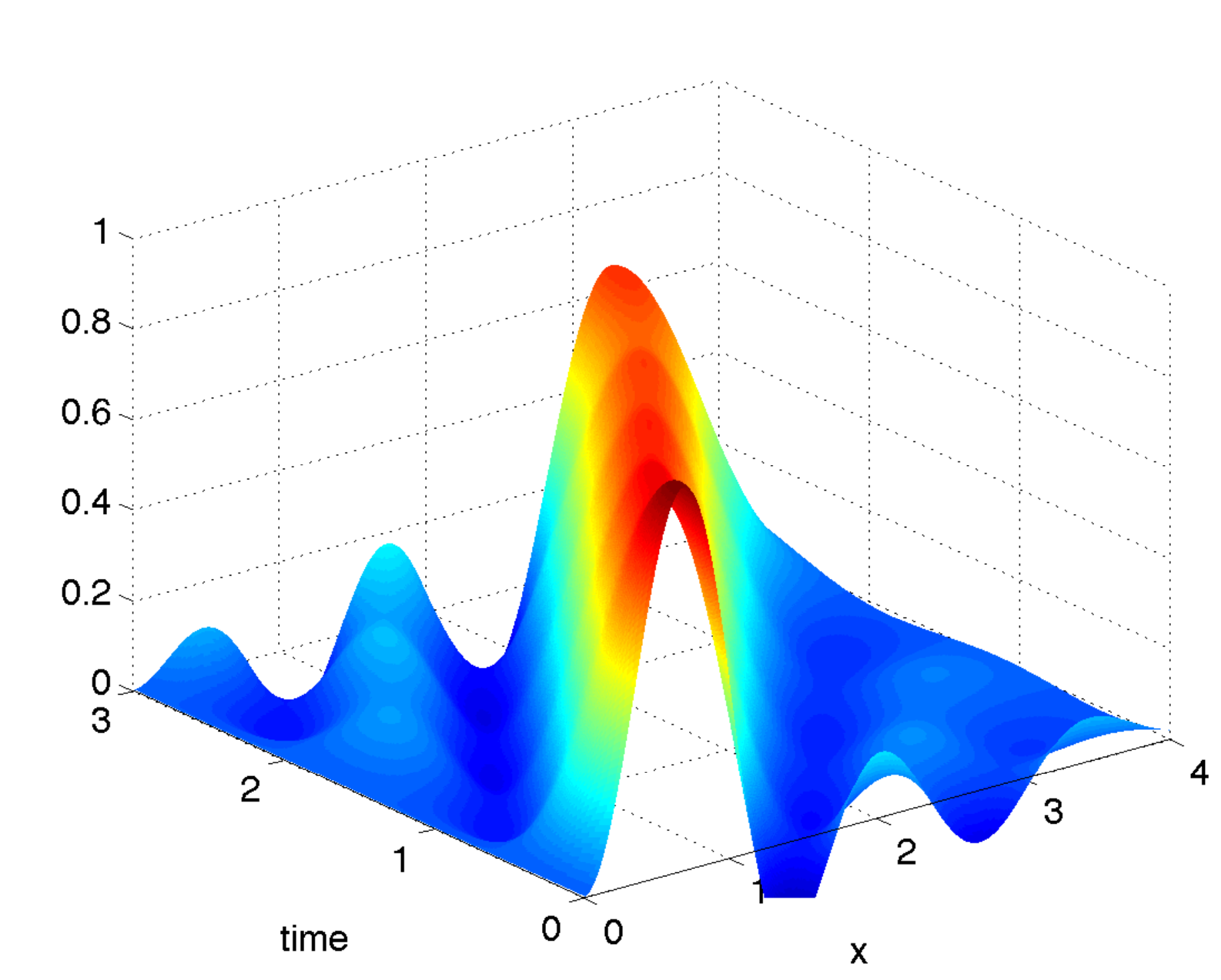}
\includegraphics[scale=0.22]{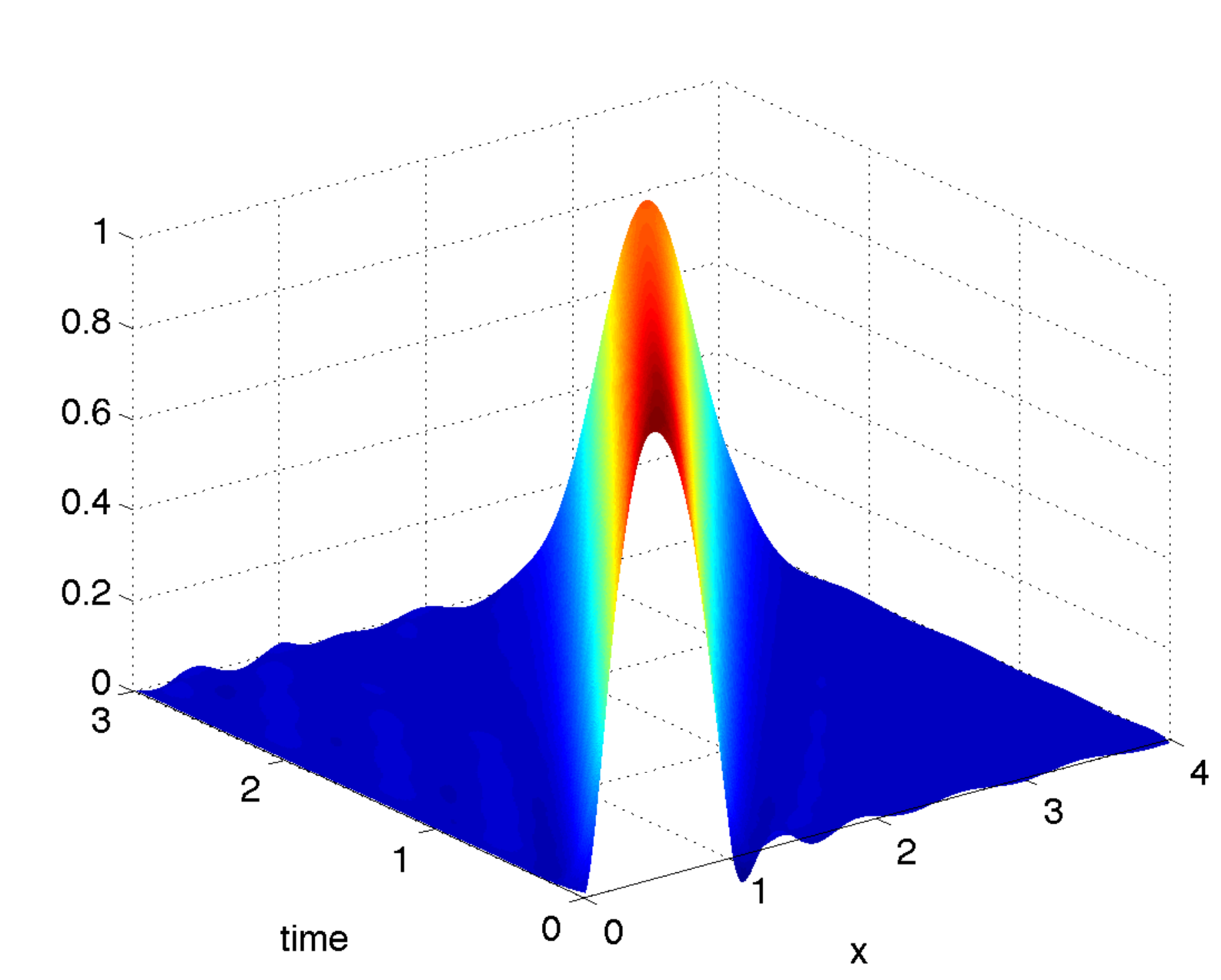}
\includegraphics[scale=0.22]{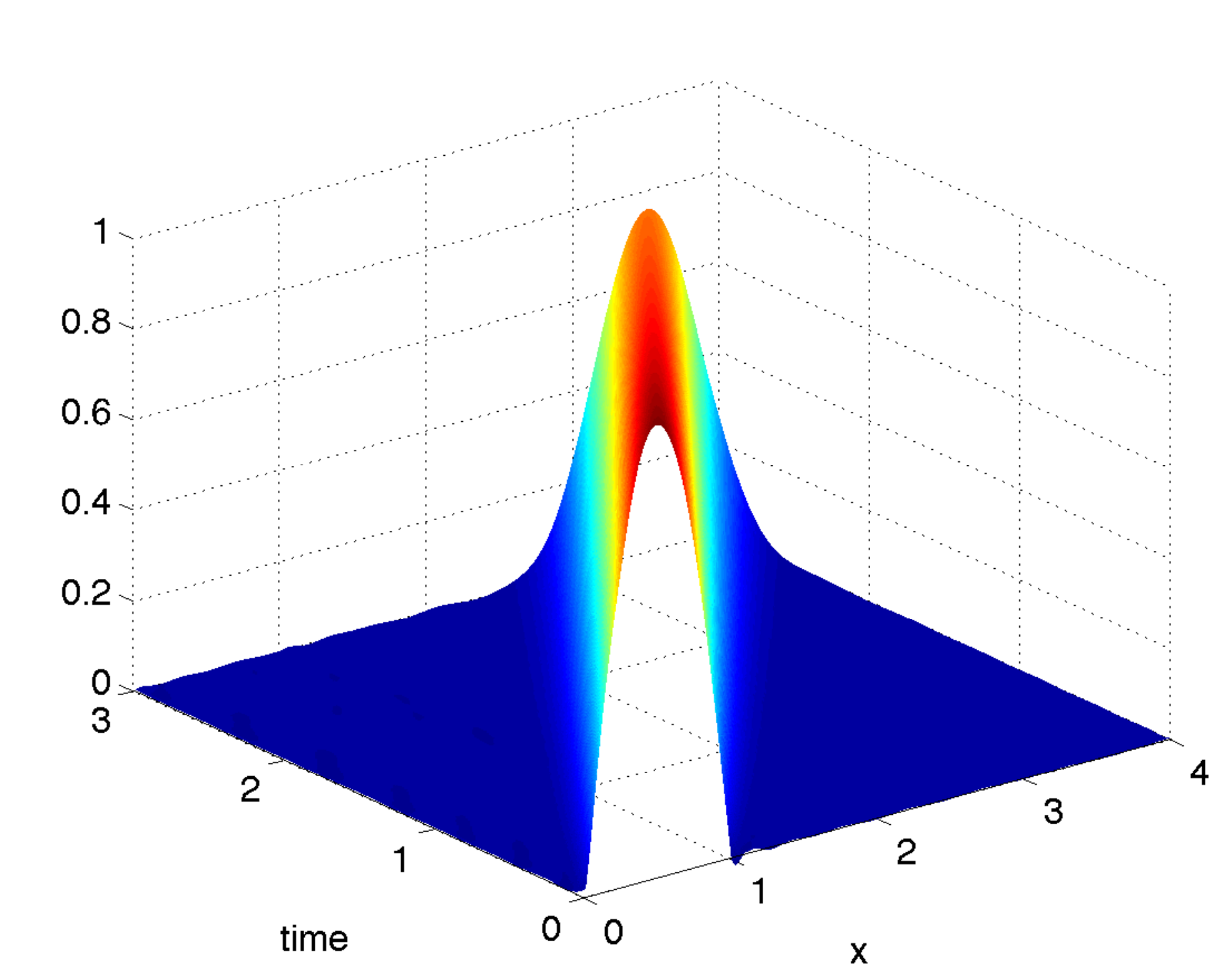}\\
\includegraphics[scale=0.22]{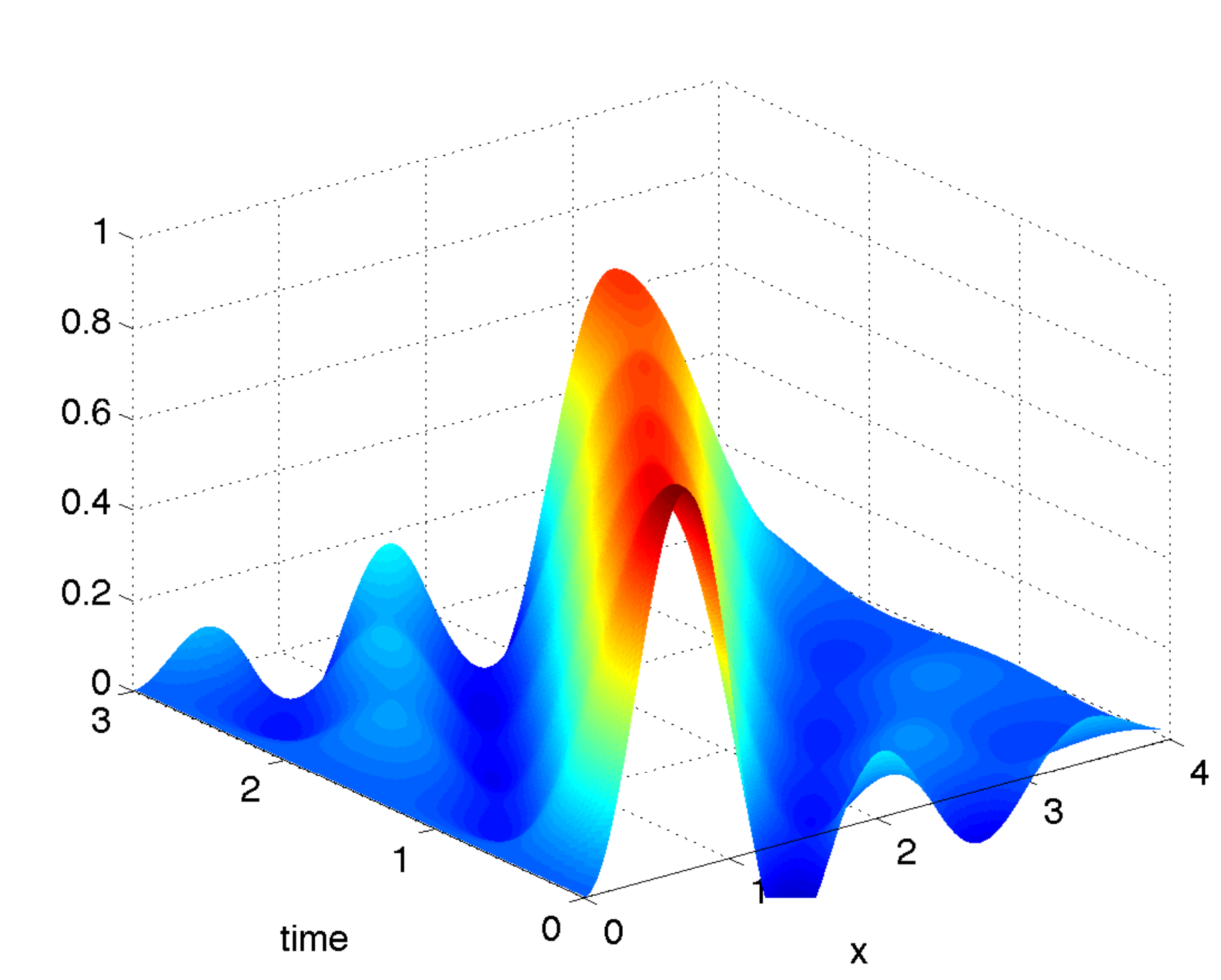}
\includegraphics[scale=0.22]{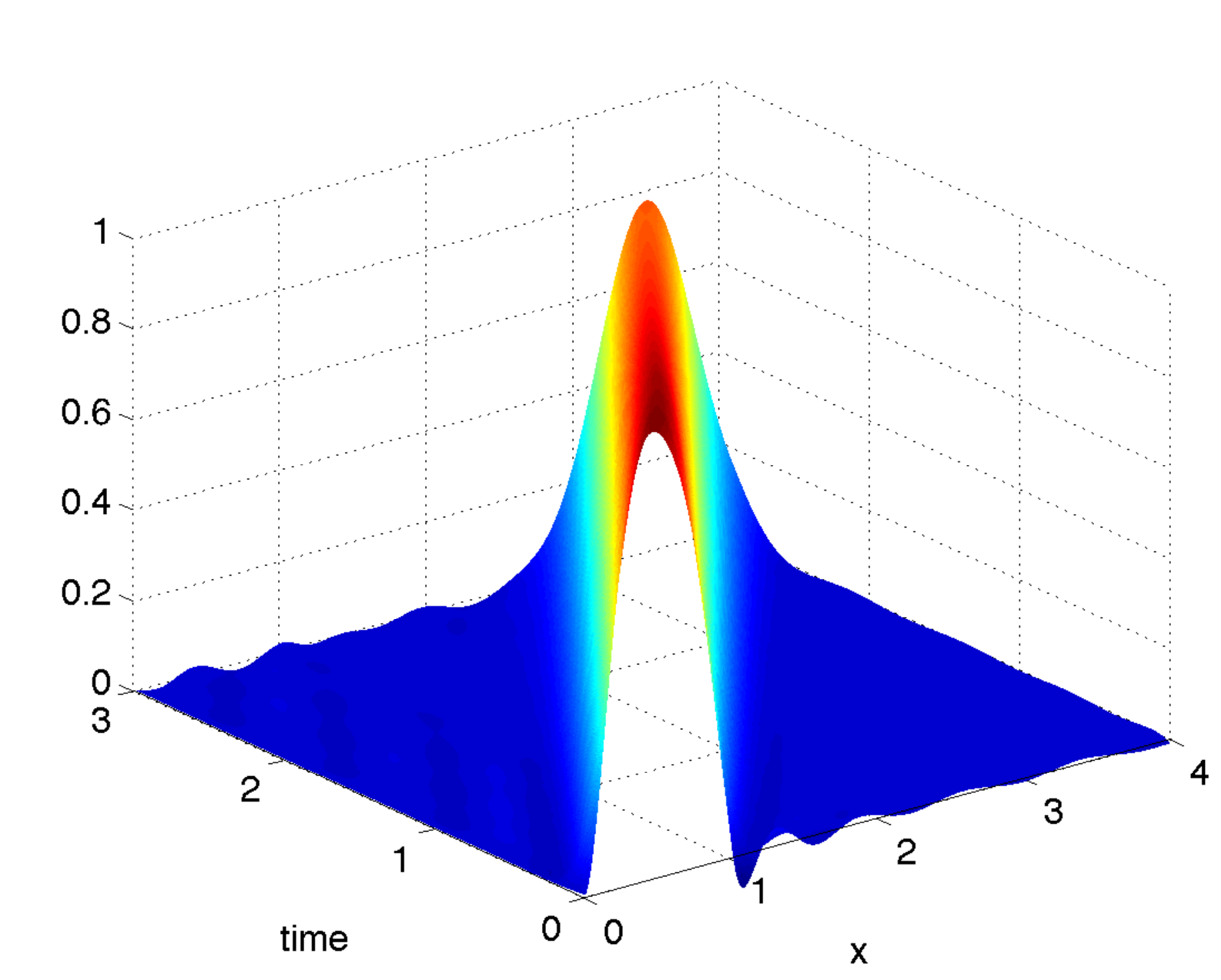}
\includegraphics[scale=0.22]{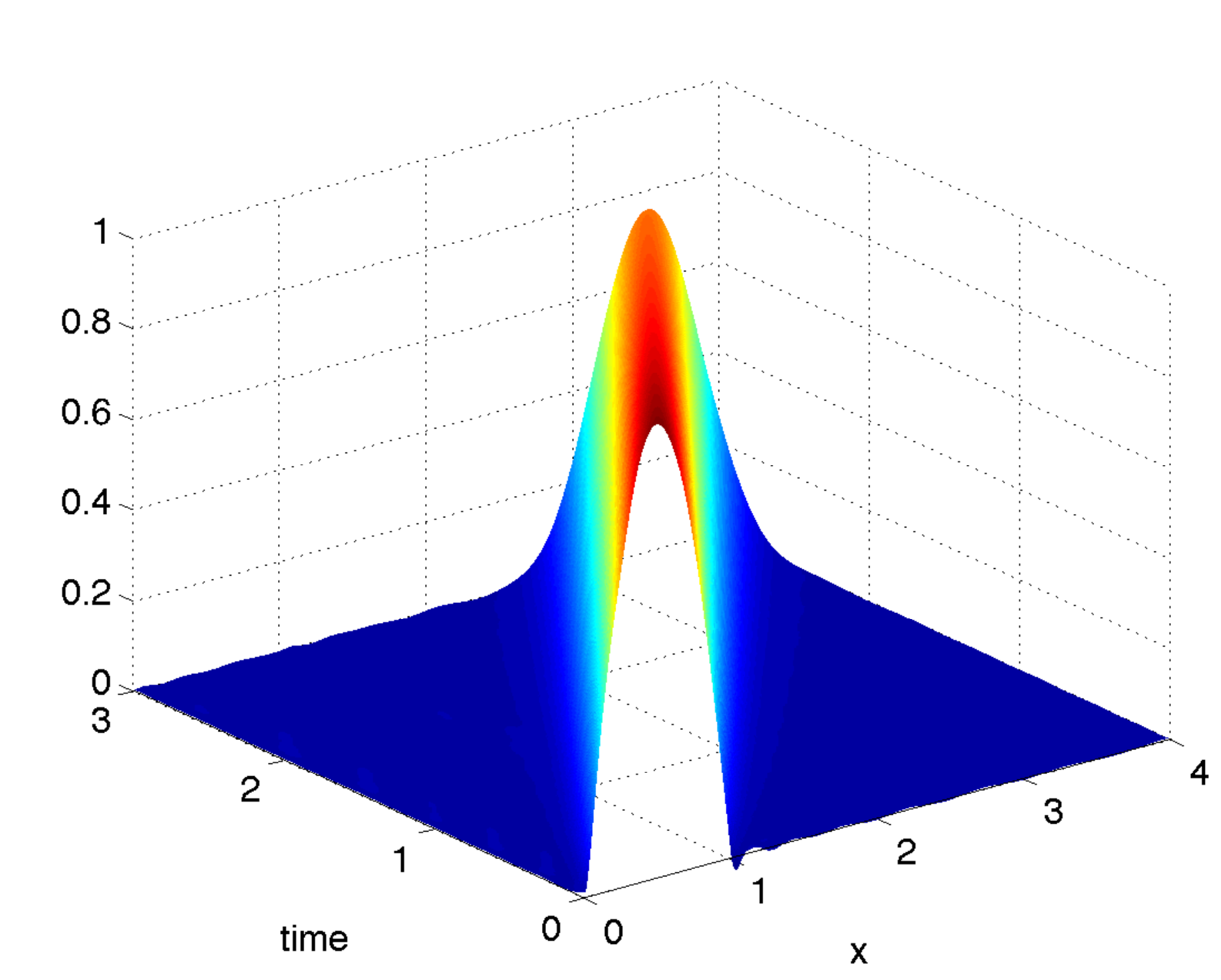}\\
\includegraphics[scale=0.22]{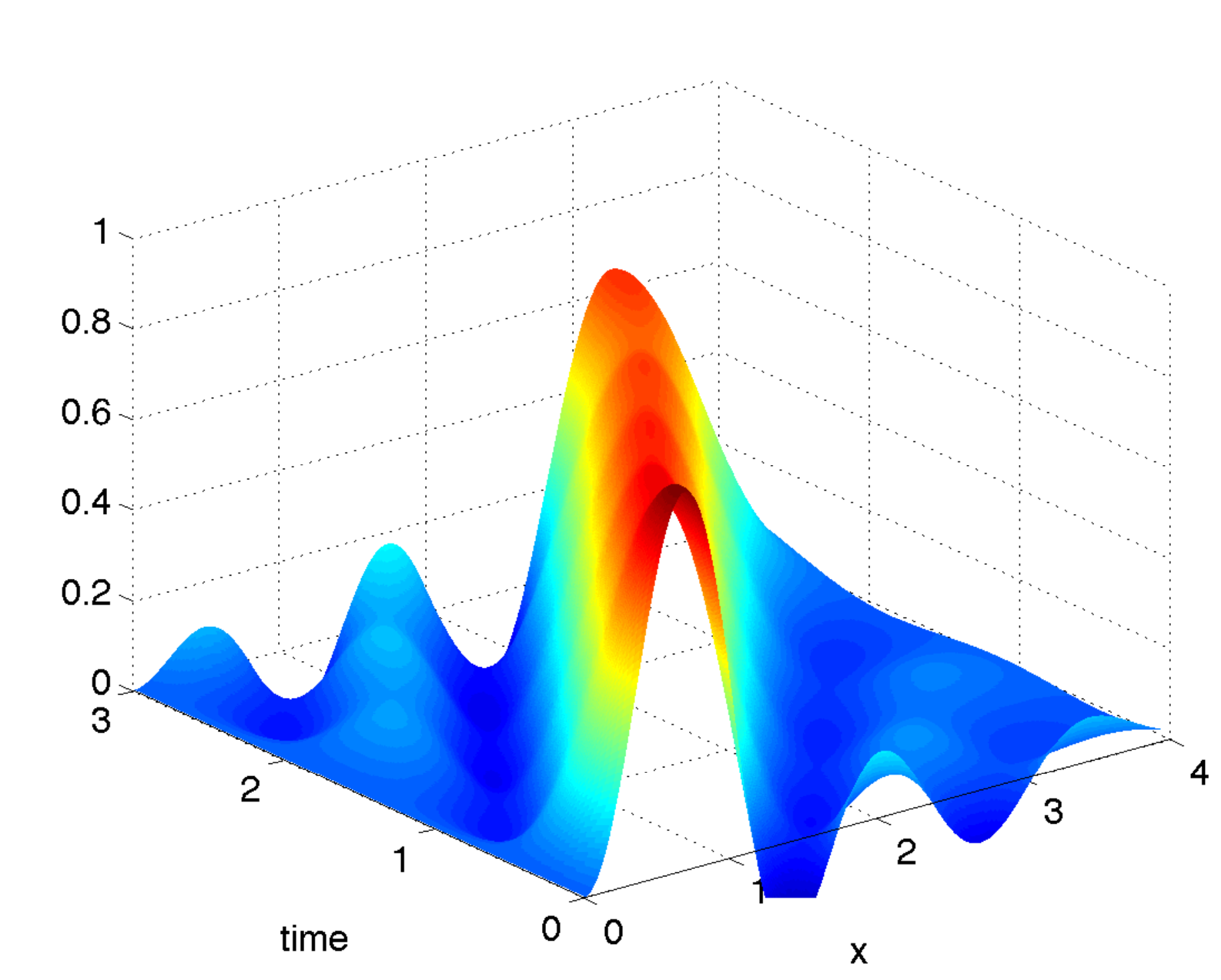}
\includegraphics[scale=0.22]{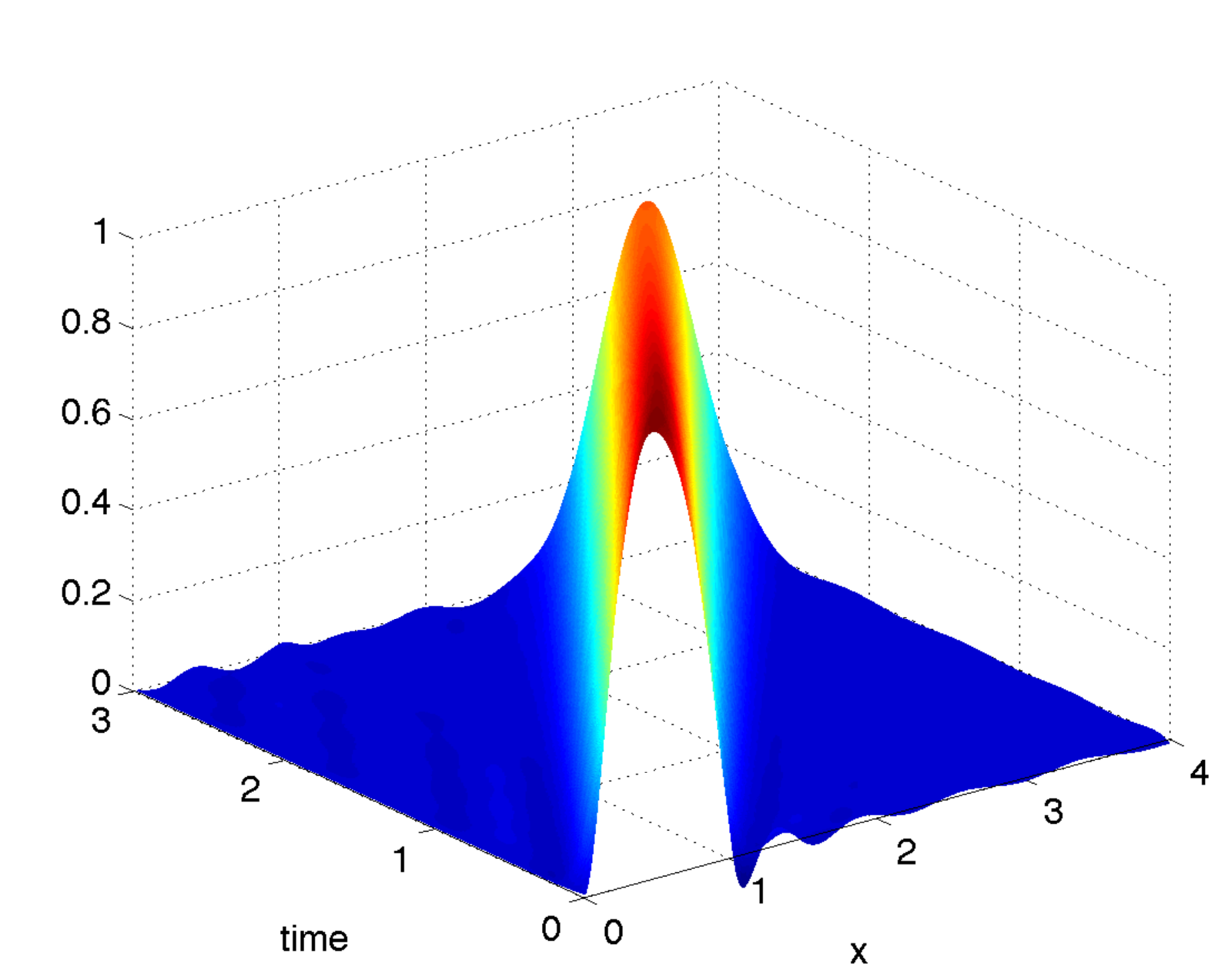}
\includegraphics[scale=0.22]{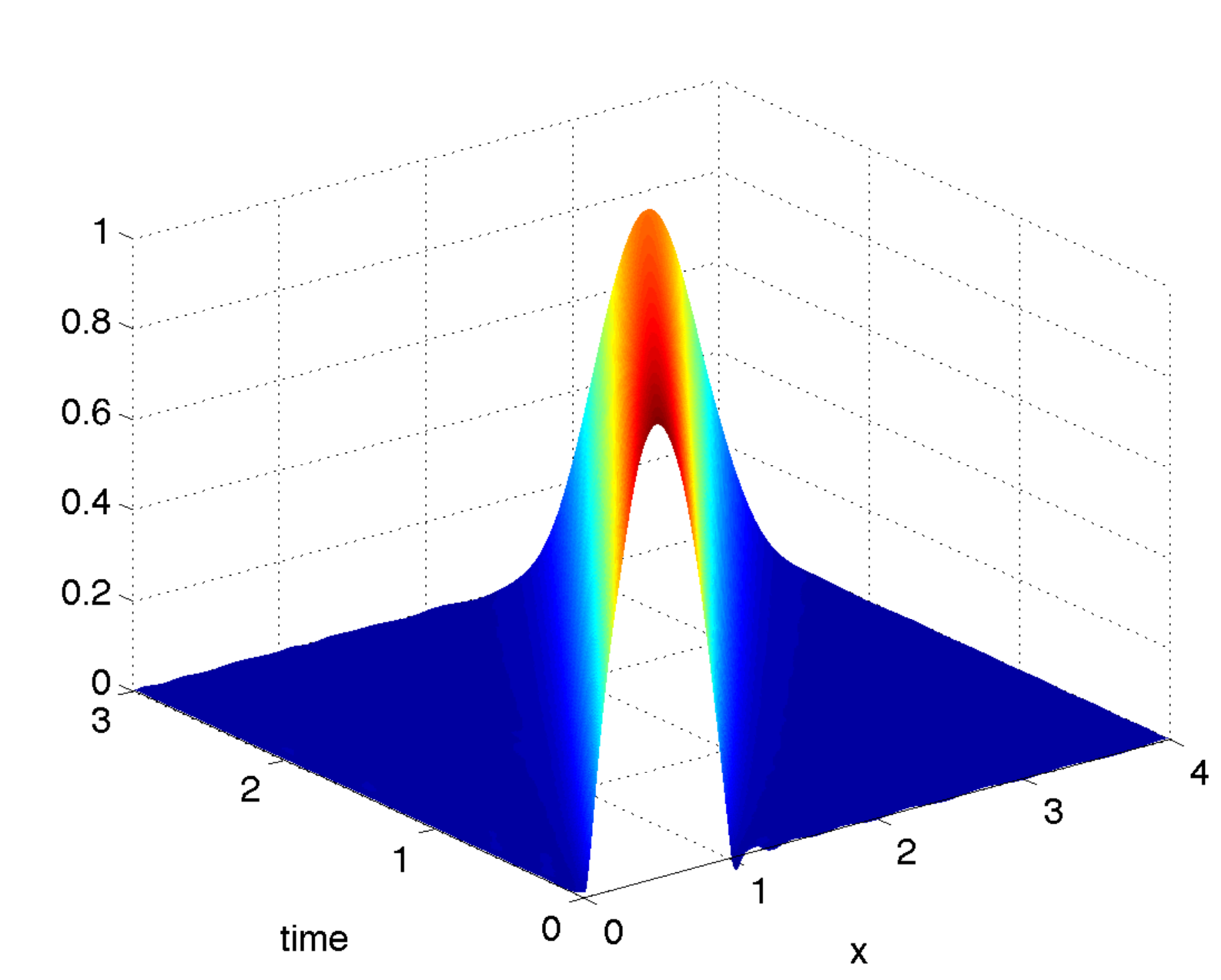}\\
\end{center}
\caption{Test 2: Reduced approximation with rank=$\{5,10,15\}.$ POD approximation (top), DMD-Galerkin (middle), DMD (data-driven) (bottom)}
\label{fig1:comp}
\end{figure}

The error analysis in the right panel of Figure \ref{fig1:err} confirms our heuristic expectations. Here we compute the relative error with respect to the Frobenius-norm where we consider as truth the solution of the governing
equations approximated by a Finite Difference scheme. The error decays as soon as we increase the dimension of the reduced model, in particular, the POD method always performs better. On the other hand, it is, in general, hard to see significant difference between the DMD-Galerkin and the DMD data-driven method as shown in Figure \ref{fig1:err}. This phenomena has been observed even for both linear or nonlinear problems.  In fact, this error analysis brings us to the idea of working with a data-driven method for the approximation of nonlinear dynamical systems. For the sake of clarity, we also show the first two modes in Figure \ref{fig1:err} where we can clearly observe that the DMD basis functions oscillate more than the POD modes.

\begin{figure}[t]
\begin{center}
\includegraphics[scale=0.22]{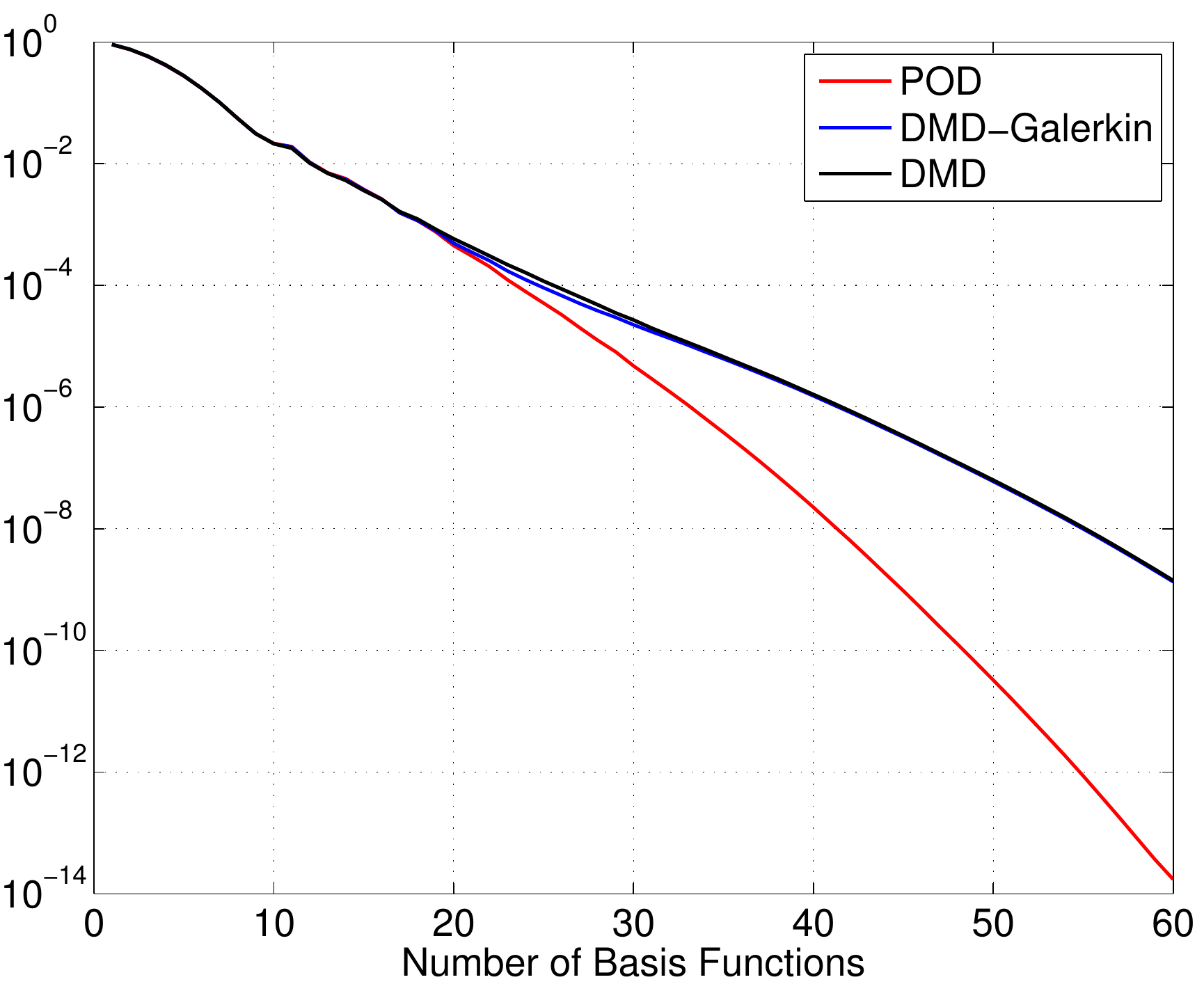}
\includegraphics[scale=0.22]{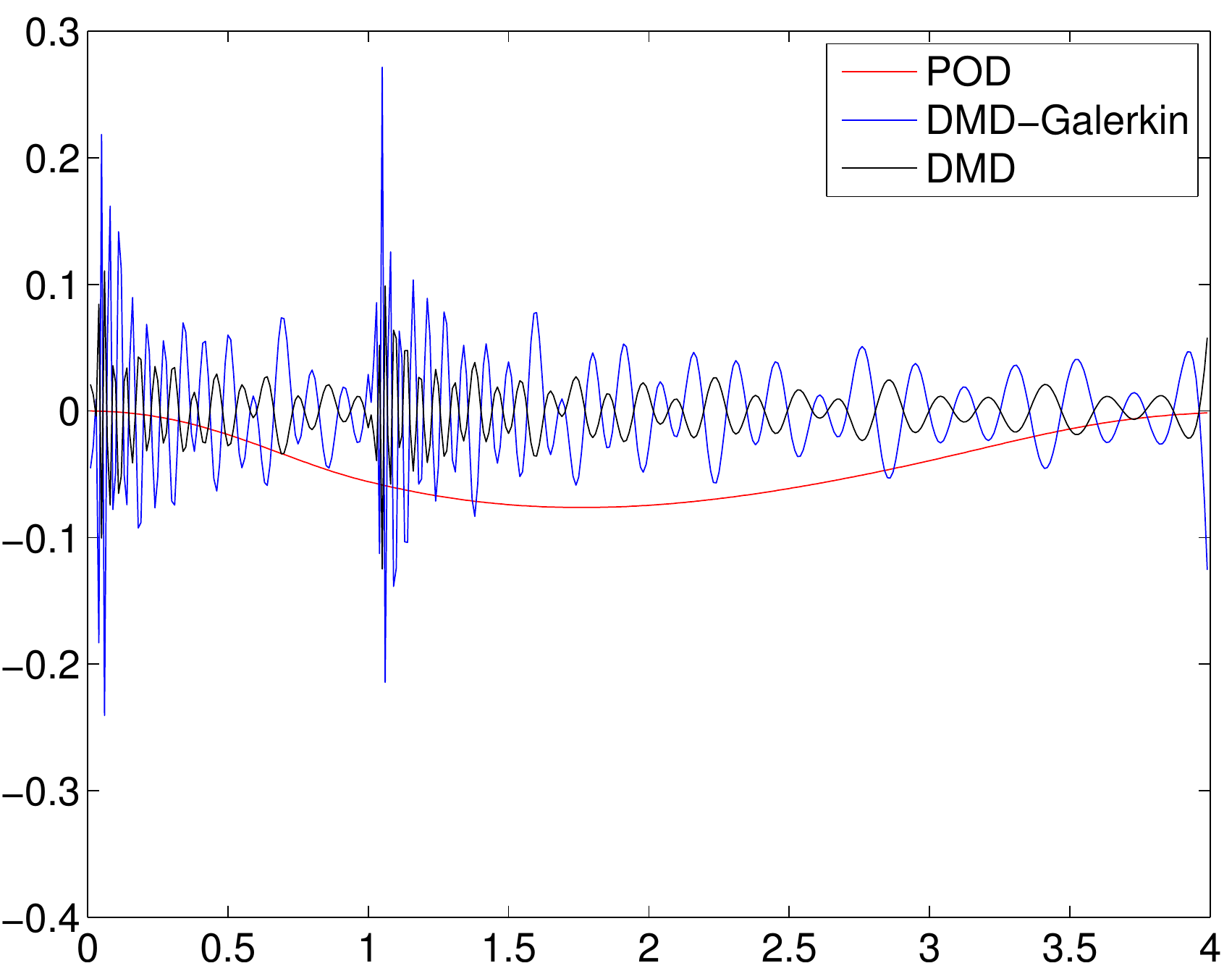}
\includegraphics[scale=0.22]{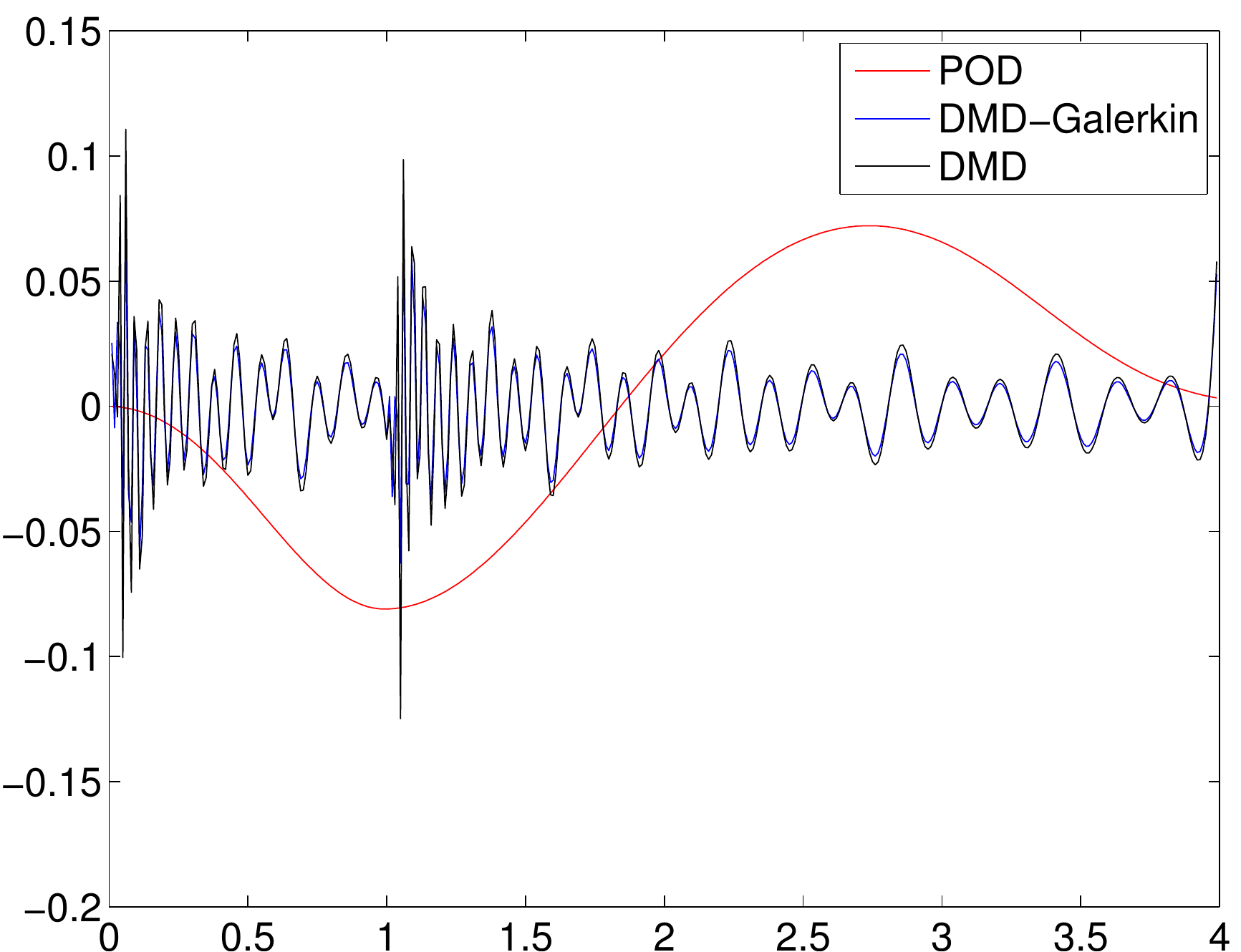}\\
\end{center}
\caption{ Error analysis  with respect to the Frobenius norm (left), first mode (middle), second mode (right).}
\label{fig1:err}
\end{figure}


\section{Coupling POD and DMD for nonlinear problem}
\label{Section4}
\setcounter{section}{4}
\setcounter{equation}{0}
\setcounter{theorem}{0}
\renewcommand{\theequation}{\arabic{section}.\arabic{equation}}

This section focuses on the approximation of a nonlinear problem by means of model order reduction. As discussed in Section 2, the use of POD basis functions does not lead to a surrogate model which is independent of the full dimension of the problem (see \eqref{pod_sys}).  We advocate an alternative method to EIM/DEIM by working with the DMD algorithm for evaluating the nonlinear term in \eqref{ode}. 
As already discussed, the snapshot measurements used in DMD approximate the dynamics and predict the future state. The use of DMD, in this work, concerns the approximate of the nonlinearity ${\bf f}(t, {\bf y}(t))$ of the dynamical system \eqref{ode}.

To begin with let us collect snapshots from the system $\{{\bf y}(t_0),\dots, {\bf y}(t_m)\}$ for some given time instances $\{t_0,\dots,t_m\}$ and compute the POD basis functions of rank $\ell$. Then, we need to collect snapshots for the nonlinearity $\{{\bf f}(t_0,{\bf y}(t_0)),\dots, {\bf f}(t_m, {\bf y}(t_m))\}$ and divide them into two different sets as explained in Section 3. We apply the DMD algorithm (see Algorithm \ref{Alg_DMD}) to the nonlinear measurements. The DMD approximation of the nonlinearity reads:

\begin{equation}
{\bf f}^{\mbox{\tiny DMD}}(t,y(t))=  \sum_{i=1}^k b_i {\boldsymbol \psi}^{\mbox{\tiny DMD}}_i \exp(\omega_i t) \, ,
\label{eq:omegaf}
\end{equation}
where ${\boldsymbol \psi}_i^{\mbox{\tiny DMD}}$ are the DMD basis functions of rank $k$ related to the nonlinear function ${\bf f}(t,{\bf y}(t))$, $b_i$ is the initial condition and $\omega_i$ are the eigenvalues of the linear matrix {\bf$\tilde{A}_{\bf y}$}. With compact notation we obtain:
\begin{equation}
\tilde{\bf f}^{\mbox{\tiny DMD}}(t,y(t))\approx  {\bf \Psi}^{\mbox{\tiny DMD}}\diag (e^{\omega^{\mbox{\tiny DMD}}t})b ,
\label{eq:omegafcompac}
\end{equation}
where $b=({\bf \Psi}^{\mbox{\tiny DMD}})^\dag {\bf f}(t_1,{\bf y}(t_1))\in\R^k$, $\diag (e^{\omega^{\mbox{\tiny DMD}}t})b\in\R^k$ represents the reduced approximation of the data in terms of the DMD modes. As we can see from \eqref{eq:omegafcompac} the nonlinearity is approximated by a DMD representation and
no further evaluation of the nonlinearity is required. This circumvents the DEIM selection of the interpolation points.
If we plug the approximation of the nonlinearty \eqref{eq:omegafcompac} into the POD system \eqref{pod_sys} we get the following reduced system:

\begin{equation}\label{pod_sysdmd}
\left\{\begin{array}{l}
{\bf M}^\ell\dot{{\bf y}}^\ell(t)={\bf A}^\ell {\bf y}^\ell(t)+{\bf \Psi}^T{\bf \Psi}^{\mbox{\tiny DMD}}\diag (e^{\omega^{\mbox{\tiny DMD}}t})b \\
{\bf y}^\ell(0)={\bf y_0}^\ell.
\end{array}\right.
\end{equation}
Let us analyze the dimension of the new reduced dynamical system \eqref{pod_sysdmd}. The matrix ${\bf M}^\ell, {\bf A}^\ell\in\R^{\ell\times\ell}$ have the same dimension of the POD system. The quantity $({\bf \Psi}^{\mbox{\tiny POD}})^T{\bf \Psi}^{\mbox{\tiny DMD}}\in\R^{\ell\times k}$ is independent of the dimension of the full system, and $\diag (e^{\omega^{\mbox{\tiny DMD}}t})b\in\R^k$. Even for this formulation we are able to build a surrogate model which does not depend on the dimension of the original system. Moreover, in this formulation we do not have to evaluate the nonlinearity further, which gives an important speed up in the efficiency of the formulation. As in the DEIM case some quantities can be precomputed offline. Of course, this method is closely related to the snapshot set, and approximate the  nonlinear term with a linear regression operator. The algorithm is summarized in \ref{Alg_PODDMD}

\begin{algorithm}
\caption{POD-DMD}
\label{Alg_PODDMD}
\begin{algorithmic}[1]
\REQUIRE Snapshots $\{{\bf y}_(t_0),\ldots,{\bf y}(t_m)\}$, $\ell$ number of POD modes, $k$ number of DMD modes
\STATE Compute the POD basis function $\{{\boldsymbol{\psi}}_i\}_{i=1}^\ell$ of rank of $\ell$.
\STATE Compute nonlinear snapshots $\{{\bf f}(t_0,{\bf y}(t_0)),\dots, {\bf f}(t_m, {\bf y}(t_m))\}$
\STATE Set ${\bf Y}=[{\bf f}(t_0,{\bf y}(t_0)),\dots, {\bf f}(t_{m-1}, {\bf y}(t_{m-1}))]$
\STATE Set ${\bf Y}'=[{\bf f}(t_1,{\bf y}(t_1)),\dots, {\bf f}(t_m, {\bf y}(t_m)]$,
\STATE Compute DMD modes following Algorithm \ref{Alg_DMD}
\STATE Set and integrate equation \eqref{pod_sysdmd}
\STATE Project back full solution
\end{algorithmic}
\end{algorithm}

The POD-DMD method has one significant advantage:  computational speed. Indeed, as we will show, the
POD-DMD algorithm is significantly faster than the POD-DEIM method in approximating the nonlinear terms
in the model reduction. Indeed, the computational efficiency for this task is improved by an order of magnitude
or more. The drawback of the method is that the DMD modes produced for the low-rank projection are not orthogonal. Thus the convergence of the solution to the original high-dimensional system plateaus and the error is not reduced beyond a prescribed point. We hope to fix this problem in future work by potentially orthogonalizing the DMD modes.

%
%
%

\section{Numerical Tests}
\label{Section5}
\setcounter{section}{5}
\setcounter{equation}{0}
\setcounter{theorem}{0}
\setcounter{algorithm}{0}
\renewcommand{\theequation}{\arabic{section}.\arabic{equation}}

In this section we present our numerical tests. In our numerical computations we use the finite difference method to reduce a partial differential equation into the form \eqref{ode} and integrate the system with a semi-implicit scheme. All the numerical simulations reported in this paper are performed on an iMac with an Intel Core i5, 2.7Ghz and 8GB RAM using MATLAB R2013a. 

In the following numerical examples we built different surrogate models, such as POD, POD-DEIM and POD-DMD and compared their performance in term of CPU time and the error with respect to a reference solution computed by the Finite Difference approach. The tests consider three types of equations.

\paragraph{\bf Test 3: Semi-Linear Parabolic Equation}

Let us consider the following equation:
\begin{equation}\label{prb_test1}
\begin{aligned}
y_t(x,t)-\theta \Delta y(x,t)+\mu(y(x,t)-y^3(x,t))&=0&& (x,t)\in \Omega\times[0,T],\\
y(x,0)&=y_0(x)&&x\in \Omega,\\
y(a,t)&=0=y(b,t)&& t\in [0,T],
\end{aligned}
\end{equation}
where $\Omega=[0,1]\times[0,1], T=3, y_0(x)=0.1$ if $0.1\leq x_1x_2\leq0.6$ and $0$ elsewhere. The POD basis vectors are built upon 100 equidistant snapshots. The FD discretization yields a system of ODEs of the same form as \eqref{ode}. The solution of this equation generates a stationary solution $y(x,t)\equiv 1$ for large $t$ as shown in Figure \ref{test1:sol}.
Figure \ref{test1:svd} shows the singular values of the snapshot set, and of the nonlinear term in \eqref{prb_test1}.

\begin{figure}[H]
\begin{center}
\includegraphics[scale=0.25]{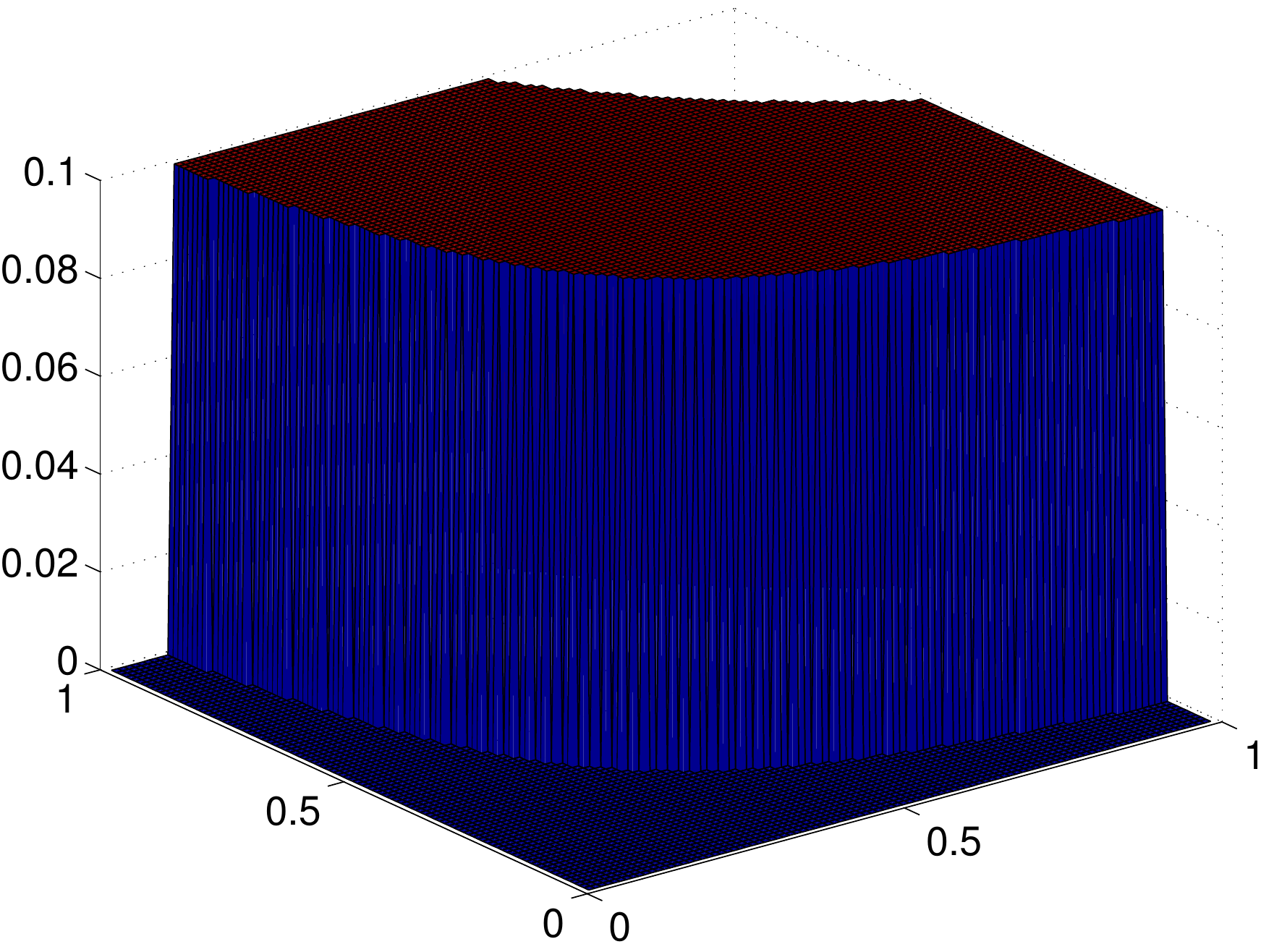}
\includegraphics[scale=0.25]{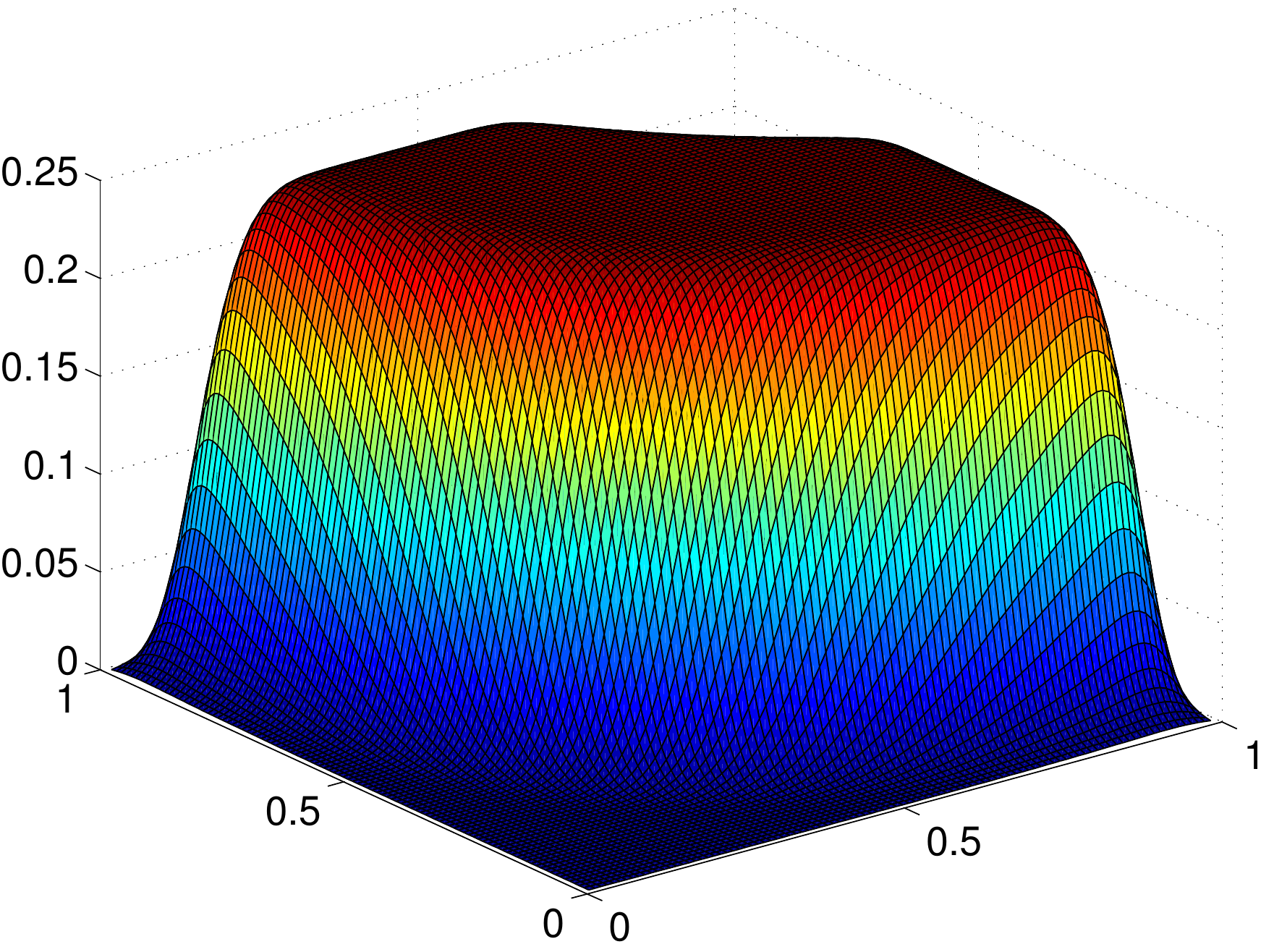}\\
\includegraphics[scale=0.25]{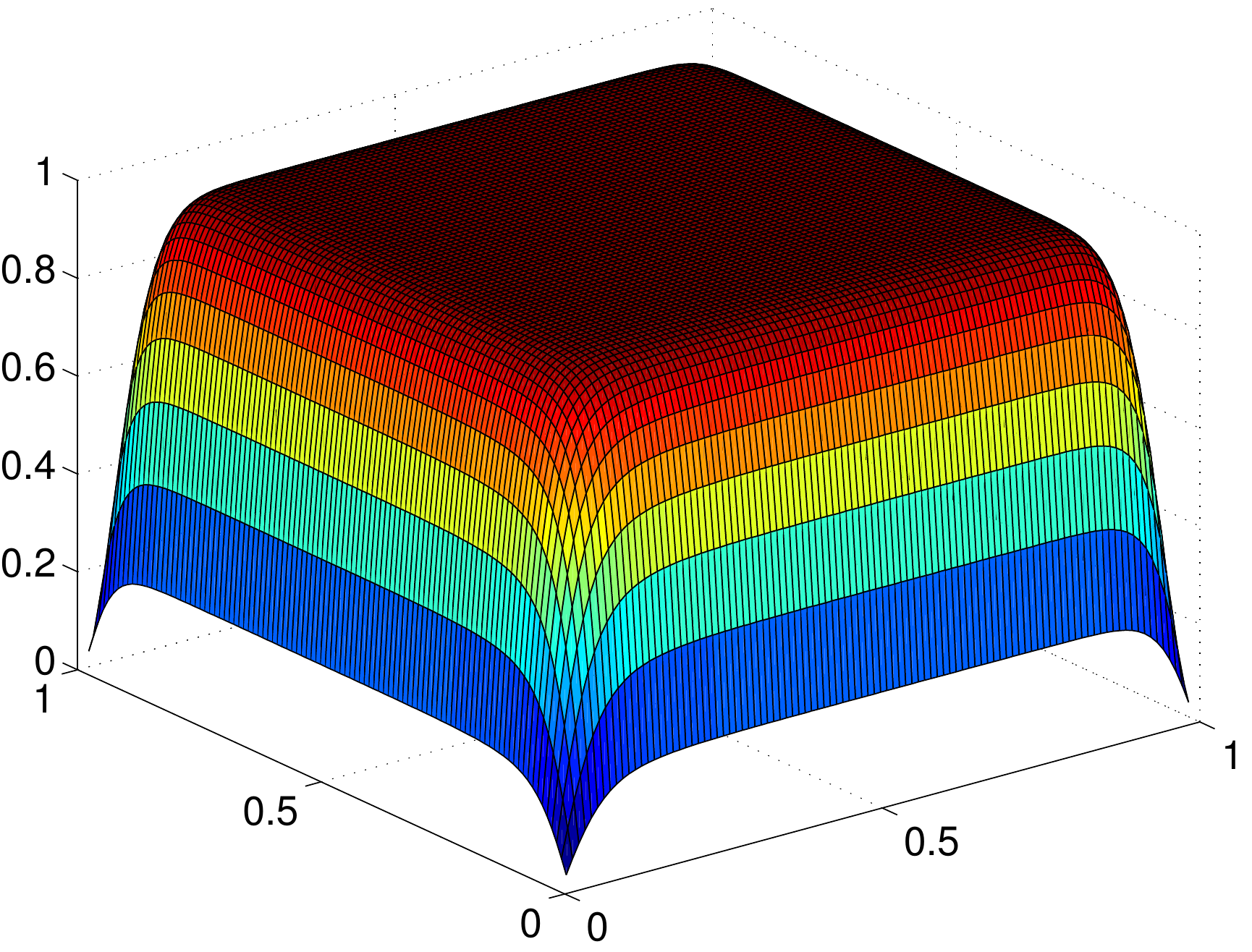}
\includegraphics[scale=0.25]{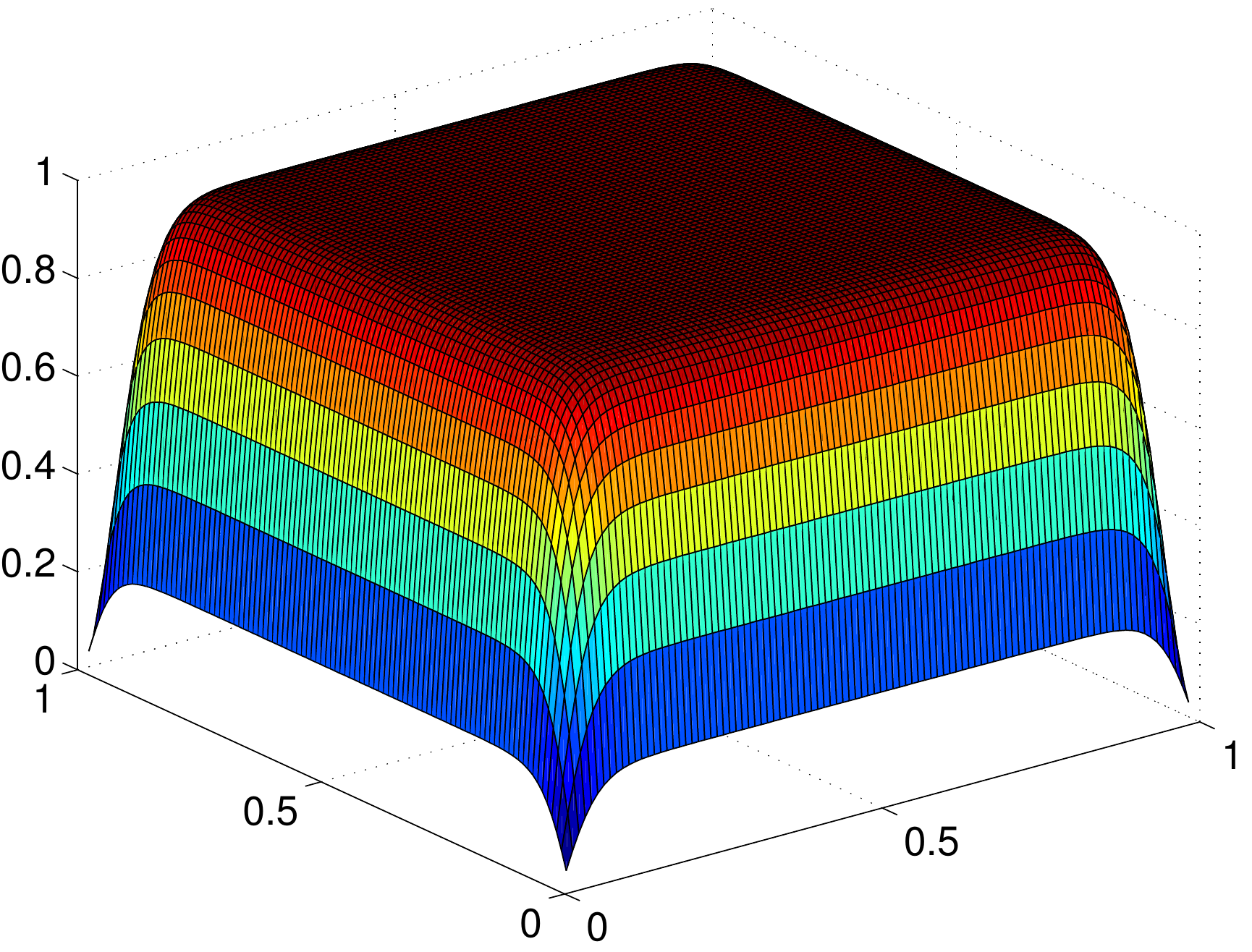}\\
\end{center}
\caption{Test 3: Solution of equation \eqref{prb_test1} at time $t=\{0,0.1\}$ (top) and $t=\{1.5,3\}$ (bottom)}.
\label{test1:sol}
\end{figure}

The complexity of problem \eqref{prb_test1} is reduced by model order reduction. When dealing with model order reduction, it is relevant to consider the CPU time of the simulation and the error. In general it is important to have a trade-off between the two quantities. Figure \ref{test1:an} considers the CPU time.  The POD-DMD approximation is always faster than any other approximation for any dimension of the reduced system. Increasing the number of POD basis functions, the POD-DEIM turns out to be even more expensive than POD. The strength of the POD-DMD is the fact we do not have to evaluate the nonlinearity after collecting snapshots. We note that the number of POD, DEIM and DMD are always the same in Figure \ref{test1:an}. On the right of Figure \ref{test1:an} we compute the relative error with respect to the Frobenius norm. It is clear that the POD provides the best approximation. As expected, all the methods decrease their error when increasing the number of basis functions.

\begin{figure}[htbp]
\begin{center}
\includegraphics[scale=0.25]{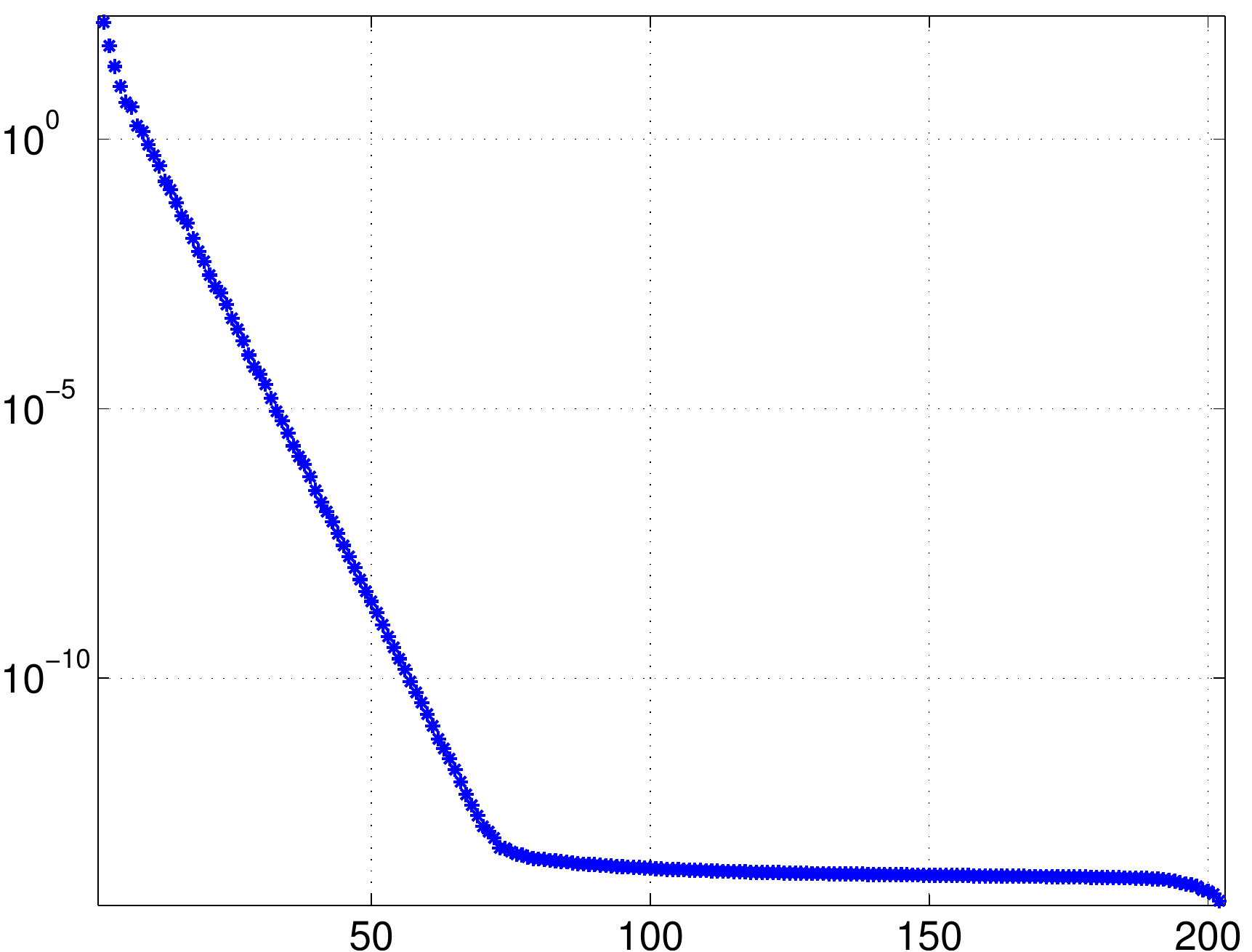}
\includegraphics[scale=0.25]{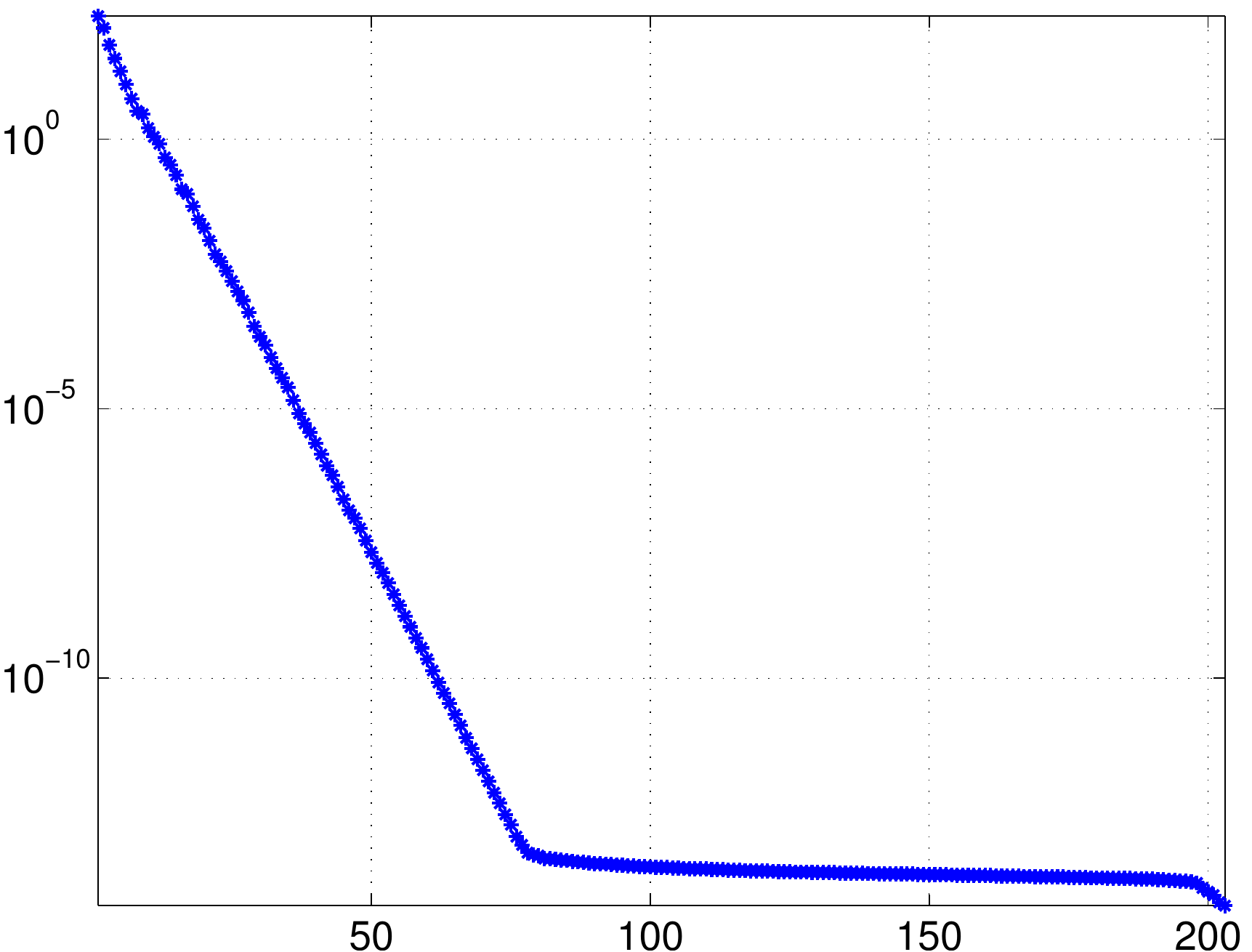}\\
\end{center}
\caption{Test 3: Singular values of the solution (left) and the nonlinearity (right) of \eqref{prb_test1}}
\label{test1:svd}
\end{figure}

\begin{figure}[htbp]
\begin{center}
\includegraphics[scale=0.25]{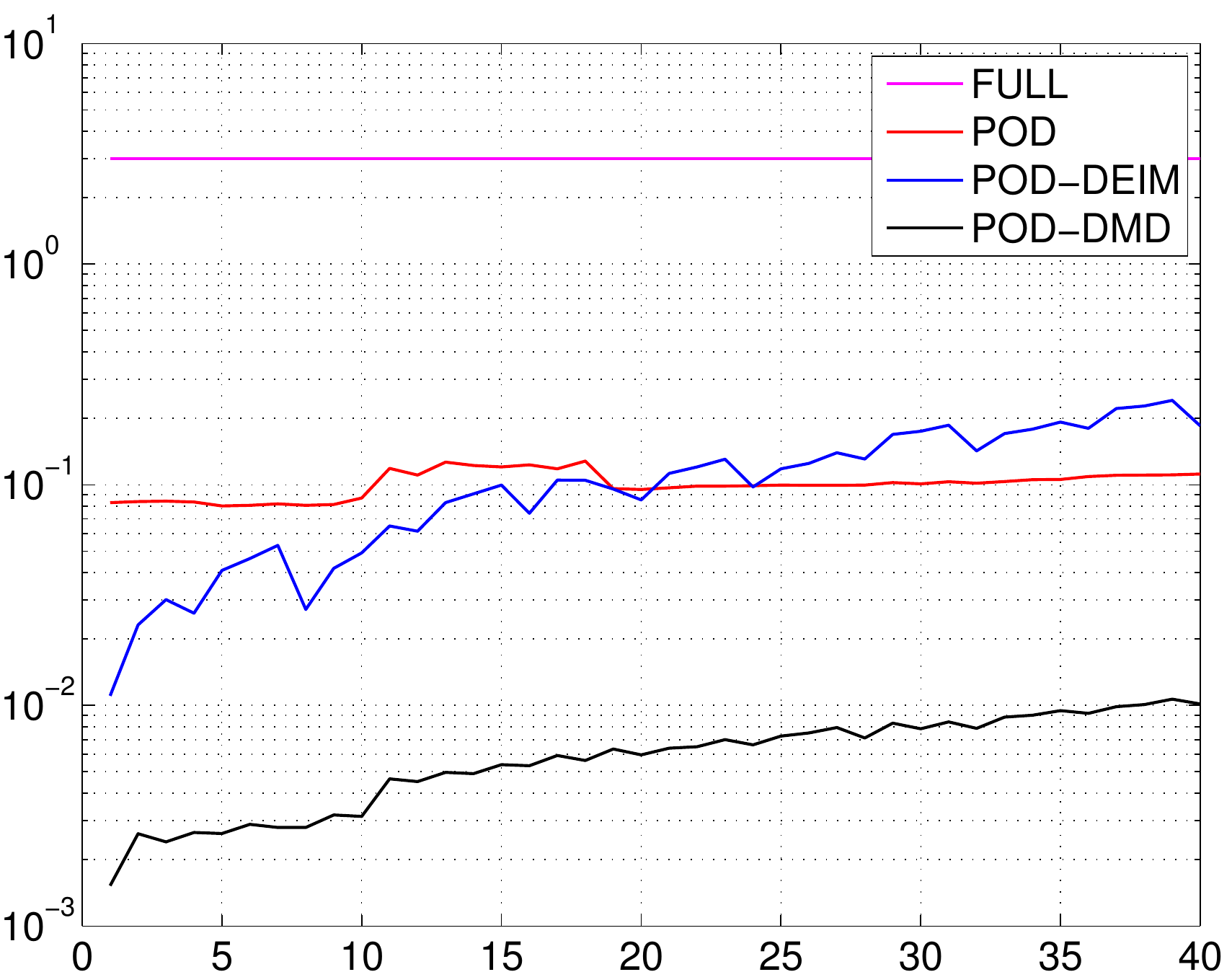}
\includegraphics[scale=0.25]{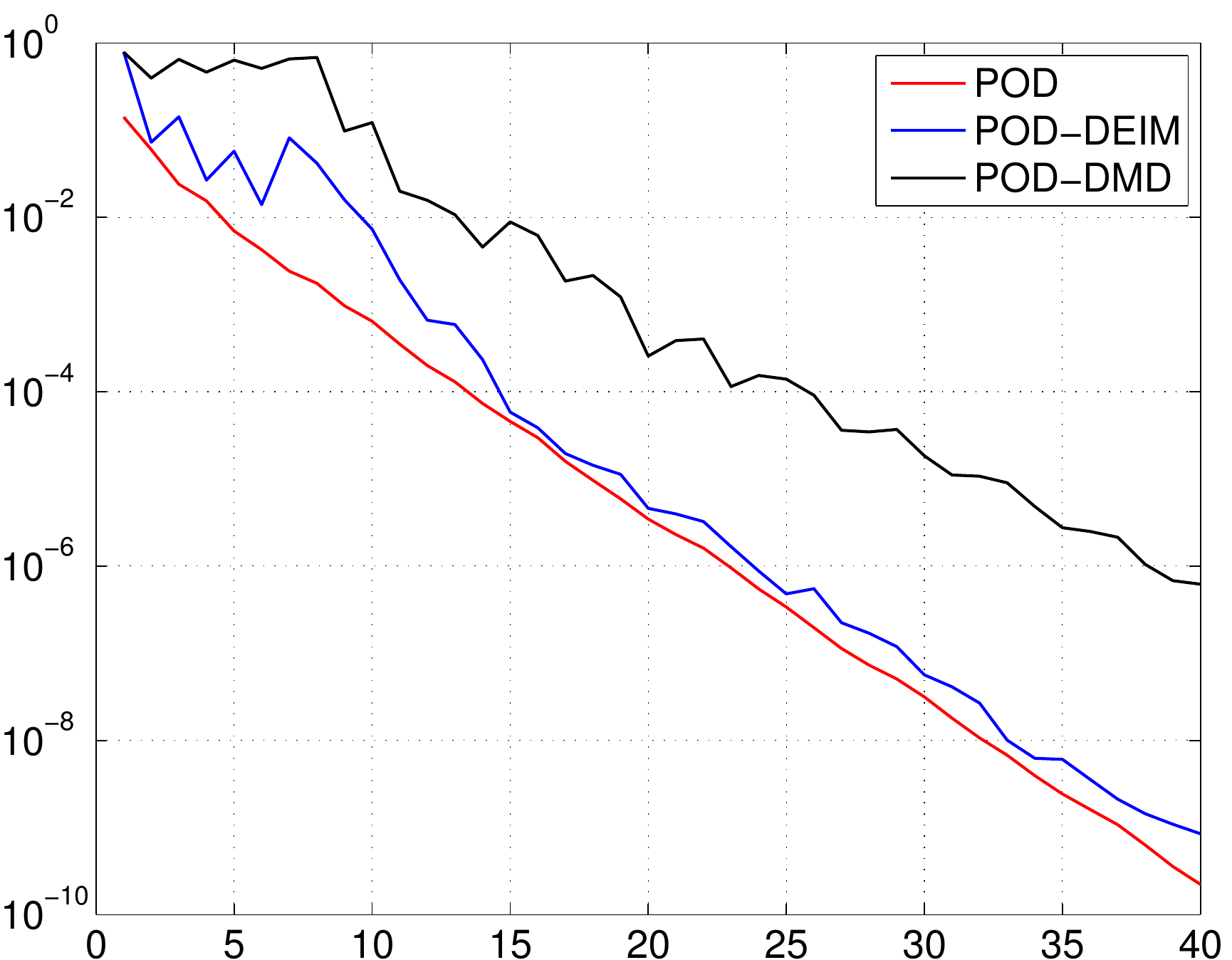}\\
\end{center}
\caption{Test 3: CPU-time (left) and Relative Error in Frobenius norm. Number of POD modes and DEIM/DMD points are the same}
\label{test1:an}
\end{figure}

Since the POD-DMD is faster than other method it is natural to look at the performance with different number of basis functions. Figure \ref{test1:an2} shows the error for a fixed number of basis functions $\ell=\{5,10,15\}$ and $k\in[1,40]$. As we can see the POD-DMD performs better than POD-DEIM when $\ell=5,10$. In this case increasing the number of DMD basis functions lead to more accurate solutions of POD-DEIM. Moreover, we can observe a monotone decay of the error for the POD-DMD approach.

\begin{figure}[htbp]
\begin{center}
\includegraphics[scale=0.22]{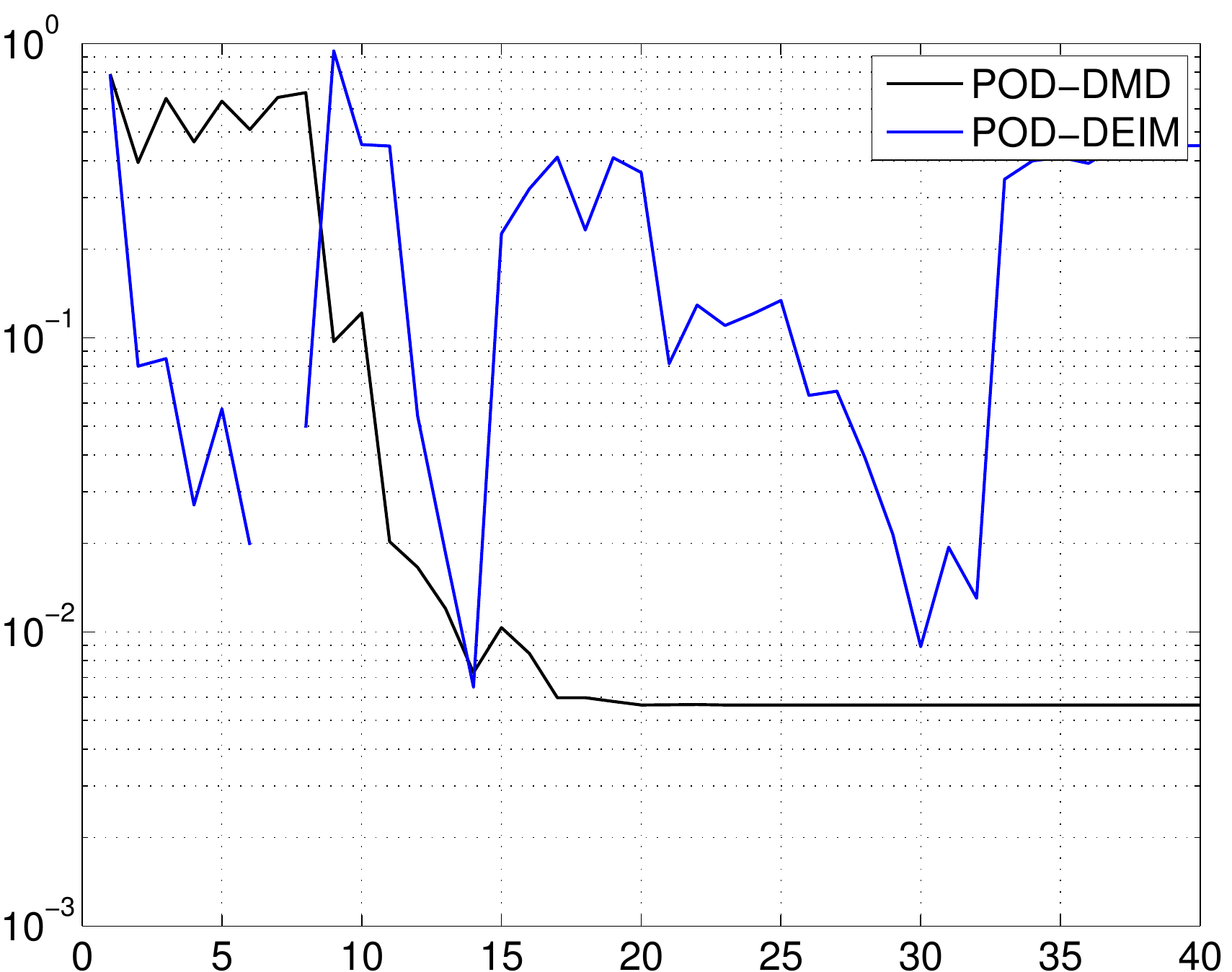}
\includegraphics[scale=0.22]{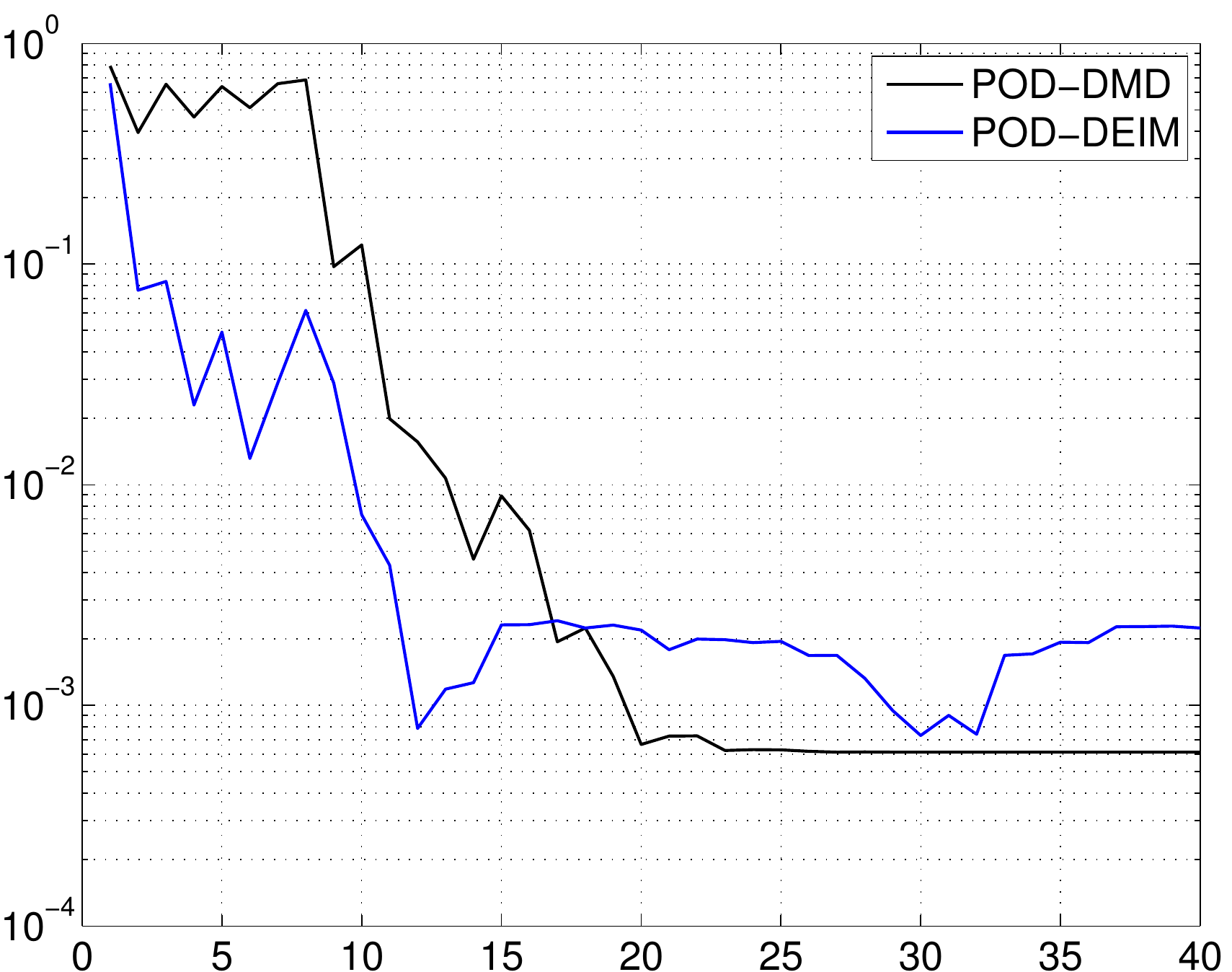}
\includegraphics[scale=0.22]{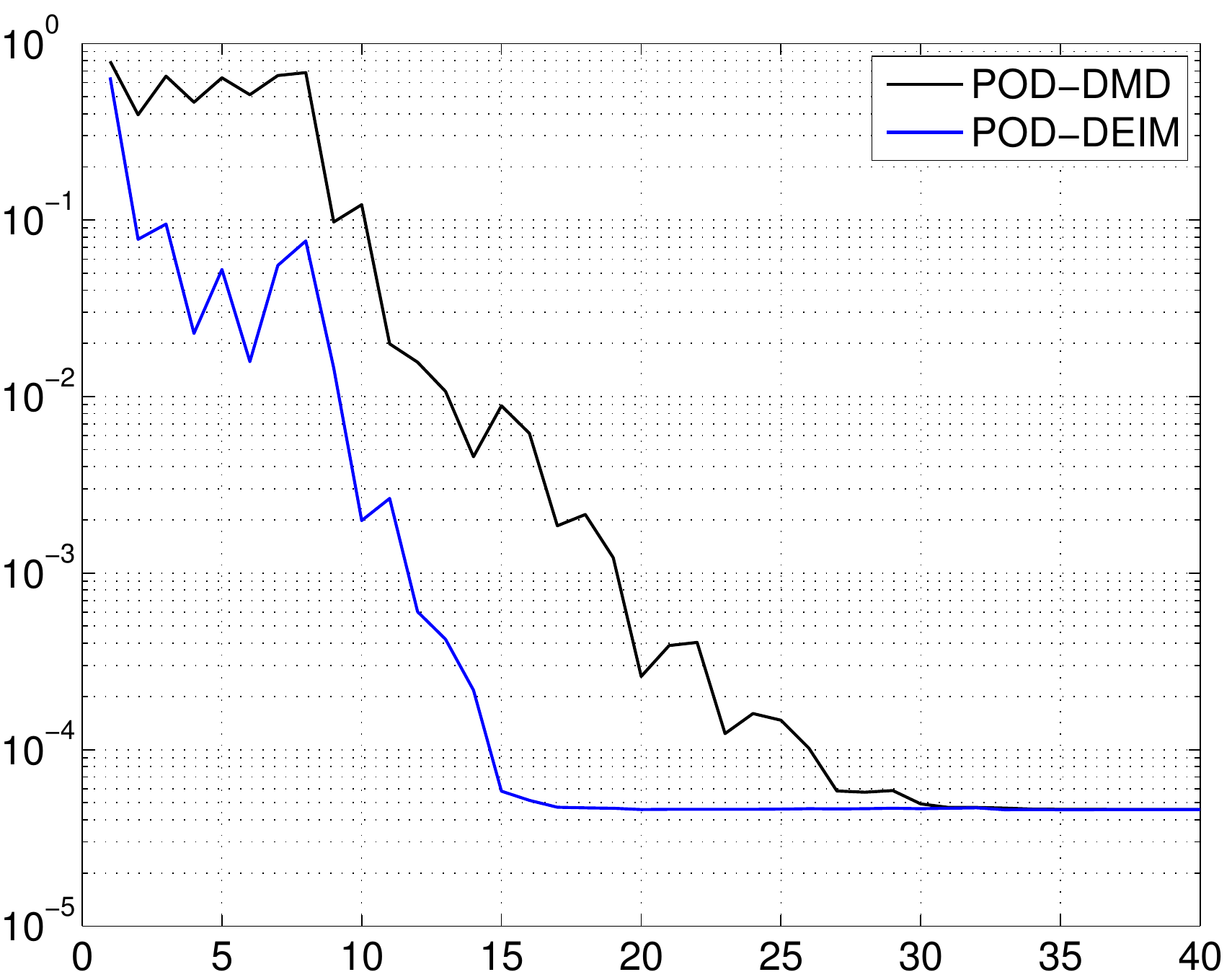}
\end{center}
\caption{Test 3: Relative Error for 5 POD basis functions (left), 10 POD basis (middle), 15 POD basis (right)}

\label{test1:an2}
\end{figure}

\paragraph{\bf Test 4: Burgers' Equation}
Let us consider the following 1D Burger's equation:
\begin{equation}\label{prb_test2}
\begin{aligned}
y_t(x,t)-\theta y_{xx}(x,t)+y(x,t)y_x(x,t)&=0&& (x,t)\in[a,b]\times[0,T],\\
y(x,0)&=y_0(x)&& x\in [a,b],\\
y(a,t)&=0=y(b,t)&&  t\in [0,T],
\end{aligned}
\end{equation}
where $a=0, b=1,T=1, \theta=0.01, y_0(x)=\mbox{sgn}(x)$ 
In Figure \ref{test2:svd} we visualize the full solution of \eqref{prb_test2} and the decay of the singular values. We note that the decay is very similar for $y$ and its nonlinearity.
The results of the model reduction are shown in Figure \ref{test2:sol}. In the left panel, the POD approximation with 20 basis function is demonstrated.  The POD-DEIM is visualized in the middle panel, and finally, the POD-DMD approach is in the right panel.  With a 20-rank truncation of the DEIM and DMD approximation,  it is difficult to see any differences in the solution.

\begin{figure}[htbp]
\begin{center}
\includegraphics[scale=0.22]{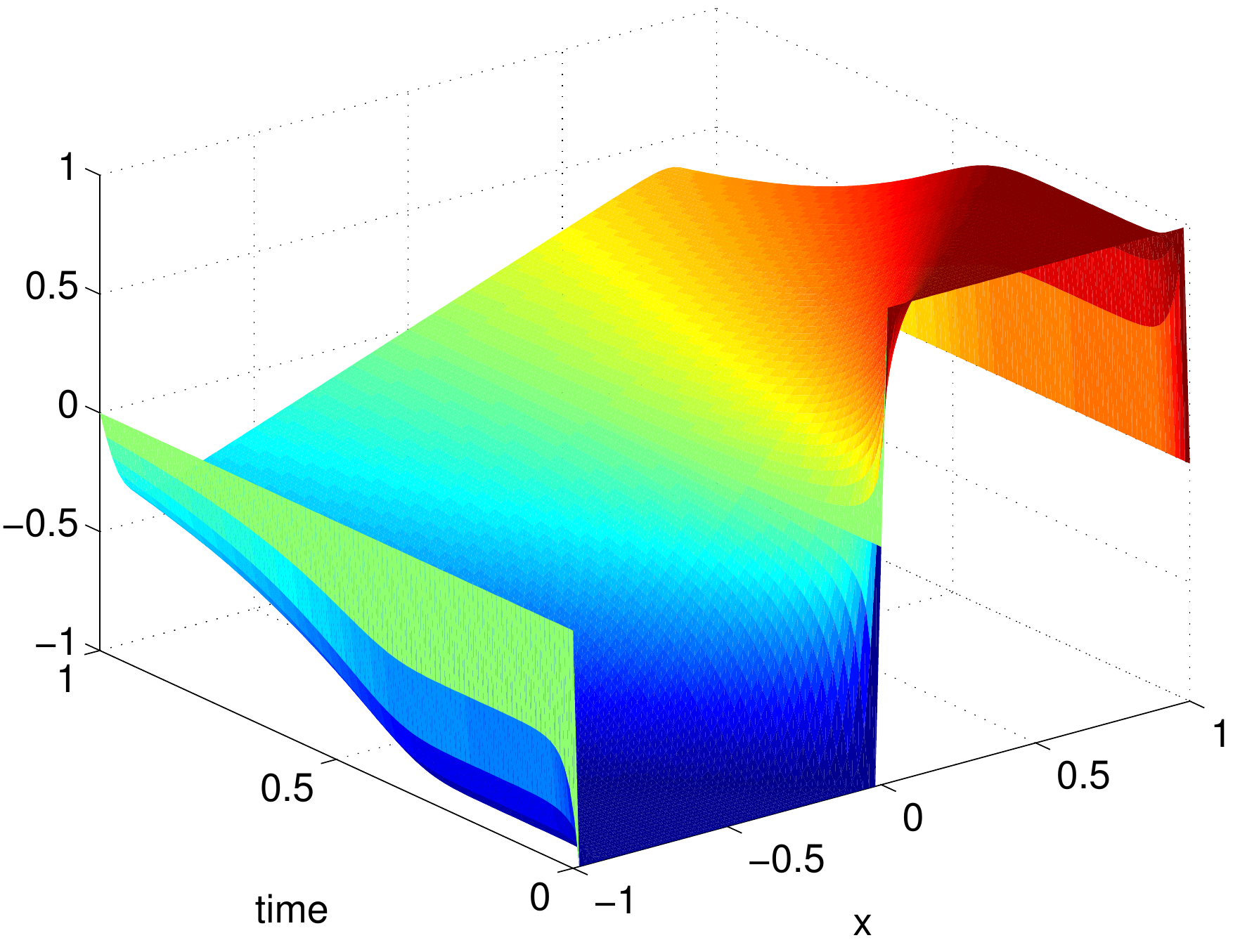}
\includegraphics[scale=0.22]{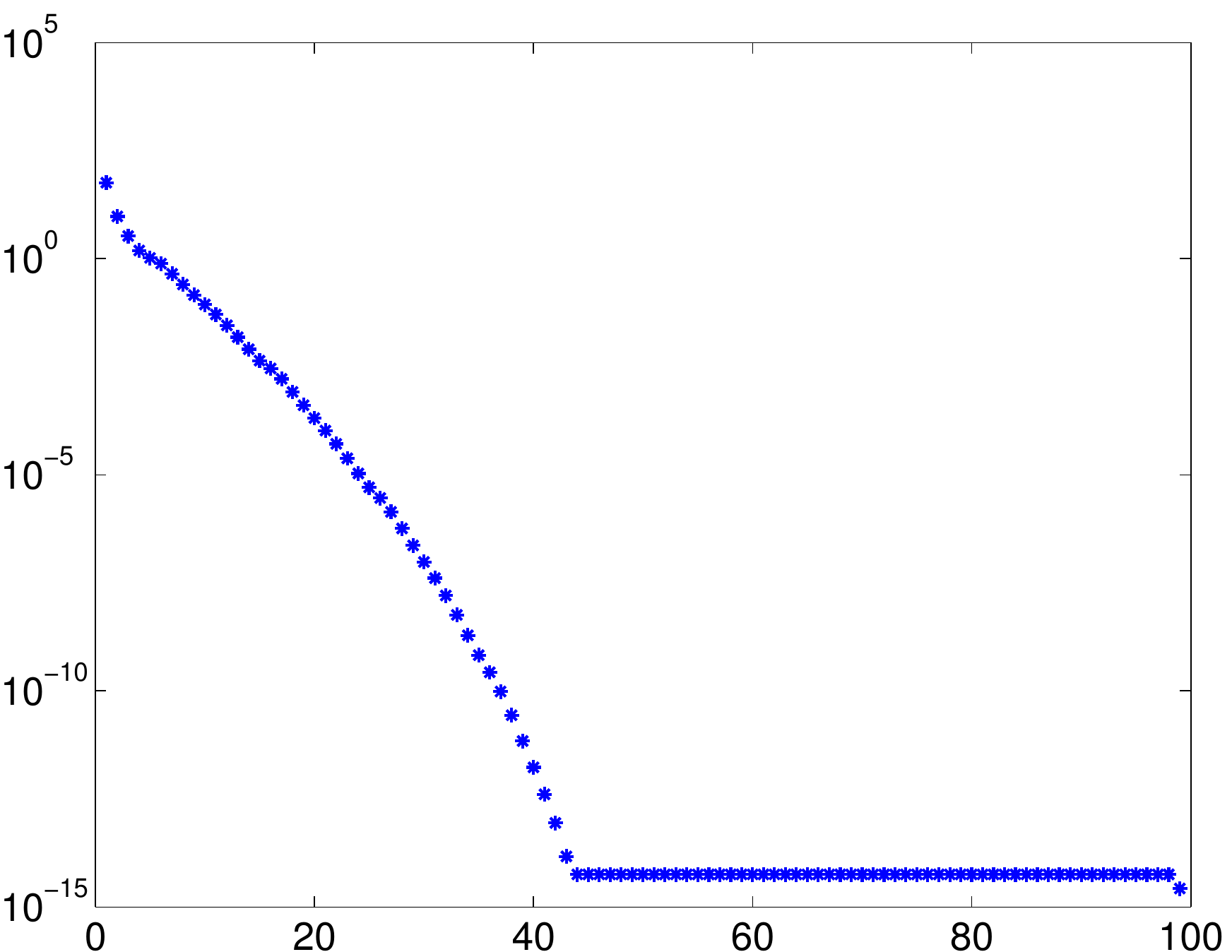}
\includegraphics[scale=0.22]{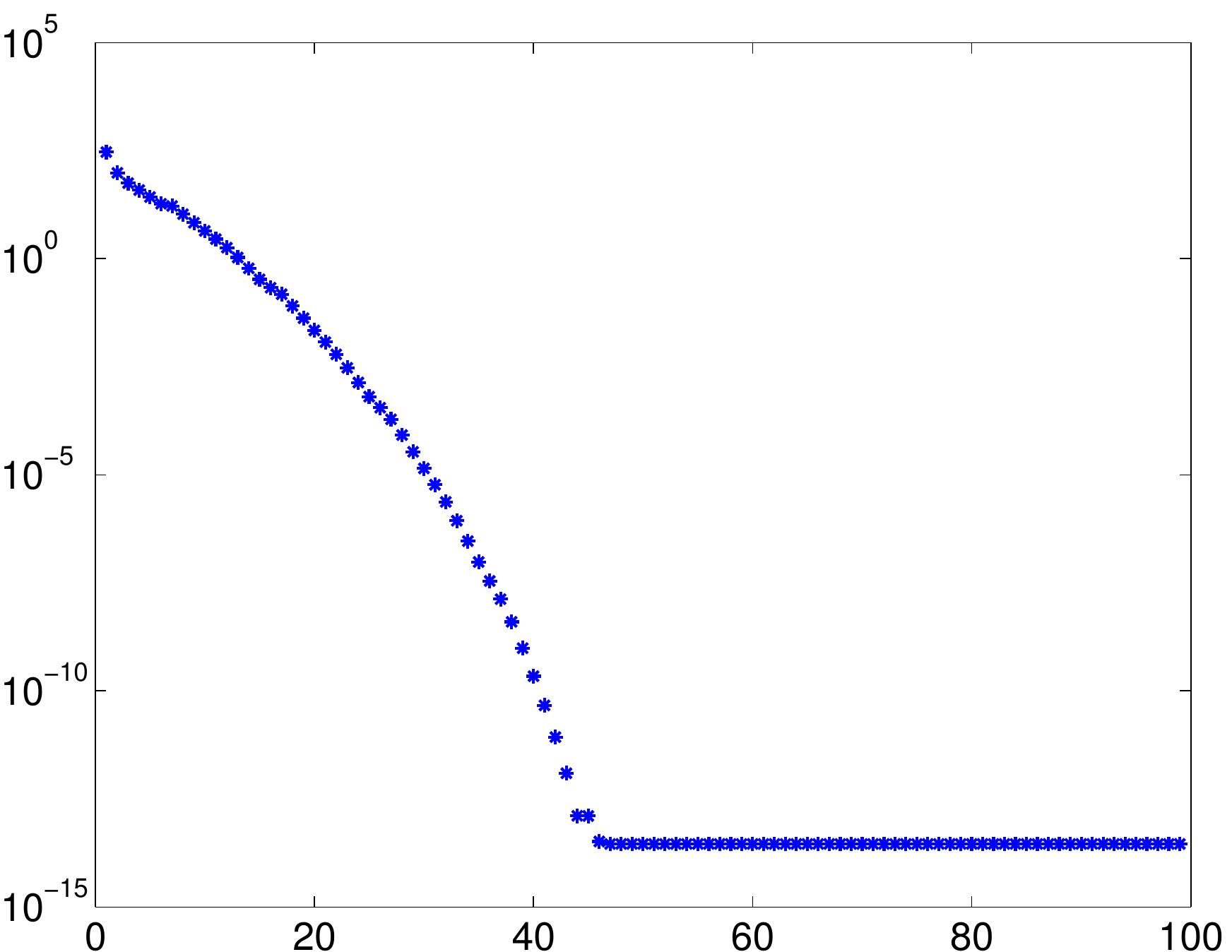}
\end{center}
\caption{Test 4: Full approximation (left), Singular values of the solution (middle) and of the nonlinearity (right) of \eqref{prb_test2}}
\label{test2:svd}
\end{figure}

\begin{figure}[htbp]
\begin{center}

\includegraphics[scale=0.22]{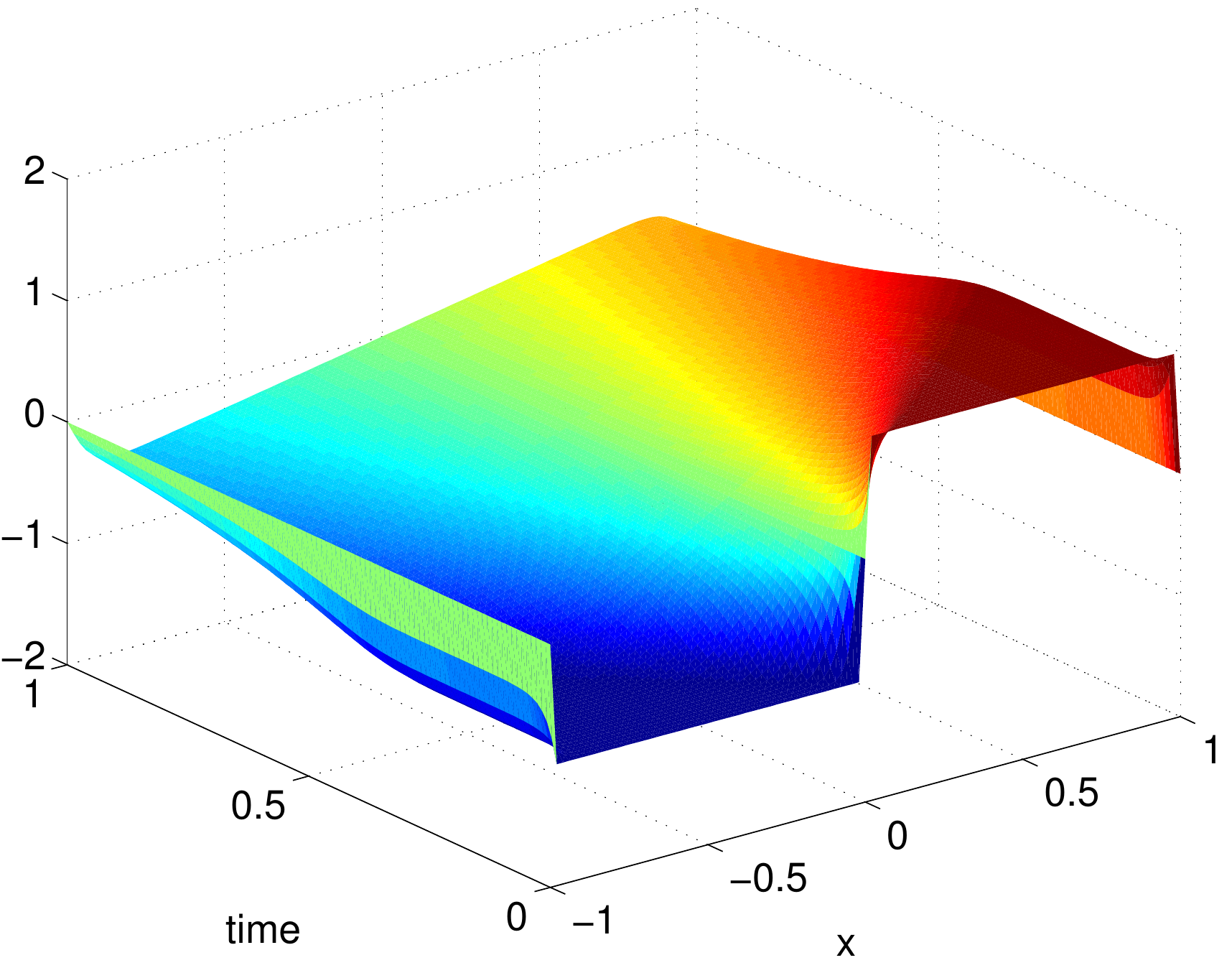}
\includegraphics[scale=0.22]{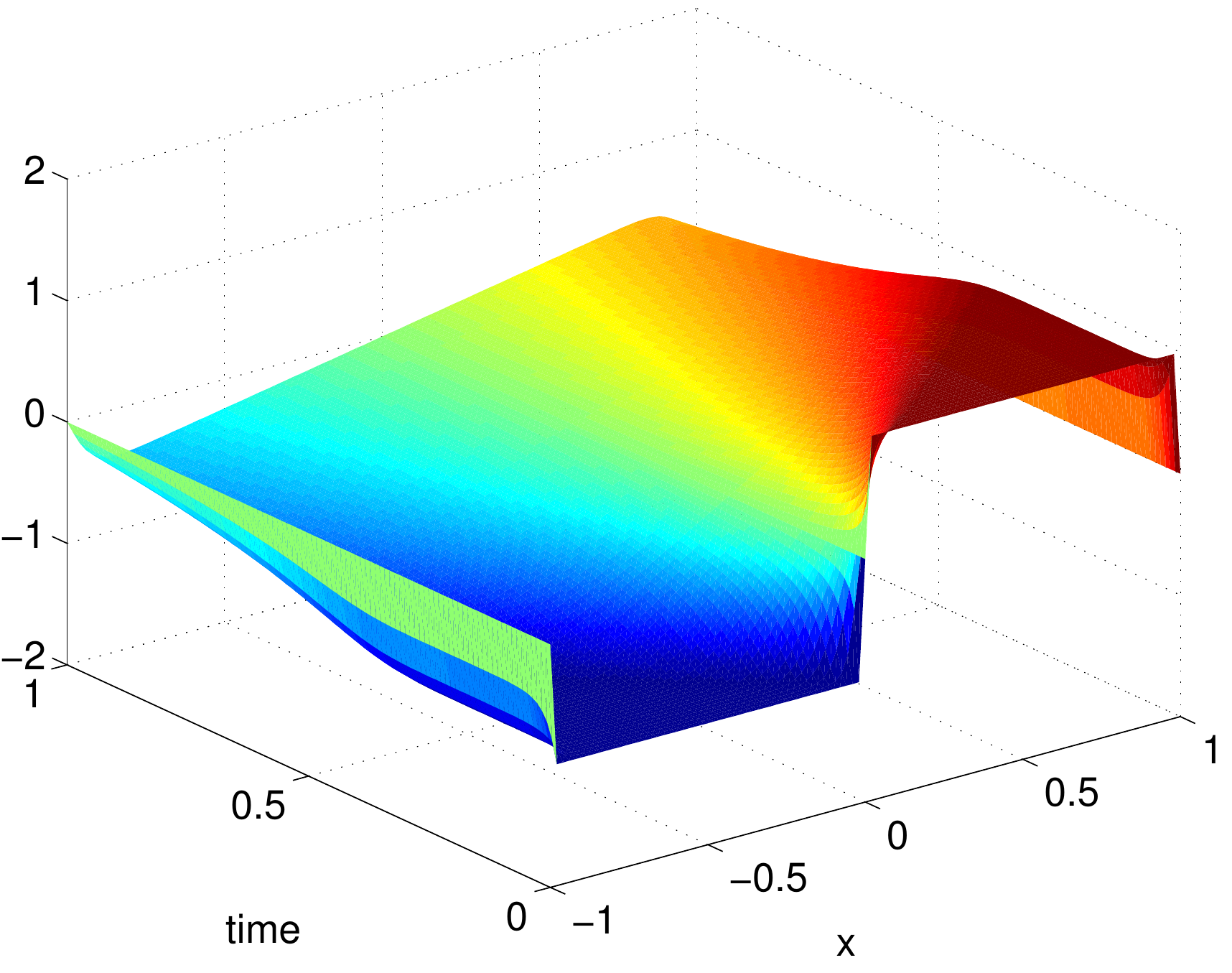}
\includegraphics[scale=0.22]{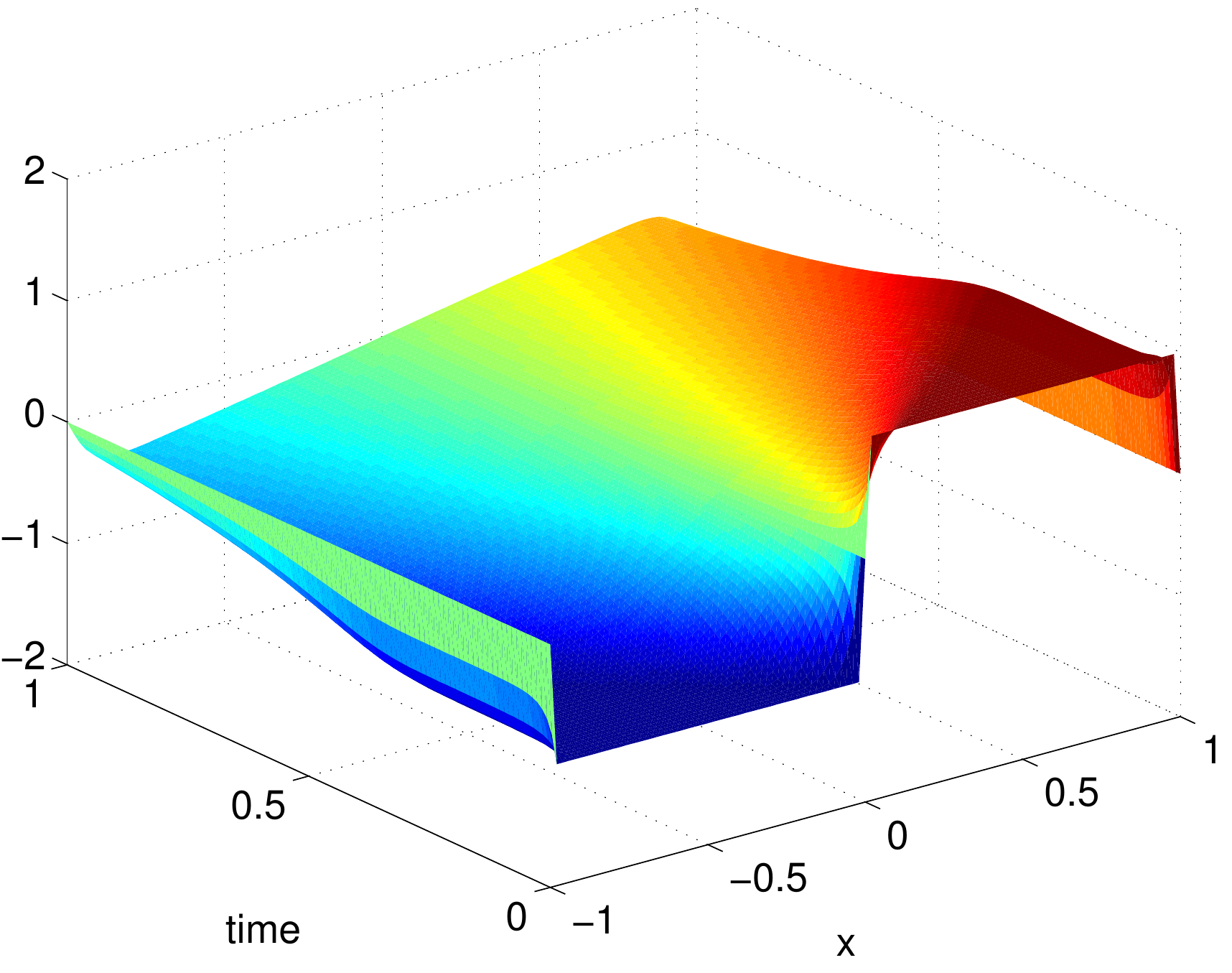}
\end{center}
\caption{Test 4: Approximation with 20 POD basis (left), Approximation with 20 POD basis function and 20 DEIM (middle), Approximation with 20 POD basis function and 20 DMD (right)}
\label{test2:sol}
\end{figure}

The CPU time is expressed in Figure \ref{test2:err} and we can see, as already discussed in the previous example, that the CPU time of the POD-DMD is always below the other two approximations but we lose some accuracy, as expected. The relative error is in Figure \ref{test2:err} we can see that the POD-DMD does not decrease even if we increase the number of POD basis functions. Let us remember that the Burger's equation has a hyperbolic structure which is more complicated to capture, especially in the nonlinear term.

\begin{figure}[htbp]
\begin{center}
\includegraphics[scale=0.25]{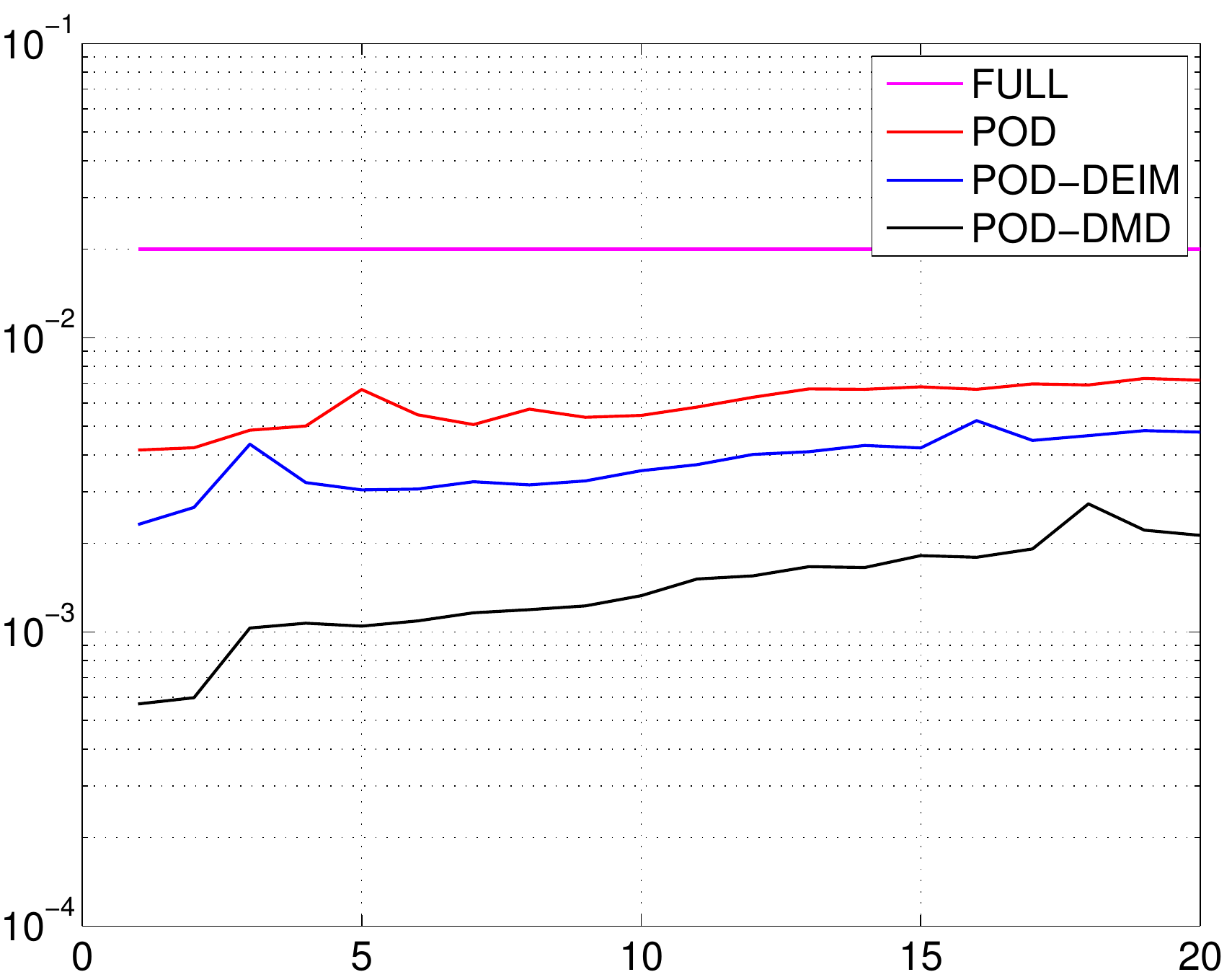}
\includegraphics[scale=0.25]{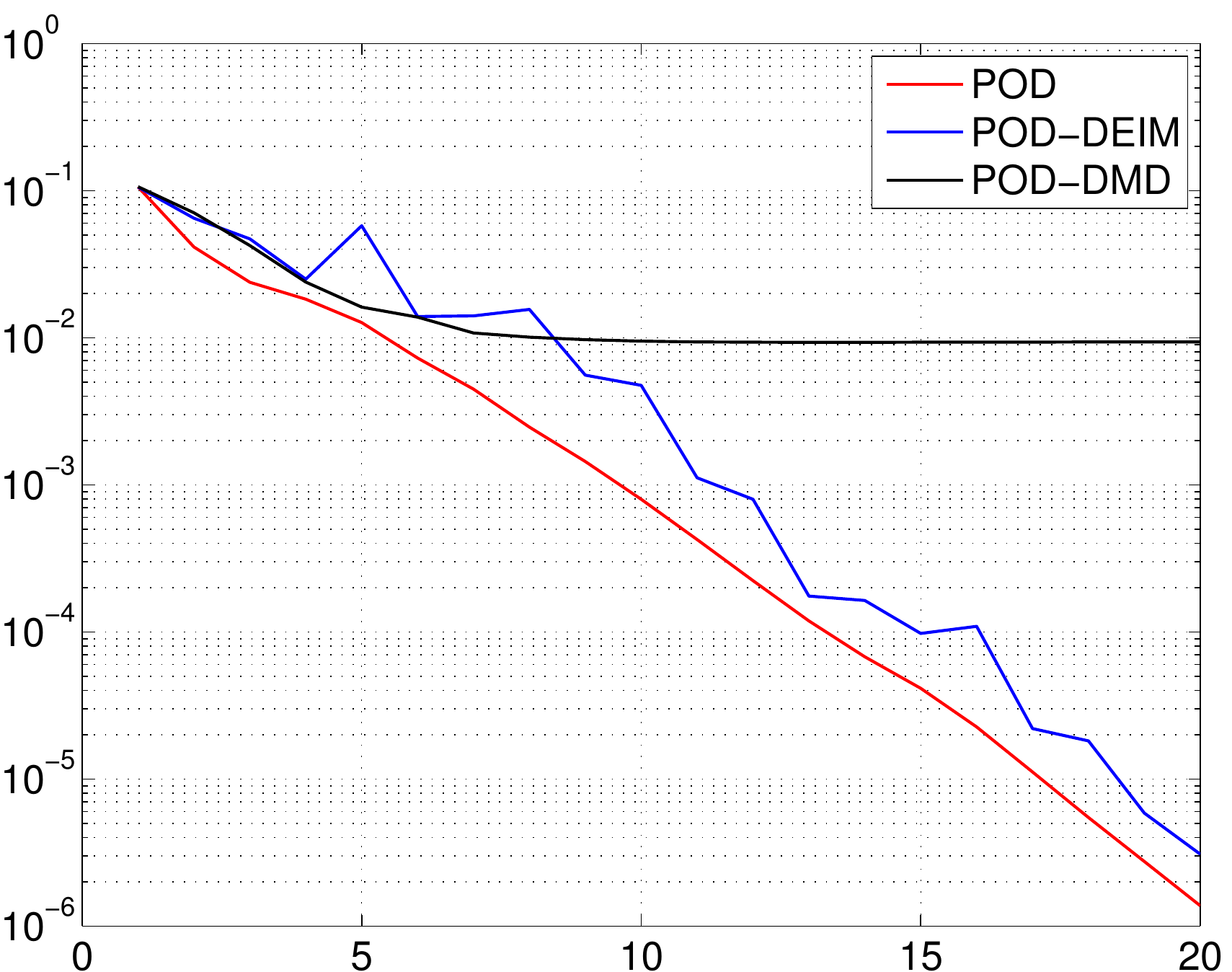}\\
\end{center}
\caption{Test 4: CPU-time (left) and Relative Error in Frobenius norm. Number of POD modes and DEIM/DMD points are the same}
\label{test2:err}
\end{figure}

It is difficult to directly compare POD-DMD and POD-DEIM since the meaning of the rank in DEIM is different from DMD. For this reason we also computed the error varying the rank $k$ for a fixed number of POD basis functions. Figure \ref{test2:err2} shows that POD-DMD is more stable when the number of POD basis functions is not large (left picture). The POD-DEIM is in general more accurate, especially for large $k$. Since the POD-DMD is always faster, this is not a big issue, and in fact, one could work with a low-dimensional structure of the POD basis functions and consider a larger number of DMD basis functions.

\begin{figure}[htbp]
\begin{center}
\includegraphics[scale=0.22]{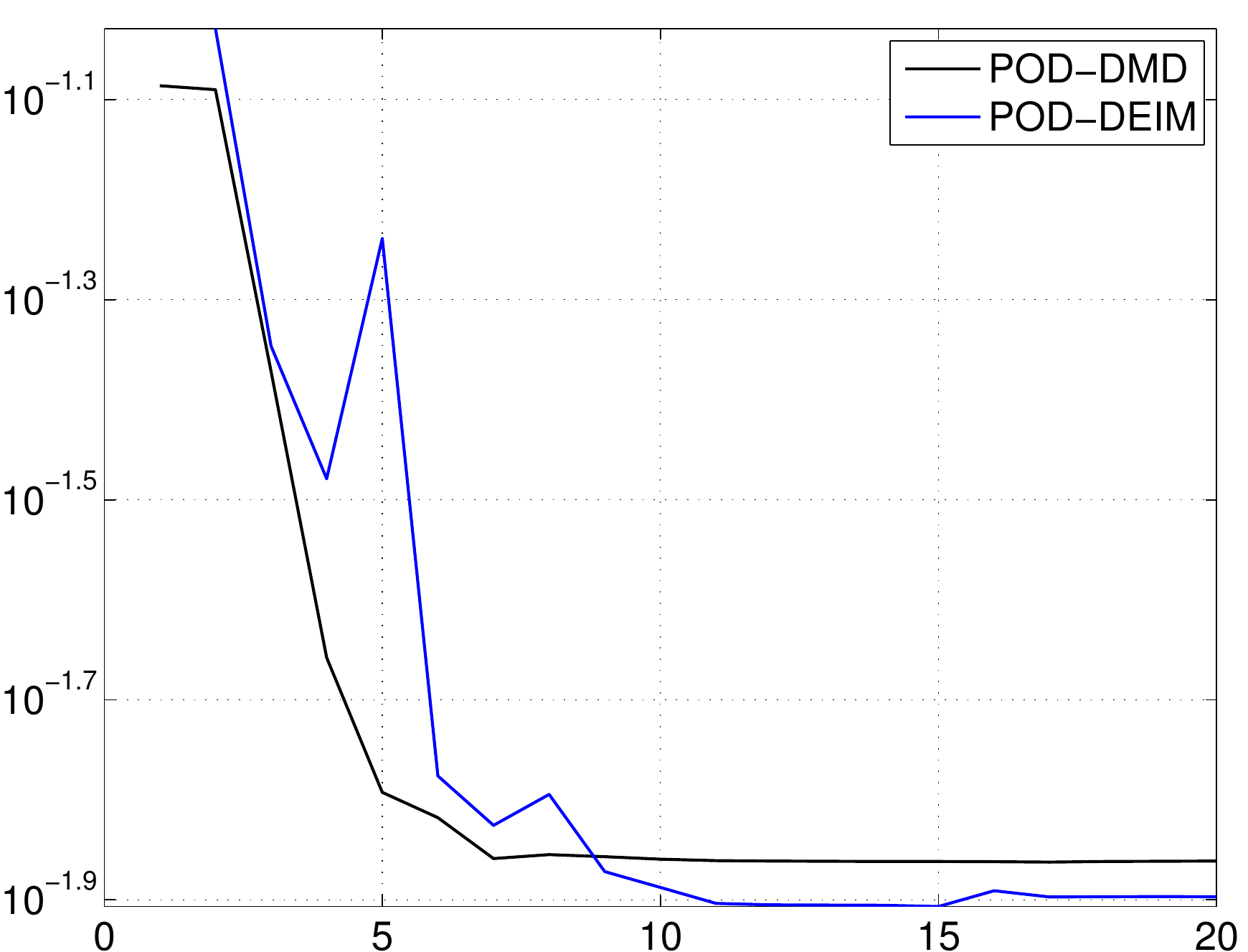}
\includegraphics[scale=0.22]{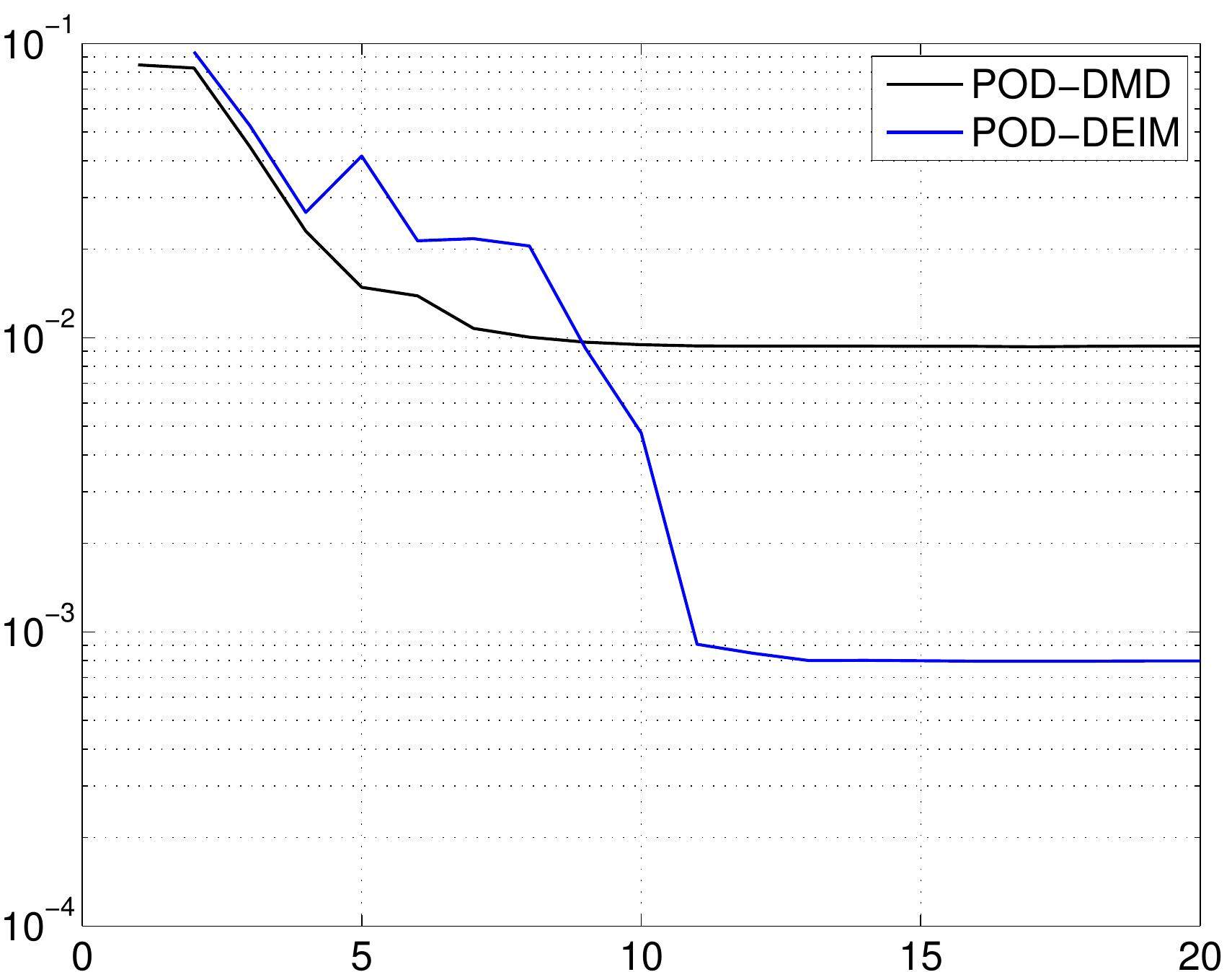}
\includegraphics[scale=0.22]{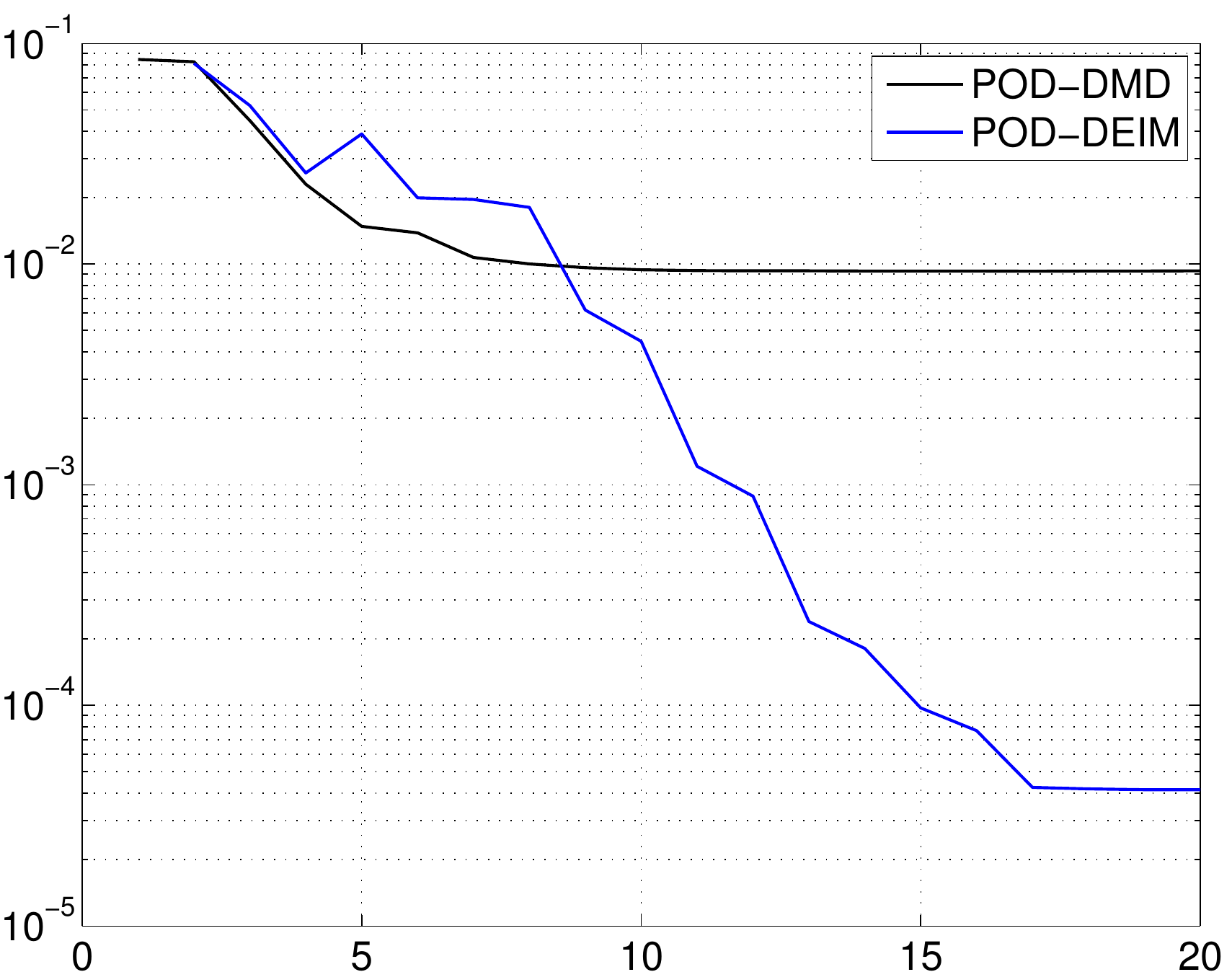}
\end{center}
\caption{Test 4: Relative Error for 5 POD basis functions (left), 10 POD basis (middle), 15 POD basis (right)}
\label{test2:err2}
\end{figure}


\paragraph{\bf Test 5: Nonlinear Schr\"odinger equation}
Let us consider the following Schr\"odinger equation
\begin{equation}\label{prb_test3}
\begin{aligned}
y_t-i \theta y_{xx}-i |y|^2y&=0&& (x,t)\in[-L,L]\times[0,T],\\
y(x,0)&=y_0(x)&& x\in [a,b],\\
y(-L,t)&=0=y(L,t)&& t\in [0,T],
\end{aligned}
\end{equation}
where $L=15, T=2, \theta=1$ and $y_0(x)=sech(x)$. The solution of \eqref{prb_test3} is shown in Figure \ref{test3:svd} on the left. The singular values of the solution are shown in the middle panel while the singular values of the nonlinearity are in the right panel. It is well-known that Schr\"odinger's equation generates waves functions in its solution and therefore it is difficult to capture this behavior with  only a few modes.

\begin{figure}[htbp]
\begin{center}
\includegraphics[scale=0.22]{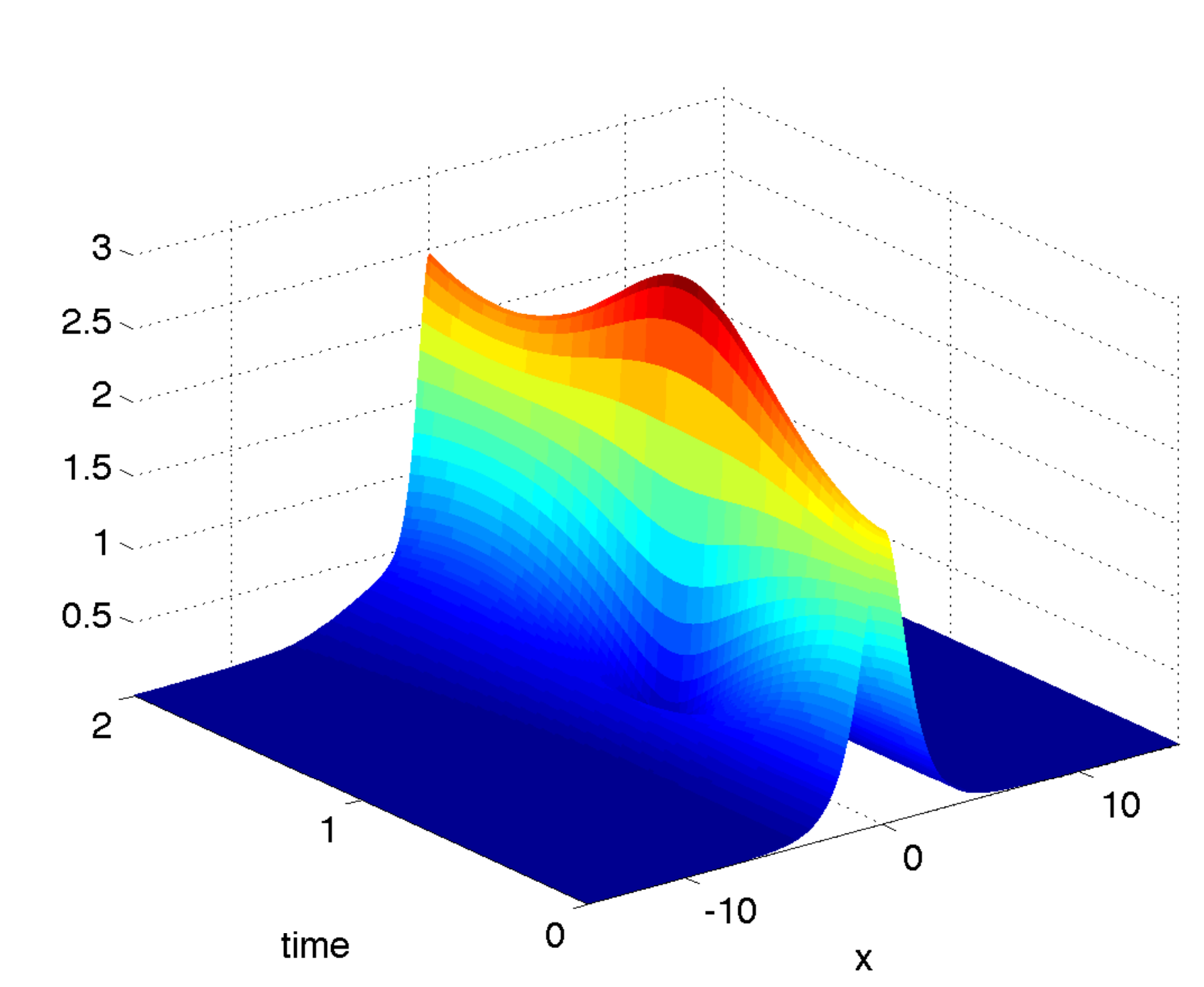}
\includegraphics[scale=0.22]{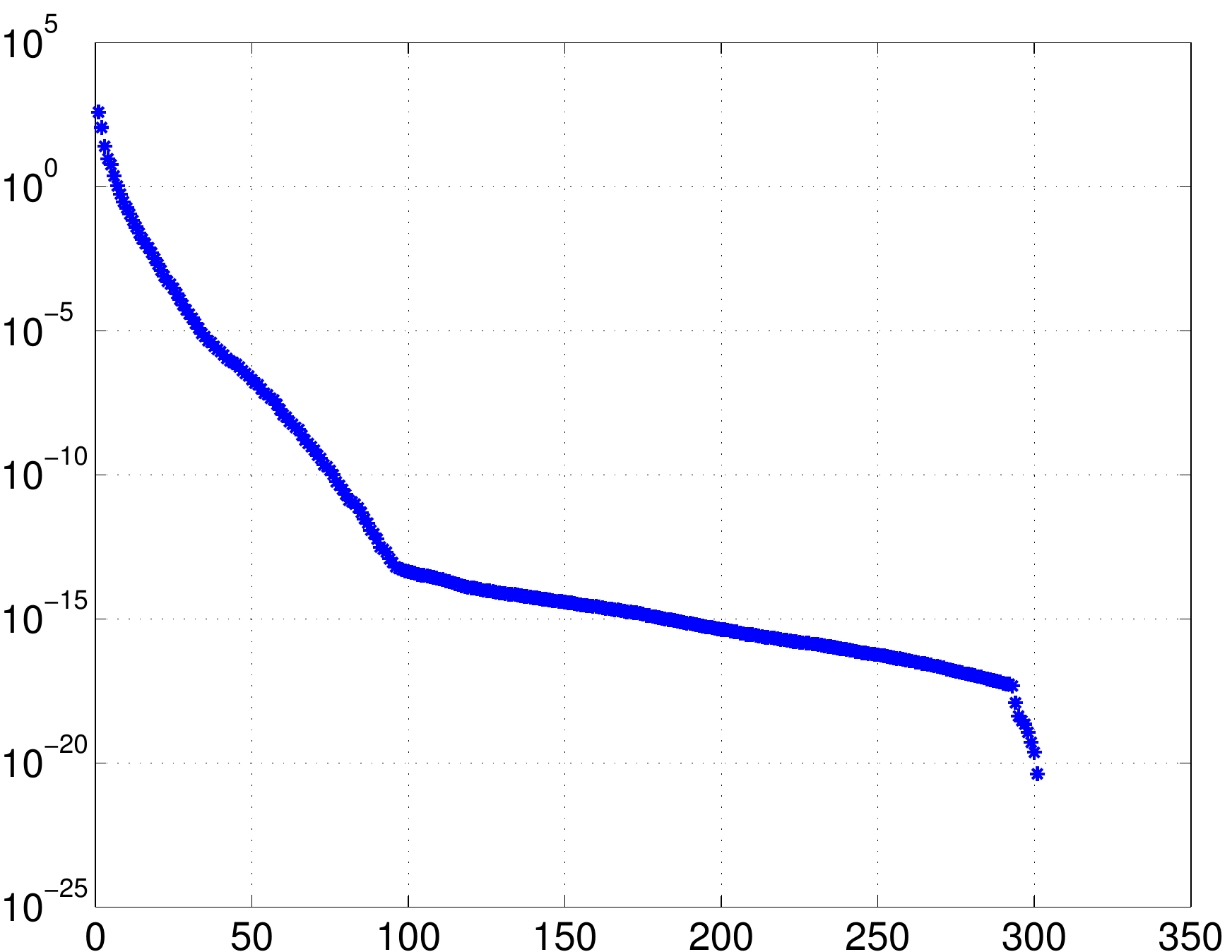}
\includegraphics[scale=0.22]{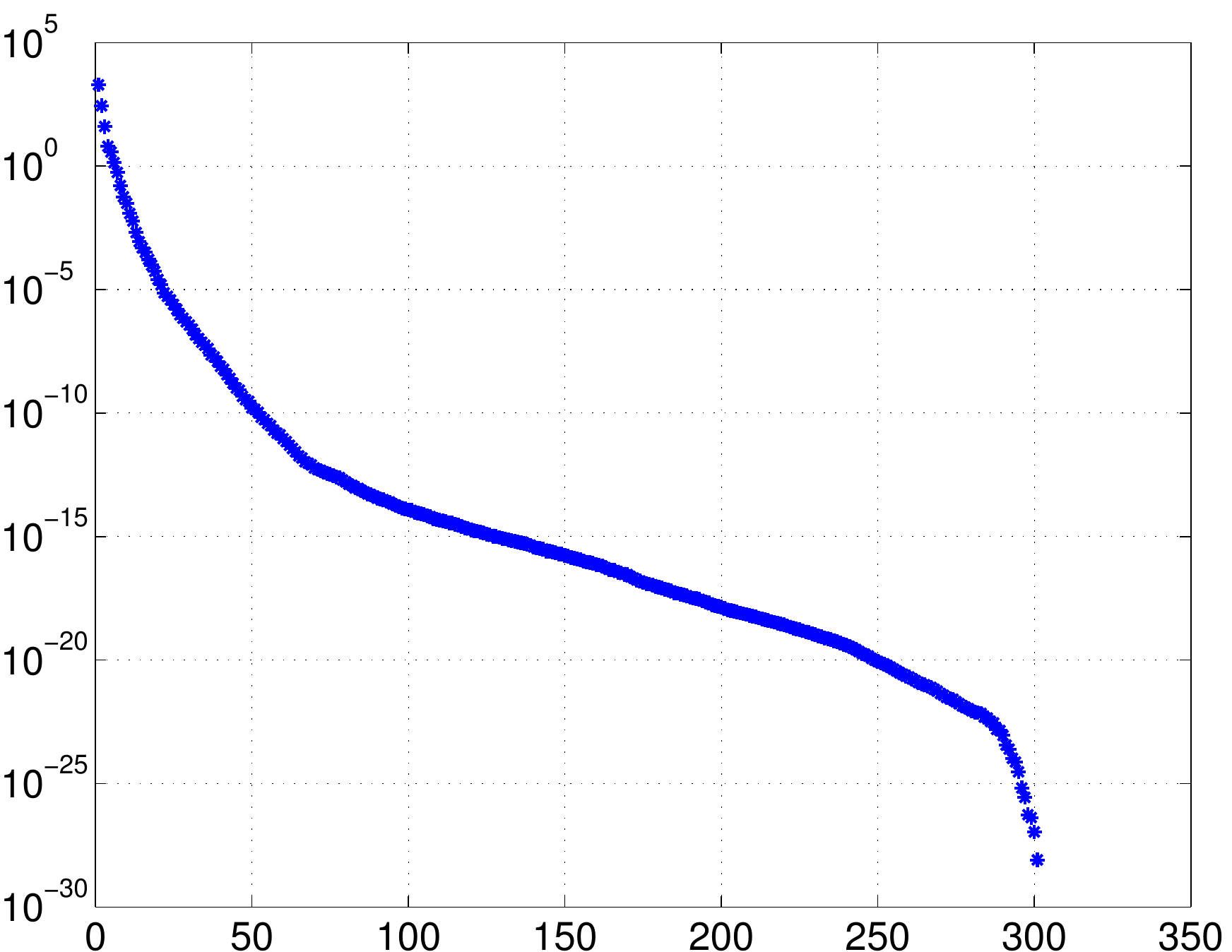}
\end{center}
\caption{Test 5: Full approximation (left), Singular values of the solution (middle) and of the nonlinearity (right) of \eqref{prb_test3}}
\label{test3:svd}
\begin{center}
\includegraphics[scale=0.22]{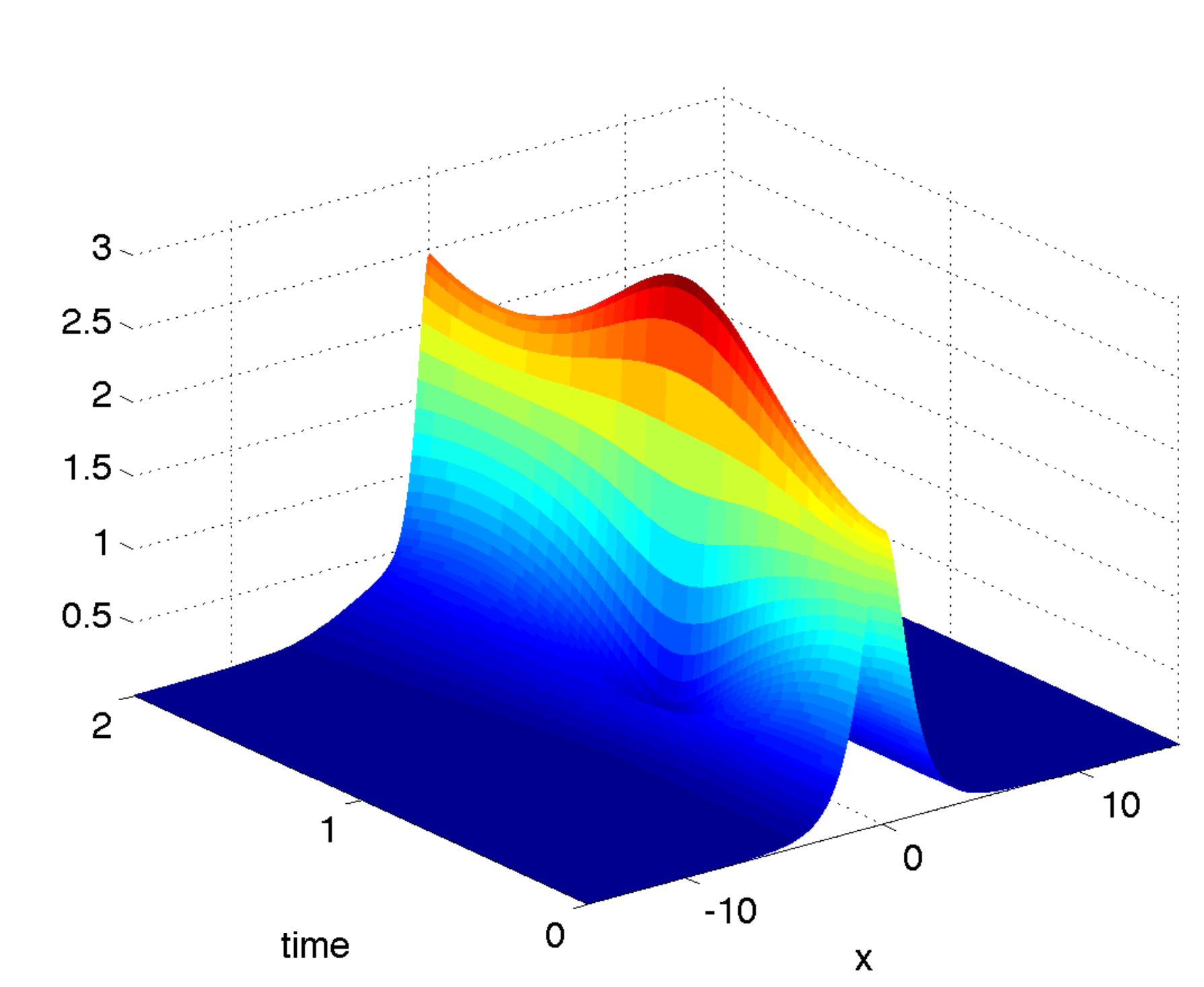}
\includegraphics[scale=0.22]{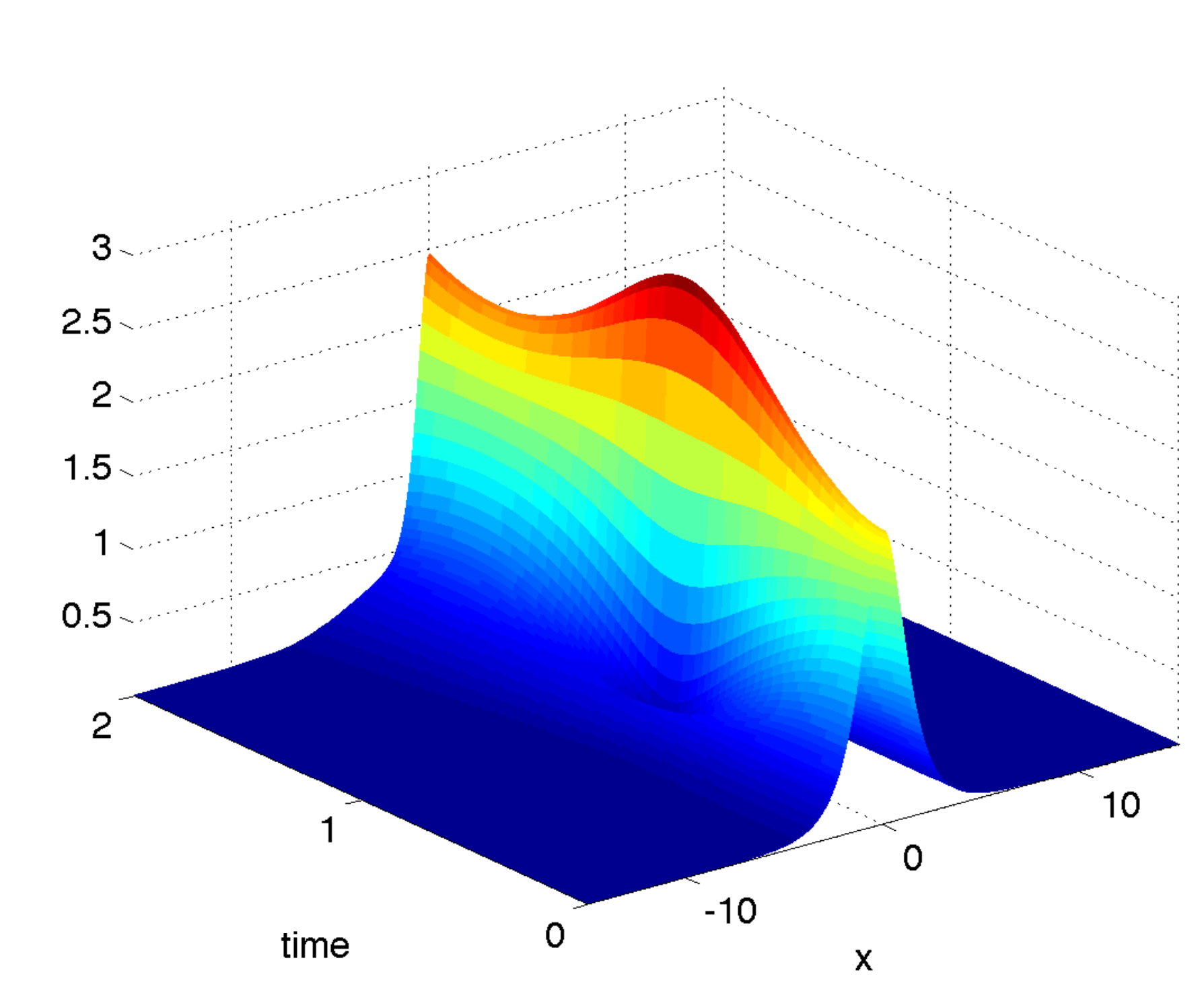}
\includegraphics[scale=0.22]{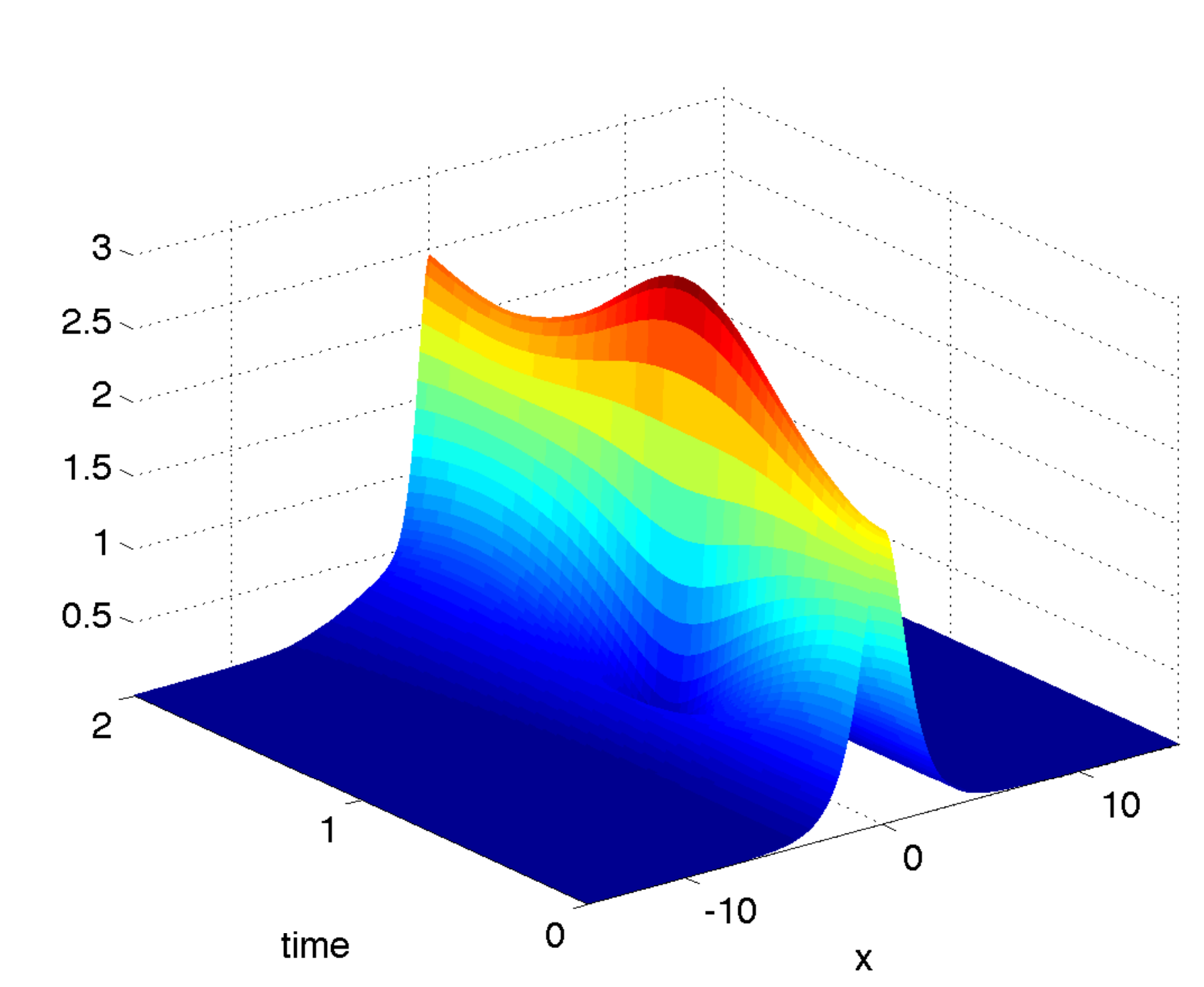}
\end{center}
\caption{Test 5: Approximation with 10 POD basis (left), Approximation with 10 POD basis function and 10 DEIM (middle), Approximation with 10 POD basis function and 10 DMD (right)}
\label{test3:sol}
\end{figure}

%

\begin{figure}[H]
\begin{center}
\includegraphics[scale=0.25]{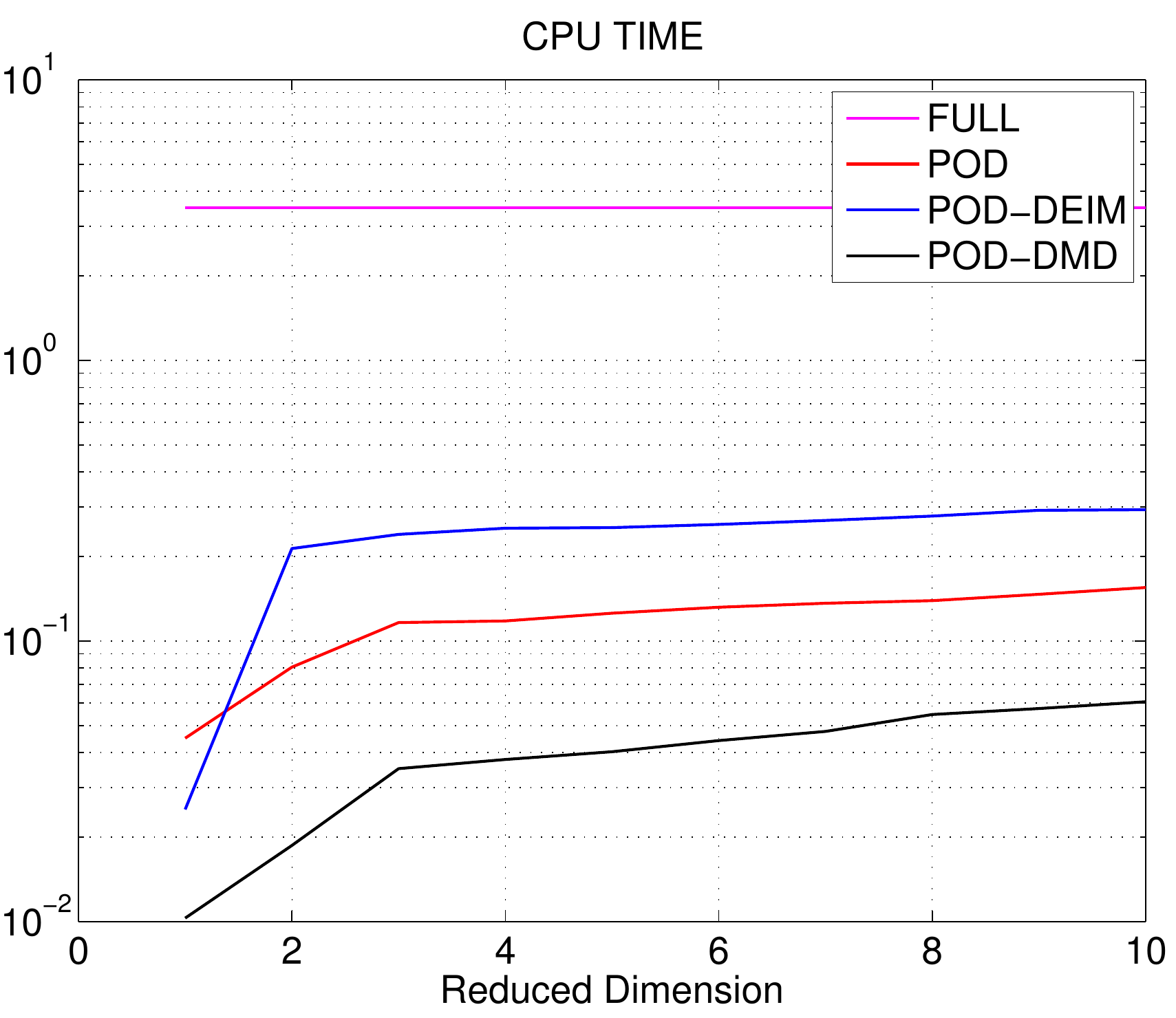}
\includegraphics[scale=0.25]{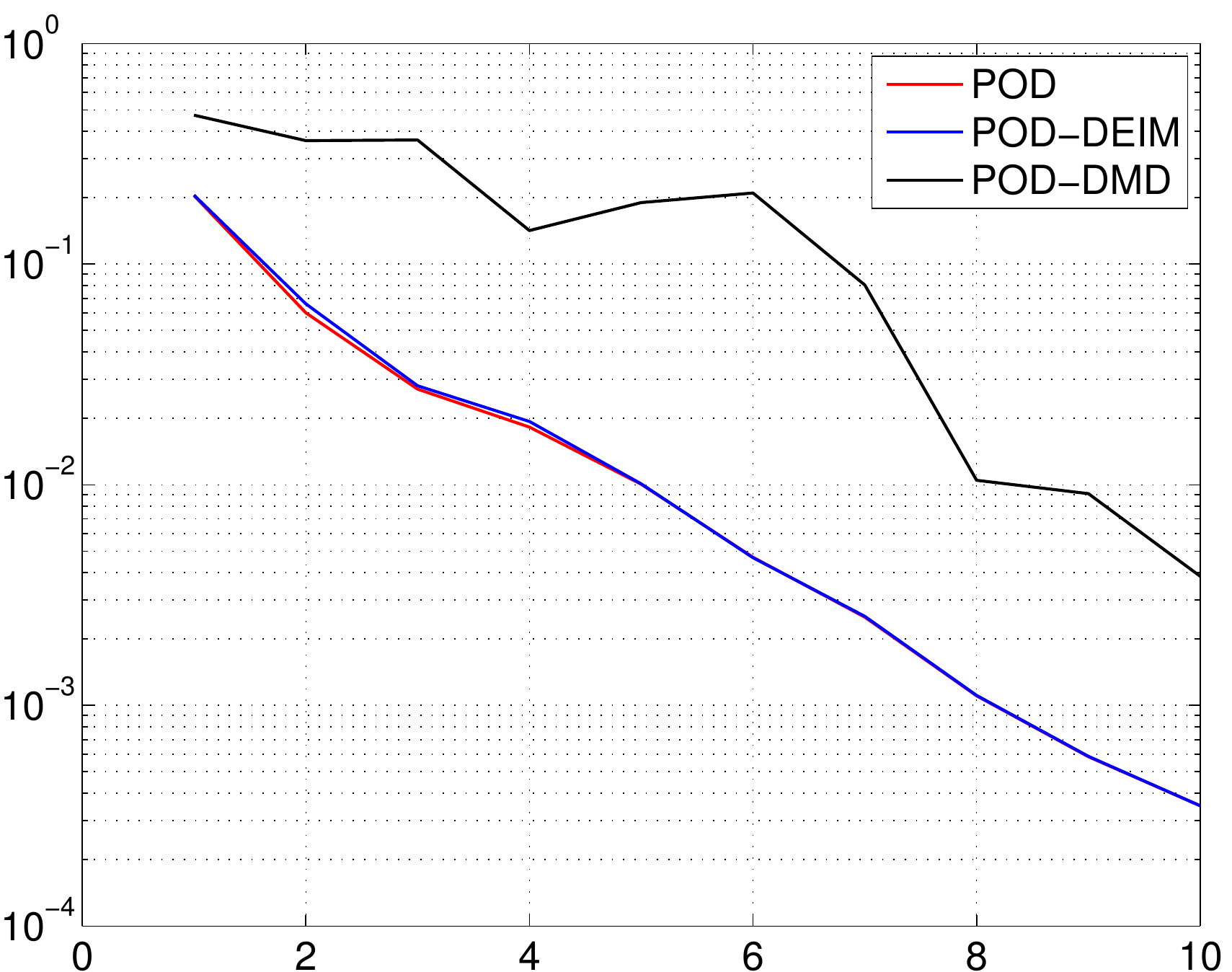}\\
\end{center}
\caption{Test 5: CPU-time (left) and Relative Error in Frobenius norm. Number of POD modes and DEIM/DMD points are the same}
\label{test3:err}
\begin{center}
\includegraphics[scale=0.22]{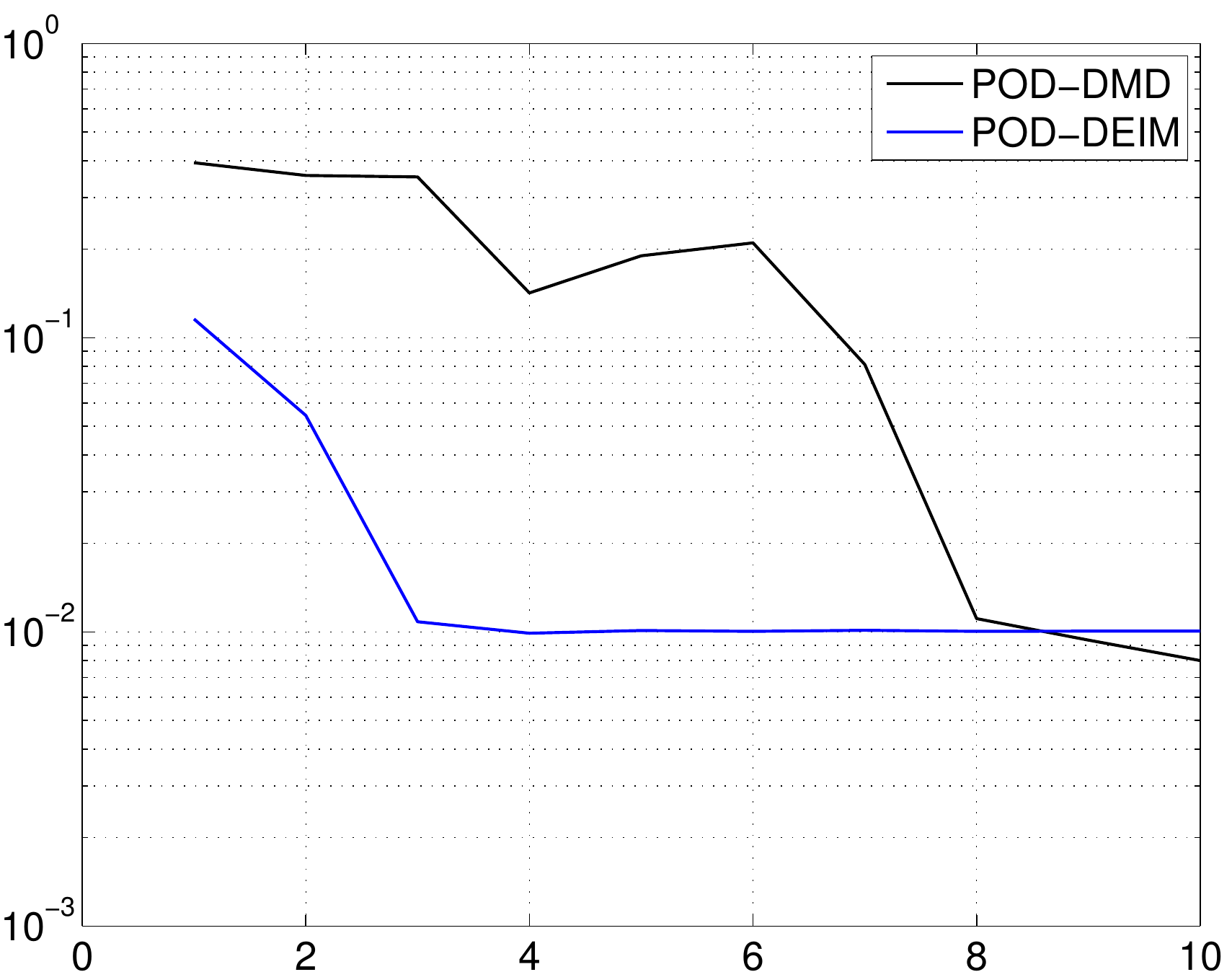}
\includegraphics[scale=0.22]{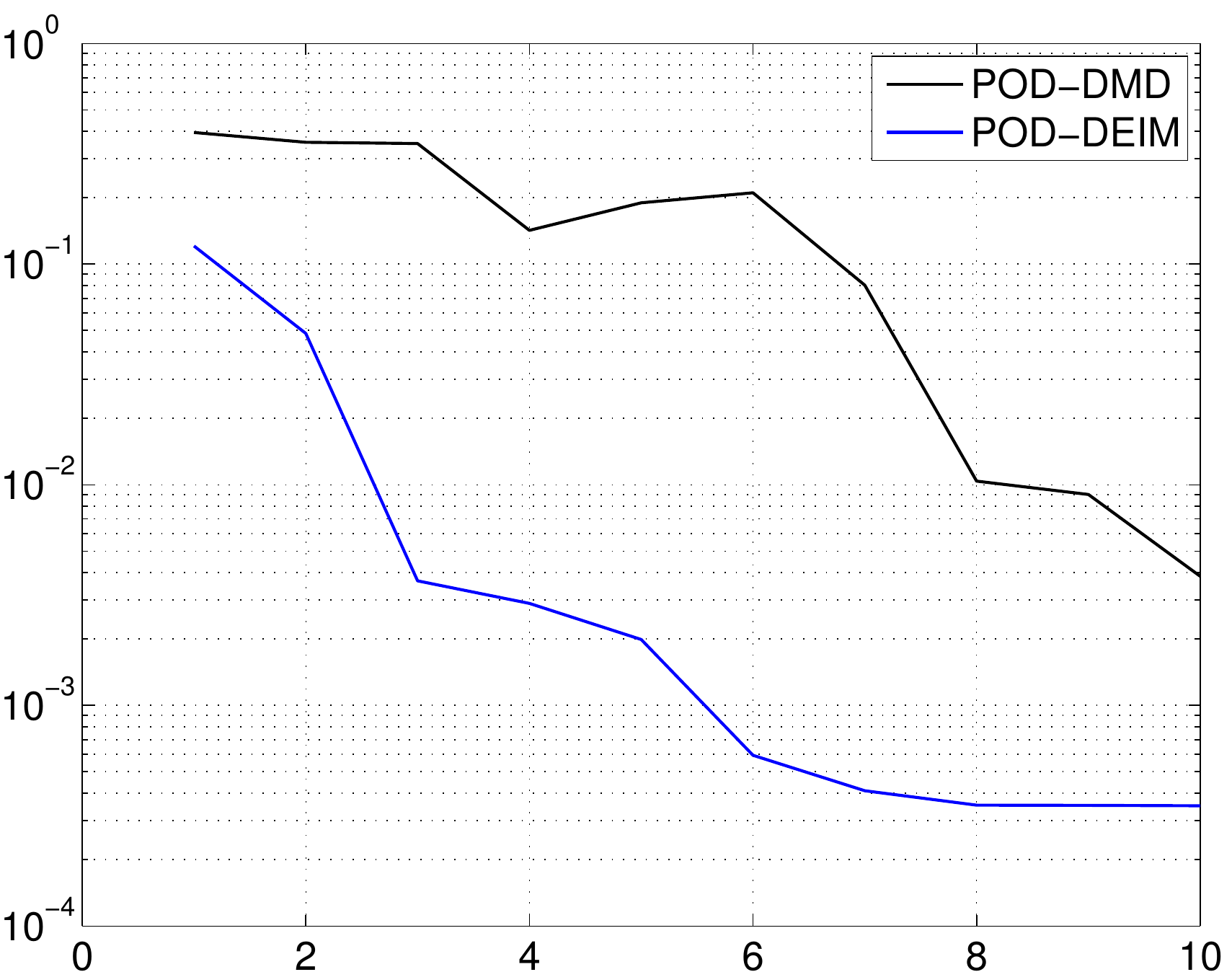}
\includegraphics[scale=0.22]{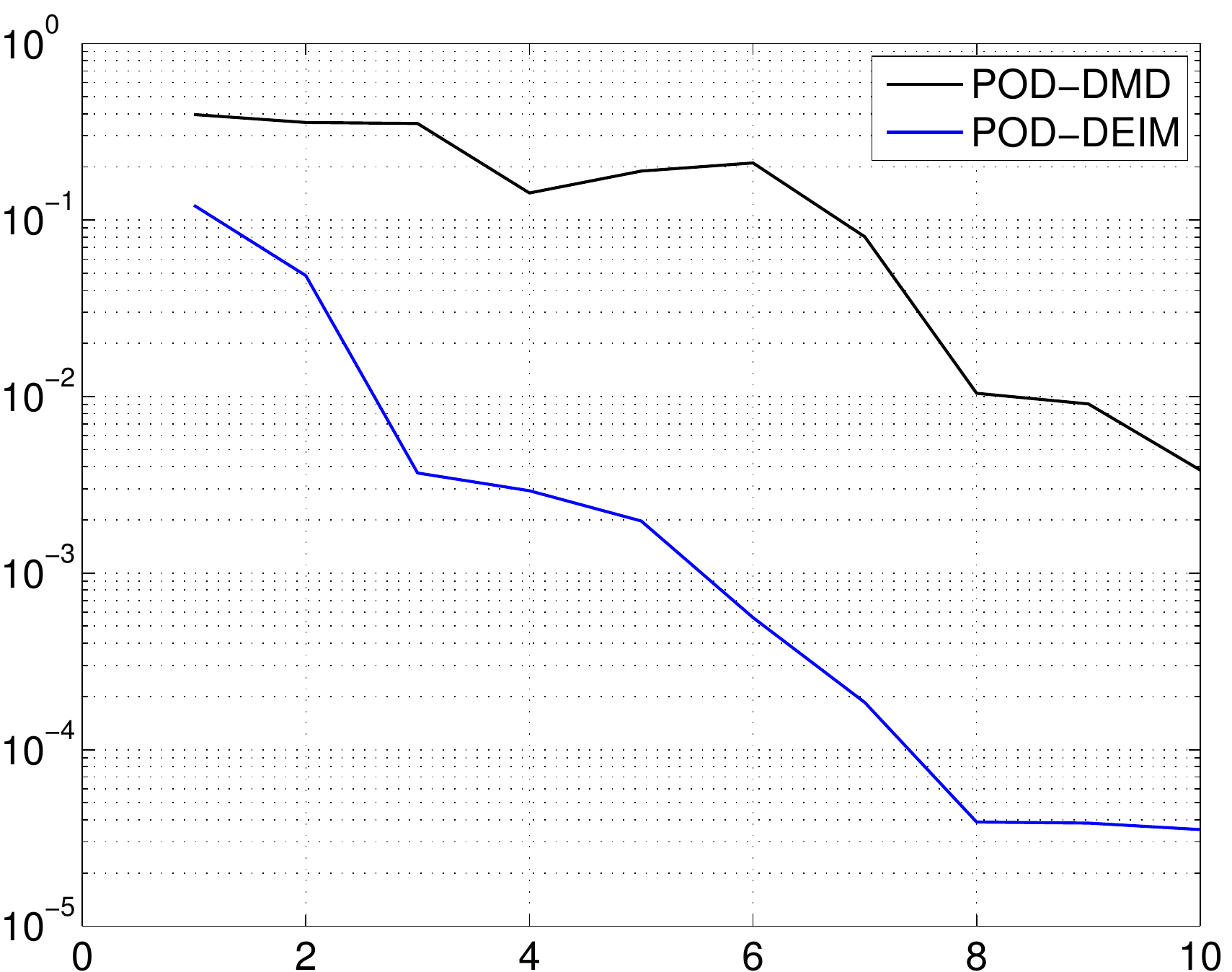}
\end{center}
\caption{Test 5: Relative Error for 5 POD basis functions (left), 10 POD basis (middle), 15 POD basis (right)}
\label{test3:err2}

\end{figure}
In this case the POD-DEIM performs always better then POD-DMD, but it is extremely expensive.

%

Approximation by means of model order reduction technique is in Figure \ref{test3:sol}, we can see it is hard to visualize distinguish any difference between the reduced solutions.
Again, we emphasize the speed up of the POD-DMD method with respect to the other methods. In this example, POD-DEIM turns out to be more expensive than POD itself as we can see in Figure \ref{test3:err}.
The POD-DMD is the least accurate.

\section{Conclusions and future work}
\label{Section6}
\setcounter{section}{6}
\setcounter{equation}{0}
\setcounter{theorem}{0}
\setcounter{algorithm}{0}
\renewcommand{\theequation}{\arabic{section}.\arabic{equation}}

In order to make model reduction methods such as POD computationally efficient, innovative methods for evaluating the nonlinear terms of the governing equations \eqref{ode} must be used. Previous successful techniques use sparse sampling to evaluate the nonlinearity. Indeed, the discrete empirical interpolation method identifies through a greedy algorithm a limited number of spatial sampling locations that can allow for reconstruction of the nonlinear terms in a low-dimensional manner. Such sparse sampling of the nonlinearity is directly related to compressive sensing strategies whereby a small number of
sensors can be used to characterize the dynamics of the high-dimensional {{nonlinear}} system.
In this paper we present a new model order reduction approach for nonlinear dynamical systems. The method couples the POD method for the projection of the system and the DMD algorithm for the approximation of the nonlinear term. In particular, DMD provides a significant reduction of the system in terms of the CPU time since the nonlinearity is never evaluated online. 
{The method is effective for {{nonlinear dynamical}} systems where POD approximations are relevant. 

The advantages of the POD-DMD method over traditional POD-DEIM have been demonstrated in numerous computational examples. Specifically, the method shows marked improvement in the computational speed for evaluating the nonlinearity, performing nearly an order of magnitude faster in comparison to POD-DEIM. However, a drawback of the method is related to the fact that the DMD modes produced are not orthogonal, thus limiting the performance of the method in terms of error convergence. Thus although the number of modes can be increased, the error plateaus, limiting how well one can approximate the original system with the POD-DMD low-dimensional system.
Future work will focus on error estimation of the proposed method in both the DMD and POD-DMD projection techniques. Special focus will be given to improving the error estimates. We also intent to use recent innovations in the DMD method, specifically around multi-resolution analysis~\cite{Kutz:mrdmd} and compression~\cite{brunton:cs}, to more effectively construct approximations to the nonlinear dynamics. 


\begin{thebibliography}{99}


          




\bibitem{mpe} P. Astrid, 
\newblock {\em Fast reduced order modeling technique for large scale LTV systems,}
\newblock  in Proc. 2004 Am. Control Conf. {\bf 1}, 762-767 (2004).



\bibitem{BMNP04}
M. Barrault, Y. Maday, N.C. Nguyen, A.T. Patera,
\newblock An empirical interpolation method: application to efficient reduced-basis discretization of partial differential equations
\newblock Comptes Rendus Mathematique, {\bf 339} (2004), 667–672.

 \bibitem{karen1} P. Benner, S. Gugercin and K. Willcox, 
 \newblock {\em A Survey of Projection-Based Model Reduction Methods for Parametric Dynamical Systems,}
\newblock SIAM Review , to appear, 2015.

\bibitem{brunton:cs}
S. L. Brunton, J. L. Proctor, and J. N. Kutz.
\newblock {\em Compressive sampling and dynamic mode decomposition.}
(to appear) J. Comp. Dyn. (2016)
  
\bibitem{Carlberg:2013}
K.~Carlberg, C.~Farhat, J.~Cortial, and D.~Amsallem.
\newblock {\em The {GNAT} method for nonlinear model reduction: Effective
  implementation and application to computational fluid dynamics and turbulent
  flows.}
\newblock Journal of Computational Physics, 242:623--647, 2013.




\bibitem{CS10}
S. Chatarantabut, D. Sorensen.
\newblock Nonlinear Model Reduction via Discrete Empirical Interpolation.
\newblock {\em SIAM J. Sci. Comput}, {\bf 32} (2010), 2737-2764.


\bibitem{DG15}
Z. Drmac, S. Gugercin.
\newblock A new selection operator for the discrete empirical interpolation method - improved a priori error bound and extensions
\newblock Preprint


\bibitem{gap1}  R. Everson and L. Sirovich, ``Karhunen-Lo\'eve procedure for gappy data,"  J. Opt. Soc. Am. A {\bf 12}, 1657-1664 (1995).

\bibitem{gavish}
M.~Gavish and D.~L. Donoho.
\newblock The optimal hard threshold for singular values is $4/\sqrt{3}$.
\newblock {\em ArXiv e-prints}, 2014.


   
\bibitem{HLBR_turb}
P.~J. Holmes, J.~L. Lumley, G.~Berkooz, and C.~W. Rowley.
\newblock {\em Turbulence, coherent structures, dynamical systems and
  symmetry}.
\newblock Cambridge Monographs in Mechanics. Cambridge University Press,
  Cambridge, England, 2nd edition, 2012.
  
  \bibitem{hotellingJEdPsy33}
H.~Hotelling.
\newblock Analysis of a complex of statistical variables into principal
  components.
\newblock {\em J. Educ.\ Psychol.}, 24:417--441 (1933).



\bibitem{koopman}
B.~O. Koopman.
\newblock \textsl{Hamiltonian Systems and Transformation in Hilbert Space}.
\newblock {\em PNAS}, 17:315--318, 1931.


 
 \bibitem{KV01} 
K. Kunisch, S. Volkwein.
\newblock Galerkin proper orthogonal decomposition methods for parabolic problems.
\newblock {\em Numer. Math.} 90 (2001), 117-148.


\bibitem{KV02}
K. Kunisch, S. Volkwein.
\newblock{Galerkin proper orthogonal decomposition methods for a general equation in fluid dynamics.}
\newblock{SIAM, J. Numer. Anal. {\bf 40} (2002), 492-515}.


\bibitem{Kutz:mrdmd}
J. N. Kutz, X. Fu and S. L. Brunton.
\newblock {\em Multi-Resolution Dynamic Mode Decomposition}
(accepted) SIAM J. App. Dyn. Systems (2016).

\bibitem{Kutz:2013}
J.~N. Kutz.
\newblock {\em Data-Driven Modeling \& Scientific Computation: Methods for
  Complex Systems \& Big Data}.
\newblock Oxford University Press, 2013.


\bibitem{lorenzMITTR56}
E.~N. Lorenz.
\newblock Empirical orthogonal functions and statistical weather prediction.
\newblock Technical report, Massachusetts Institute of Technology, December
  1956.


\bibitem{Lumley:1970}
J.~L. Lumley.
\newblock {\em Stochastic Tools in Turbulence}.
\newblock Academic Press, 1970.


\bibitem{Mezic2004}
Igor Mezi{\'c} and Andrzej Banaszuk.
\newblock Comparison of systems with complex behavior.
\newblock {\em Physica D: Nonlinear Phenomena}, 197(1-2):101 -- 133, 2004.


\bibitem{Mezic2005}
Igor Mezi{\'c}.
\newblock Spectral properties of dynamical systems, model reduction and
  decompositions.
\newblock {\em Nonlinear Dynamics}, 41(1-3):309--325, 2005.

\bibitem{mezic2}
I.~Mezi\'c.
\newblock \textsl{Analysis of Fluid Flows via Spectral Properties of the
  Koopman Operator}.
\newblock {\em Annual Review of Fluid Mechanics}, 45:357--378, 2013.


\bibitem{patera} 
N.C. Nguyen, A. T. Patera, J. Peraire, 
\newblock {\em A ``best points'' interpolation method for efficient approximation of parametrized functions.}
\newblock Int. J. Num. Methods Eng. {\bf 73}, 521--543 (2008).
    


\bibitem{Pearson:1901}
K.~Pearson.
\newblock On lines and planes of closest fit to systems of points in space.
\newblock {\em Philosophical Magazine}, 2(7--12):559--572, 1901.
    



 \bibitem{rom_book}  
\newblock A. Quarteroni and G. Rozza Eds. 
\newblock  {\em Reduced Order Methods for Modeling and Computational Reduction}, (Springer, 2014)


\bibitem{DMD4}
C.~Rowley, I.~Mezi\'{c}, S.~Bagheri, P.~Schlatter, and D.~Henningson.
\newblock \textsl{Spectral analysis of nonlinear flows}.
\newblock {\em Journal of Fluid Mechanics}, 641:115--127, 2009.



\bibitem{DMD0}
P.~J. Schmid and J.~Sesterhenn.
\newblock Dynamic mode decomposition of numerical and experimental data.
\newblock In {\em 61st Annual Meeting of the APS Division of Fluid Dynamics}.
  American Physical Society, November 2008.

\bibitem{DMD1}
P.~Schmid.
\newblock \textsl{Dynamic mode decomposition of numerical and experimental
  data}.
\newblock {\em Journal of Fluid Mechanics}, 656:5--28, 2010.



\bibitem{Sir87}
 L. Sirovich. 
\newblock {\em Turbulence and the dynamics of coherent structures. Parts I-II,}
\newblock Quarterly of Applied Mathematics, {\bf XVL} (1987), 561-590.


\bibitem{DMD5}
J.~Tu, C.~Rowley, D.~Luchtenberg, S.~Brunton, and J.~N. Kutz.
\newblock \textsl{On Dynamic Mode Decomposition: Theory and Applications}.
\newblock {\em Journal of Computational Dynamics}, 1:391--421, 2014.



\bibitem{Vol11}
S. Volkwein.
\newblock {Model Reduction using Proper Orthogonal Decomposition.}.
\newblock Lecure Notes, University of Konstanz, 2013



\bibitem{gap2}  
K. Willcox, 
\newblock {\em Unsteady flow sensing and estimation via the gappy proper orthogonal decomposition,} 
\newblock  Computers and Fluids {\bf 35}: 208-226 (2006).


\bibitem{karni} B. Yildirim, C. Chryssostomidis and G.E. Karniadakis, ``Efficient sensor placement for ocean measurements using low-dimensional concepts," Ocean Modeling, {\bf 273}(3-4), 160-173, (2009).







\end{thebibliography}
\end{document}